\documentclass[]{amsart}

\usepackage{amsmath}

   \usepackage[colorlinks = true,
            linkcolor = blue,
            urlcolor  = blue,
            citecolor = blue,
            anchorcolor = blue]{hyperref}

\theoremstyle{plain}
\newtheorem{thm}{Theorem}[subsection]
\newtheorem{prop}[thm]{Proposition}
\newtheorem{lemma}[thm]{Lemma}
\newtheorem{cor}[thm]{Corollary}

\theoremstyle{definition}
\newtheorem{defn}[thm]{Definition}
\newtheorem{ex}[thm]{Example}

\theoremstyle{remark}
\newtheorem{rmk}[thm]{Remark}

\newcommand{\vol}{\ensuremath{\operatorname{vol}}}

\newcommand{\real}{\ensuremath{\operatorname{Re}}}
\newcommand{\imag}{\ensuremath{\operatorname{Im}}}

\newcommand{\G}{\ensuremath{\operatorname{G_2}}}
\newcommand{\Gs}{\ensuremath{\operatorname{G_2}}{-structure}}
\newcommand{\SP}{\ensuremath{\operatorname{Spin (7)}}}
\newcommand{\SPs}{\ensuremath{\operatorname{Spin (7)}}{-structure}}

\newcommand{\ph}{\ensuremath{\varphi}}

\newcommand{\Ph}{\ensuremath{\Phi}}
\newcommand{\wzeo}{\ensuremath{\wedge^0_1}}
\newcommand{\wons}{\ensuremath{\wedge^1_7}}
\newcommand{\wtws}{\ensuremath{\wedge^2_7}}
\newcommand{\wtwf}{\ensuremath{\wedge^2_{14}}}
\newcommand{\wtho}{\ensuremath{\wedge^3_1}}
\newcommand{\wths}{\ensuremath{\wedge^3_7}}
\newcommand{\wtht}{\ensuremath{\wedge^3_{27}}}
\newcommand{\wfoo}{\ensuremath{\wedge^4_1}}
\newcommand{\wfos}{\ensuremath{\wedge^4_7}}
\newcommand{\wfot}{\ensuremath{\wedge^4_{27}}}
\newcommand{\wfis}{\ensuremath{\wedge^5_7}}
\newcommand{\wfif}{\ensuremath{\wedge^5_{14}}}
\newcommand{\wsis}{\ensuremath{\wedge^6_7}}
\newcommand{\wseo}{\ensuremath{\wedge^7_1}}
\newcommand{\wone}{\ensuremath{\wedge^1_8}}
\newcommand{\wtwt}{\ensuremath{\wedge^2_{21}}}
\newcommand{\wthe}{\ensuremath{\wedge^3_8}}
\newcommand{\wthfo}{\ensuremath{\wedge^3_{48}}}
\newcommand{\wfoth}{\ensuremath{\wedge^4_{35}}}

\newcommand{\wfie}{\ensuremath{\wedge^5_8}}
\newcommand{\wfifo}{\ensuremath{\wedge^5_{48}}}
\newcommand{\wsit}{\ensuremath{\wedge^6_{21}}}
\newcommand{\wsee}{\ensuremath{\wedge^7_8}}
\newcommand{\weio}{\ensuremath{\wedge^8_1}}
\newcommand{\Rse}{\ensuremath{\mathbb R^7}}
\newcommand{\Rei}{\ensuremath{\mathbb R^8}}
\newcommand{\wthpos}{\ensuremath{\wedge^3_{\text{pos}}}}
\newcommand{\wfopos}{\ensuremath{\wedge^4_{\text{pos}}}}

\newcommand{\st}{\ensuremath{\ast}}
\newcommand{\hk}{\ensuremath{\lrcorner}}
\newcommand{\stph}{\ensuremath{\ast \varphi}}
\newcommand{\lieg}{\ensuremath{\operatorname{\mathfrak {g}_2}}}
\newcommand{\liesp}{\ensuremath{\operatorname{\mathfrak
{spin}(7)}}}

\newcommand{\ws}{\ensuremath{w^{\flat}}}
\newcommand{\us}{\ensuremath{u^{\flat}}}
\newcommand{\vs}{\ensuremath{v^{\flat}}}
\newcommand{\as}{\ensuremath{a^{\flat}}}
\newcommand{\bs}{\ensuremath{b^{\flat}}}
\newcommand{\cs}{\ensuremath{c^{\flat}}}
\newcommand{\ds}{\ensuremath{d^{\flat}}}
\newcommand{\oph}{\ensuremath{\varphi_{\!\text{o}}}}
\newcommand{\ost}{\ensuremath{\ast_{\!\text{o}}}}
\newcommand{\ostph}{\ensuremath{\ost \oph}}
\newcommand{\oPh}{\ensuremath{\Phi_{\!\text{o}}}}
\newcommand{\ovol}{\ensuremath{\vol_{\text{o}}}}
\newcommand{\omu}{\ensuremath{\mu_{\text{o}}}}
\newcommand{\ozeta}{\ensuremath{\zeta_{\text{o}}}}
\newcommand{\otheta}{\ensuremath{\theta_{\!\text{o}}}}
\newcommand{\oh}{\ensuremath{h_{\!\text{o}}}}
\newcommand{\og}{\ensuremath{g_{\text{o}}}}
\newcommand{\lt}[1]{\ensuremath{{#1}_{\,{}_{\widetilde{}}}  }}

\newcommand{\lttt}[1]{\ensuremath{{#1}_{\:{}_{\widetilde{}}}  }}
\newcommand{\ltn}[1]{\ensuremath{{#1}_{{}_{\widetilde{}}} }}
\newcommand{\nph}{\ensuremath{\tilde \varphi}}
\newcommand{\nPh}{\ensuremath{\tilde \Phi}}
\newcommand{\nst}{\ensuremath{\tilde \ast}}
\newcommand{\nstph}{\ensuremath{\tilde \st \tilde \varphi}}
\newcommand{\nvol}{\ensuremath{\lttt{\vol}}}
\newcommand{\nmu}{\ensuremath{\tilde \mu}}
\newcommand{\nzeta}{\ensuremath{\tilde \zeta}}
\newcommand{\ntheta}{\ensuremath{\tilde \theta}}
\newcommand{\nh}{\ensuremath{\tilde h}}
\newcommand{\ointeg}{\ensuremath{d \ostph + \frac{1}{3} \otheta
\wedge \ostph }}
\newcommand{\ointdual}{\ensuremath{d \oph + \frac{1}{4} \otheta
\wedge \oph }}
\newcommand{\ninteg}{\ensuremath{d \nstph + \frac{1}{3} \ntheta
\wedge \nstph }}
\newcommand{\nintdual}{\ensuremath{d \nph + \frac{1}{4} \ntheta
\wedge \nph }}
\newcommand{\oSPclass}{\ensuremath{d \oPh + \frac{1}{7} \otheta
\wedge \oPh }}
\newcommand{\nSPclass}{\ensuremath{d \nPh + \frac{1}{7} \ntheta
\wedge \nPh }}
\newcommand{\opw}{\ensuremath{( 1 + {|w|}_{\text{o}}^2 ) }}
\newcommand{\ovw}{\ensuremath{{|v \wedge w|}}}
\newcommand{\opvw}{\ensuremath{\left( 1 + \frac{4}{7}{|v \wedge
w|}_{\text{o}}^2 \right) }}
\newcommand{\ddbar}{\ensuremath{\partial \bar{\partial}}}

\numberwithin{equation}{section}
\numberwithin{table}{section}
\numberwithin{figure}{section}

\begin{document}

\title[\G\ and \SP\ Structures on Manifolds]{Deformations of
\G\ and \SP\ \\ Structures on Manifolds}
\author{Spiro Karigiannis}
\address{Department of Mathematics and Statistics\\ McMaster
University}
\email{spiro@math.mcmaster.ca}
\subjclass{53, 58}
\date{\today}

\begin{abstract}
We consider some deformations of \Gs s on $7$-manifolds. We
discover a canonical way to deform a \Gs\ by a vector field in
which the associated metric gets ``twisted'' in some way by the
vector cross product. We present a system of partial differential
equations for an unknown vector field $w$ whose solution would
yield a manifold with holonomy \G. Similarly we consider analogous
constructions for \SP -structures on $8$-manifolds. Some of the
results carry over directly, while others do not because of the
increased complexity of the \SP\ case.

\smallskip

\noindent \textcolor{red}{Changes made after publication: There were sign errors in equations (2.26) and (2.34) and in the proof of Lemma 2.4.3 that have now been corrected in this arXiv version. Also Remark 3.3.7 was incorrect and has been removed.}
\end{abstract}

\maketitle

\tableofcontents

\section{Introduction} \label{introsec}

\subsection{Cross Product Structures} \label{crossintrosec}

Special structures on manifolds in differential geometry are
usually described by globally defined smooth sections of some
tensor bundle, satisfying some pointwise algebraic conditions, and
with potentially stronger global requirements. One such example of
an additional structure that can be imposed on a smooth Riemannian
manifold $M$ of dimension $n$ is that of an
$r$-fold {\em cross product}. This is an {\em alternating}
$r$-linear smooth map
\begin{equation*}
B : \underbrace{TM \times \ldots \times TM}_{r \text { copies}}
\to TM
\end{equation*}
that is compatible with the metric in the sense that
\begin{eqnarray*}
|B(e_1, \ldots, e_r)|^2 & = &|e_1 \wedge \ldots \wedge e_r|^2 \\
\langle B(e_1, \ldots, e_r) , e_j \rangle & = & 0 \qquad 1 \leq j
\leq r
\end{eqnarray*}
where $\langle \cdot, \cdot \rangle$ is the Riemannian metric.
Such a cross product also gives rise to an $(r+1)$-form $\alpha$
given by
\begin{equation*}
\alpha(e_1,\ldots, e_r, e_{r+1}) = \langle B(e_1,\ldots, e_r),
e_{r+1} \rangle
\end{equation*}

Cross products on real vector spaces were classified by Brown and
Gray in~\cite{BG}. Global vector cross products on manifolds
were first studied by Gray in~\cite{Gr5}. They fall into four
categories:
\begin{enumerate}
\item When $r = n-1$ and $\alpha$ is the volume form of the
manifold. Under the metric identification of vector fields and one
forms, this cross product corresponds to the {\em Hodge star}
operator on $(n-1)$-forms. This is not an extra structure beyond
that given by the metric.
\item When $n = 2 m$ and $r = 1$, we can have a one-fold cross
product $J: TM \to TM$. Such a map satisfies $J^2 = - I$ and is an
{\em almost complex structure}. The associated $2$-form is the
{\em K\"ahler form} $\omega$.
\item The first of two {\em exceptional cases} is a $2$-fold cross
product on a $7$-manifold. Such a structure is called a \Gs, and
the associated $3$-form $\ph$ is called a \G-form.
\item The second exceptional case is a $3$-fold cross product on
an $8$-manifold. This is called a \SPs, and the associated
$4$-form $\Ph$ is called a \SP-form.
\end{enumerate}
In cases 2--4 the existence of these structures is a
topological condition on $M$ given in terms of characteristic
classes (see~\cite{Gr5, J4, Sa}). One can also study the
restricted sub-class of such manifolds where the associated
differential form
$\alpha$ is {\em parallel} with respect to the Levi-Civita
connection $\nabla$. In case (1), the volume form is always
parallel. For the almost complex structures
$J$ of case (2), $\nabla J = 0$ if and only if the manifold is
K\"ahler, which is equivalent to $d \omega = 0$ {\em and} the
almost complex structure is {\em integrable}. In this case, the
Riemannian holonomy of the manifold is a subgroup of $U(m)$. For
cases (3) and (4), the condition that the differential form be
parallel is a non-linear differential equation. Manifolds with
parallel \Gs s have holonomy a subgroup of \G\ and manifolds with
parallel \SPs s have holonomy a subgroup of \SP, hence their names.
One can also show (see~\cite{Bo}) that such manifolds are all {\em
Ricci-flat}.

There is a sub-class of the K\"ahler manifolds which are
Ricci-flat. Such manifolds possess a global non-vanishing
holomorphic {\em volume form} $\Omega$ in addition to the K\"ahler
form $\omega$, and these two forms satisfy some relation. These
manifolds are called {\em Calabi-Yau} manifolds as their existence
was demonstrated by Yau's proof of the Calabi conjecture~\cite{Y}:
\begin{thm}[Calabi-Yau, 1978]
Let $M$ be a compact complex manifold with vanishing first Chern
class $c_1 = 0$. Then if $\omega$ is a K\"ahler form on $M$, there
exists a unique Ricci-flat Riemannian metric $g$ on $M$ whose
associated K\"ahler form is in the same cohomology class as
$\omega$.
\end{thm}

This theorem characterises those manifolds admitting Calabi-Yau
metrics in terms of certain topological information. The
equivalence is demonstrated by writing the Ricci-flat condition as
a partial differential equation and proving existence and
uniqueness of solutions. Calabi-Yau manifolds have holonomy a
subgroup of $SU(m)$ and are characterized by {\em two} parallel
forms, $\omega$ and $\Omega$. In fact, they posses two parallel
cross products: a $1$-fold cross product $J$, and a complex
analogue of case (1) above, where $\Omega$ plays the role of the
volume form and the $(m-1)$-fold cross product is a {\em complex}
Hodge star.

Calabi-Yau manifolds (at least in complex dimension $3$) have long
been of interest in string theory. More recently, manifolds with
holonomy \G\ and \SP\ have also been studied. (See, for
example,~\cite{Br2, BS, J1, J2, J3, Hi1, Hi2, GYZ}). It would be
useful to have an analogue of the Calabi-Yau theorem, or something
similar, in the \G\ and \SP\ cases. There is a significant
difference, however, which makes \G\ and \SP\ manifolds much more
difficult to study.
An almost complex structure $J$ does not by itself determine
a metric. If we also have a Riemannian metric, then together the
compatibility requirement yields the K\"ahler form $\omega(u,v) =
g(Ju,v)$. In contrast, a $2$-fold or $3$-fold cross product
structure {\em does} determine the metric uniquely, and thus
also determines the associated $3$-form $\ph$ or $4$-form $\Ph$.
Because the metric and complex structure are ``uncoupled'' in the
Calabi-Yau case, we can start with a {\em fixed} integrable
complex structure $J$, and then look for different metrics (which
correspond to different K\"ahler forms for the same $J$) which are
Ricci-flat and make $J$ parallel. As $J$ is already integrable, it
is parallel precisely when $\omega$ is closed, so we can simply
look at different metrics which all correspond to closed K\"ahler
forms, and from that set look for a Ricci-flat metric. Hence we
can restrict ourselves to starting with a K\"ahler manifold, and
looking at other Kahler metrics which could be Ricci-flat. The
Calabi-Yau theorem then says that there exists precisely one such
metric in each cohomology class which contains at least one
K\"ahler metric.

In the \G\ and \SP\ cases, however, we cannot fix a cross product
structure and then vary the metric to make it parallel. For a
given cross product, the metric is determined. In the Calabi-Yau
case, we can start with $U(m)$ holonomy and describe the
conditions for being able to obtain $SU(m)$ holonomy. For \G\ and
\SP, there is no intermediate starting class. A crucial ingredient
in the proof of the Calabi-Yau theorem is the \ddbar\ lemma, which
allows us to write the difference of any two K\"ahler forms in
terms of an unknown {\em function} $f$. Therefore as a first step
towards an analogous result in the \G\ and \SP\ cases, we would
like to determine the simplest data required to describe the
relations between any two \G\ or \SP\ forms. 
\subsection{Overview of New Results} \label{resultssec}

If we start with only a \Gs, not necessarily parallel, this gives
us a $3$-form which satisfies some ``positive-definiteness''
property, since it determines a Riemannian metric. In~\cite{FG},
Fern\'andez and Gray classified such manifolds by looking at the
decomposition of $\nabla \ph$ into \G-irreducible components.
There are $16$ such classes, with various inclusion relations
between them. There is a similar decomposition in~\cite{GH} of
almost complex manifolds into subclasses. Some of these classes
are: integrable (complex), symplectic, almost K\"ahler, and nearly
K\"ahler. Thus these $16$ subclasses of manifolds with a \Gs\
are analogues of these ``weaker than K\"ahler'' conditions.
Similar studies by Fern\'anadez in~\cite{F1} of the \SP\ case yield
$4$ subclasses of manifolds with a \SPs.

As a first step in trying to determine an analogue for the Calabi
conjecture in the \G\ case, we can study these various
weaker subclasses and their deformations. If we start in one
class, and change the $3$-form $\ph$ in some way (which changes
the metric too) we would like to know under what conditions this
subclass is preserved, or more generally what subclass the new
\Gs\ now belongs to. The space of $3$-forms on a manifold with a
\Gs\ decomposes into a direct sum of irreducible
\G-representations:
\begin{equation*}
\wedge^3 = \wtho \oplus \wths \oplus \wtht
\end{equation*}
where $\wedge^3_k$ is a $k$-dimensional vector space at each point
on $M$. This decomposition depends on our initial $3$-form $\oph$,
however. This again is in stark contrast to the decomposition on a
complex manifold into forms of type $(p,q)$, which depends only on
the complex structure and does not change as we vary the K\"ahler
(or metric) structure. We can consider a deformation
$\nph = \ph_0 + \eta$ of the \Gs, for $\eta \in
\wedge^3_k$ and determine conditions on $\oph$ and $\eta$ which
preserve the subclass or change it in an interesting way.

If $\eta \in \wtho$, this corresponds to a conformal scaling
of the metric, and one can explicitly describe which of the
$16$ classes are conformally invariant. (These results were
already known to Fern\'andez and Gray but here they are reproduced
in a different way.) A new result in this case is the
following:
\begin{thm}
Let $\theta_o = \ost ( \ost d \oph \wedge \oph )$ be the canonical
$1$-form arising from a \Gs\ $\ph_o$. Then if $\nph = f^3 \oph$
for some non-vanishing function $f$, the new canonical $1$-form
$\tilde \theta$ differs from the old $\theta_o$ by an exact form:
\begin{equation*}
\tilde \theta = -12 d (\log (f)) + \theta_o
\end{equation*}
Thus in the classes where $\theta$ is closed, (there are some and
they are conformally invariant classes), we get a well-defined
cohomology class in $H^1(M)$, invariant under conformal changes
of metric. A similar result also holds in the \SP\ case.
\end{thm}

If, however, we deform $\oph$ by an element $\eta \in \wths$, 
then $\eta = w \hk \ostph$ for some vector field $w$, and 
in Section~\ref{vfdeformsec} we prove the following:
\begin{thm}
Under such a deformation $\nph$ is again a \Gs\, and the new
metric on vector fields $v_1$ and $v_2$ is given by
\begin{equation*}
\lt{\langle v_1, v_2 \rangle} = \frac{1}{{\opw}^{\frac{2}{3}}}
\left( {\langle v_1, v_2 \rangle}_{\! \text{o}} + {\langle w
\times v_1, w \times v_2 \rangle}_{\! \text{o}} \right)
\end{equation*}
where $\times$ is the vector cross product associated to the
original \Gs\ $\oph$.
\end{thm}
From this one can write down non-trivial differential equations on
the vector field $w$ for certain subclasses to be preserved. It
would be interesting to solve some of these equations for the
unknown vector field $w$. This would mean that there were certain
distinguished vector fields on some classes of manifolds with \Gs
s. The important result here, however, is that the new $3$-form
$\nph$ is {\em always} positive-definite. That is, it always
corresponds to a \Gs. This gives information about the structure
of the open set $\wedge^3_+(M)$ of positive definite $3$-forms on
$M$.

If instead we deform $\ph$ in the $\wths$ direction infinitesmally
by the flow equation
\begin{equation*}
\frac{\partial}{\partial t} \ph_t = w \hk \st_t \ph_t
\end{equation*}
then we show in Section~\ref{infvfdeformsec} that the metric $g$
{\em does not change} and also:
\begin{thm}
The solution is given by
\begin{equation*}
\ph (t) = \ph_0 + \frac{1 - \cos(|w|t)}{{|w|}^2}(w \hk \st (w \hk
\st \ph_0 ) ) + \frac{\sin(|w|t)}{|w|} (w \hk \st \ph_0 )
\end{equation*}
Hence the solution exists for all time and is a {\em closed path}
in $\wedge^3(M)$. Also, the path only depends on the unit vector
field $\pm \frac{w}{|w|}$, and the norm $|w|$ only affects the
speed of travel along this path.
\end{thm}
In~\cite{BS} the fact that the space of \Gs s which correspond to
the {\em same} metric as a fixed \Gs\ yields an $\mathbb R \mathbb
P^7$ bundle over $M$ is mentioned. This is the content of the
above theorem, and we provide an explicit description of these \Gs
s in terms of vector fields on $M$. In addition, in the special
cases of $M = N \times S^1$ or $M = L \times T^3$, where $N$ is a
Calabi-Yau $3$-fold and $L$ is a $K3$ surface, we show that
this closed path of \Gs s corresponds to the freedom of changing
the {\em phase} of the holomorphic volume form $\Omega \mapsto
e^{it} \Omega$ on $N$ or performing a {\em hyperK\"ahler rotation}
on $L$. Thus this theorem can be seen as a generalization of these
two situations.

The same kind of analysis can be done in the
\SP\ case. Similar but more complicated results hold in this
case and are presented in Section~\ref{SPdeformsec}. Here there are
only $4$ subclasses but the decomposition of
$\wedge^4$ into irreducible \SP-representations is more
complicated:
\begin{equation*}
\wedge^4 = \wfoo \oplus \wfos \oplus \wfot \oplus \wfoth
\end{equation*}
In this case it is the space $\wfos$ which infinitesmally gives a
closed path of \SPs s all corresponding to the same metric.
However, perhaps initially somewhat surprisingly, this time
non-infinitesmal  deformations in the $\wfos$ direction {\em do
not} yield a new \SPs. This is explained in detail in
Section~\ref{w47deformsec}. Much of the construction does
indeed carry over, however, and it may be possible to alter it
somehow to make it work.

\subsection{Notation and Conventions} \label{notationsec}

Many of the calculations that follow use various relations between
the interior product $\hk$, the exterior product $\wedge$, and the
Hodge star operator $\st$. Readers unfamiliar with this can
refer to Appendix~\ref{identitiessec}. The appendix also contains
some preliminary results about determinants that are used
repeatedly in many of the proofs that follow.

In much of the computations there are two metrics present: an old
metric $g_{\text{o}}$ and a new metric $\tilde g$. Their
associated volume forms, induced metrics on differential forms,
and Hodge star operators are also identified by a subscript
${}_{\text{o}}$ for old or a \ $\widetilde {}$\ \ for new. We also
often use the metric isomorphism between vector fields and
one-forms, and denote this isomorphism by $\ws$ for the one-form 
associated to the vector field $w$ and $\alpha^{\sharp}$ for the
vector field associated to the one-form $\alpha$. In the presence
of two metrics, this isomorphism is {\em always} only used for the
old metric $g_{\text{o}}$.

\subsection{Acknowledgements} \label{acksec}

A significant portion of this paper was the author's doctoral
thesis at Harvard University under the supervision of
S.T. Yau. For helpful discussions the author would
like to thank R.L.\ Bryant, H.B.\ Cheng, S.L.\ Kong, C.C.\ Liu, J.\
Rasmussen, S.T.\ Yau, and X.W.\ Wang. The author's research was
partially funded by the Natural Sciences and Engineering Research
Council of Canada and le Fonds pour la Formation de Chercheurs et
l'Aide \`a la Recherche.

\section{Manifolds with a \Gs} \label{G2manifoldssec}

\subsection{\Gs s} \label{G2structuressec}

Let $M$ be an oriented $7$-manifold with a global $2$-fold cross
product structure. Such a structure will henceforth be called a
\Gs. Its existence is a topological condition, given by the
vanishing of the second Stiefel-Whitney class $w_2 = 0$.
(See~\cite{Gr5, J4, Sa} for details.) This cross product
$\times$ gives rise to an associated {\em Riemannian metric} $g$
and an alternating $3$-form $\ph$ which are related by:
\begin{equation} \label{Gcompatibleeq}
\ph(u,v,w) = g(u \times v, w).
\end{equation}
This should be compared to the relation between a K\"ahler metric
$\omega$ and a compatible almost complex structure $J$:
\begin{equation*}
\omega(u,v) = g(Ju, v).
\end{equation*}
Note that in the K\"ahler case, the metric and the almost complex
structure can be prescribed independently. This is {\em not} true
in the case of manifolds with a \Gs, and this leads to
some complications (and the inherit non-linearity of the problem).
For a \Gs\ $\ph$, near a point $p \in M$ we can choose local
coordinates $x^1, \ldots, x^7$ so that {\em at the point $p$}, we
have:
\begin{equation} \label{phicoordinateseq}
\ph_p = dx^{123} - dx^{167} - dx^{527} - dx^{563} + dx^{415} +
dx^{426} + dx^{437}
\end{equation}
where $dx^{ijk} = dx^i \wedge dx^j \wedge dx^k$. In these
coordinates the metric at $p$ is the standard Euclidean metric
\begin{equation*}
g_p = \sum_{k = 1}^7 dx^k \otimes dx^k
\end{equation*}
and the Hodge star dual $\stph$ of $\ph$ is
\begin{equation} \label{stphicoordinateseq}
{( \stph )}_p = dx^{4567} - dx^{4523} - dx^{4163} - dx^{4127} +
dx^{2637} + dx^{1537} + dx^{1526}
\end{equation}
\begin{rmk} \label{Gsignsrmk}
Different conventions exist in the literature
for~\eqref{phicoordinateseq} and~\eqref{stphicoordinateseq},
which may or may not differ from our choice by renumbering of
coordinates and/or a change of orientation.
\end{rmk}
The $3$-forms on $M$ that arise from a \Gs\ are called {\em
positive} $3$-forms or {\em non-degenerate}. We will denote this
set by $\wthpos$. The subgroup of $\operatorname{SO(7)}$ that
preserves $\ph_p$ is the exceptional Lie group $\G$. This can be
found for example in~\cite{Br2, H}. Hence at each point $p$, the
set of \Gs s at $p$ is isomorphic to $\operatorname{GL(7, \mathbb
R)} / \G$, which is $49 - 14 = 35$ dimensional. Since $\wedge^3
(\Rse)$ is also $35$ dimensional, the set $\wthpos (p)$ of
positive $3$-forms at $p$ is an {\em open subset} of $\wedge^3_p$.
We will determine some new information about the structure of
$\wthpos$ in Section~\ref{vfdeformsec}.

\begin{rmk}
Note that in the \SP\ case the situation is very different. The set
of $4$-forms on an $8$-manifold $M$ that determine a \SPs\ is {\em
not} an open subset of $\wedge^4 (M)$. This is discussed in
Section~\ref{SPstructuressec}.
\end{rmk}

\subsection{Decompostion of $\bigwedge^{*}(M)$ into irreducible \G
-representations} \label{G2reptheorysec}
All of the facts collected in this section are well known and
more details can be found in~\cite{FG,J4,Sa}.

The group $\G$ acts on $\Rse$, and hence acts on the spaces
$\wedge^\st$ of differential forms on $M$. One can decompose each
space $\wedge^k$ into irreducible $\G$-representations. The results
of this decomposition are presented below (see~\cite{FG, Sa}). The
notation $\wedge^k_l$ refers to an $l$-dimensional irreducible
$\G$-representation which is a subspace of $\wedge^k$. Also,
``$\vol$'' will denote the volume form of $M$ (determined by the
metric $g$), and $w$ is a vector field on $M$.

\begin{eqnarray*}
\wzeo = \{f \in C^{\infty}(M) \} & \qquad \qquad & \wons = \{\alpha
\in \Gamma(\wedge^1 (M) \} \\ \wedge^2 = \wtws \oplus \wtwf &
\qquad \qquad & \wedge^3 = \wtho \oplus \wths \oplus \wtht \\
\wedge^4 = \wfoo \oplus \wfos \oplus \wfot & \qquad \qquad &
\wedge^5 = \wfis \oplus \wfif \\ \wsis = \{ w \hk \vol \} & \qquad
\qquad & \wseo = \{f \vol; f \in C^{\infty}(M) \}
\end{eqnarray*}

Since $\G \subset \operatorname{SO(7)}$, the decomposition respects
the {\em Hodge star} $\st$ operator, and $\st \wedge^k_l =
\wedge^{7-k}_l$. Nevertheless, we will still describe the
remaining cases explicitly, in several ways, as all the
descriptions will be useful for us.

Before we describe $\wedge^k_l$ for $k = 2,3,4,5$, let us describe
some isomorphisms between these subspaces:
\begin{prop} \label{isomorphisms}
The map $\alpha \mapsto \ph \wedge \alpha$ is an isomorphism
between the following spaces:
\begin{eqnarray*}
\wzeo \cong  \wtho & \qquad \qquad & \wons \cong \wfos \\ \wtws
\cong \wfis & \qquad \qquad & \wtwf \cong \wfif \\ \wths \cong
\wsis & \qquad \qquad & \wfoo \cong \wseo
\end{eqnarray*}
The map $\alpha \mapsto \stph \wedge \alpha$ is an isomorphism
between the following spaces:
\begin{eqnarray*}
\wzeo \cong \wfoo & \qquad \qquad & \wons \cong \wfis \\ \wtws
\cong \wsis & \qquad \qquad & \wtho \cong \wseo 
\end{eqnarray*}
In addition, if $\alpha$ is a $1$-form, we have the following
identities:
\begin{eqnarray} \label{isopheq} \st \left( \ph \wedge \st (\ph
\wedge \alpha) \right) & = & - 4 \alpha \\ \nonumber \stph \wedge
\st \left( \ph \wedge \alpha \right) & = & 0 \\
\label{isostpheq} \st \left( \stph \wedge \st (\stph \wedge \alpha)
\right) & = & 3 \alpha \\ \label{iso2eq} \ph \wedge \st \left(
\stph \wedge \alpha \right) & = & 2 \left( \stph \wedge \alpha
\right)
\end{eqnarray}
\end{prop}
\begin{proof}
Since the statements are pointwise, it is enough to check them in
local coordinates using~\eqref{phicoordinateseq}
and~\eqref{stphicoordinateseq}. This is tedious but
straightforward.
\end{proof}

From these identities, we can prove the following lemma:
\begin{lemma} \label{onewedgelemma}
If $\alpha$ is a $1$-form on $M$, then we have:
\begin{eqnarray} \label{onewedgepheq} {|\ph \wedge \alpha|}^2 & =
& 4 {| \alpha |}^2 \\ \label{onewedgestpheq} {|\stph \wedge
\alpha|}^2 & = & 3 {| \alpha |}^2
\end{eqnarray}
\end{lemma}
\begin{proof}
From equation~\eqref{isopheq}, we have:
\begin{eqnarray*}
\ph \wedge \st ( \ph \wedge \alpha ) & = & -4 \st \alpha \\ \alpha
\wedge \ph \wedge \st ( \ph \wedge \alpha) & = & -4 \alpha \wedge
\st \alpha \\ - {| \ph \wedge \alpha |}^2 \vol & = & - 4 {|\alpha
|}^2 \vol \\ 
\end{eqnarray*}
which proves~\eqref{onewedgepheq}. An analogous
calculation using~\eqref{isostpheq} yields~\eqref{onewedgestpheq}.
\end{proof}

We have the following relations between $\ph$, $\stph$, and an
arbitrary vector field $w$:
\begin{lemma} \label{isomorphismsv}
The following relations hold for any vector field $w$, where $\ws$
is the associated $1$-form (obtained from the metric isomorphism):
\begin{eqnarray}
\label{hkstpheq} \st \left( \ph \wedge \ws \right) & = & w \hk
\stph \\ \label{hkpheq} \st \left( \stph \wedge \ws \right) & = & w
\hk \ph \\ \label{hk4eq} \ph \wedge \left( w \hk \stph \right) & =
& - 4 \st \ws \\ \nonumber \stph \wedge \left( w \hk \stph \right)
& = & 0 \\ \label{hk3eq} \stph \wedge \left( w \hk \ph \right) & =
& 3 \st \ws \\ \label {hk2eq} \ph \wedge \left( w \hk \ph \right)
& = & 2 \st \left( w \hk \ph \right)
\end{eqnarray}
\end{lemma}
\begin{proof}
Since on a $7$-manifold $\st^2 = 1$, these results follow from
Lemma~\ref{intextidentities} and Proposition~\ref{isomorphisms}.
\end{proof}

We now explicitly describe the decomposition beginning with
$k=2,5$.
\begin{eqnarray}
\label{wtwsdesc} \wtws & = & \{ w \hk \ph; w \in \Gamma(T(M)) \} \\
\nonumber & = &
\{\beta \in \wedge^2; \st (\ph \wedge \beta) = 2 \beta \} \\
\nonumber & = & 
\{\beta \in \wedge^2; \st \left( \stph \wedge \left( \st (\stph
\wedge \beta) \right) \right) = 3 \beta \} \\
\label{wtwfdesc} \wtwf & = & \{ \beta \in \wedge^2 ; \stph
\wedge \beta = 0 \} \\
\nonumber & = &
\{\beta \in \wedge^2; \st (\ph \wedge \beta) =  -\beta \} \\
\nonumber & = & \{\sum a_{ij} e^i \wedge e^j ; (a_{ij}) \in \lieg
\}
\end{eqnarray}

\begin{eqnarray}
\label{wfisdesc} 
\wfis & = & \{ \alpha \wedge \stph ; \alpha \in \wons \} \\
\nonumber & = & \{\gamma \in \wedge^5; \ph \wedge \st \gamma = 2
\gamma \} \\ \nonumber & = &
\{\gamma \in \wedge^5;  \stph \wedge \left( \st (\stph
\wedge \st \gamma) \right)  = 3 \gamma \} \\
\label{wfifdesc} \wfif & = & \{ \gamma \in \wedge^5 ; \ph
\wedge \st \gamma = -\gamma \} \\
\nonumber & = & \{ \gamma \in \wedge^5 ; \stph
\wedge \st \gamma = 0 \}
\end{eqnarray}

Notice that these subspaces are $+2$ and $-1$ eigenspaces of the
operators $L(\beta) = \st( \phi \wedge \beta)$ on $\wedge^2$ and
$M(\gamma) = \phi \wedge \st \gamma$ on $\wedge^5$. From this fact
we get the following useful formulas for the projections $\pi_k$
onto the $k$-dimensional representations, for
$\beta \in \wedge^2$ and $\gamma \in \wedge^5$:

\begin{eqnarray} \label{projectioneq} \st (\ph \wedge \beta) & = &
2 \pi_7(\beta) - \pi_{14}(\beta) \\
\nonumber \pi_7(\beta) & = & \frac{\beta + \st(\ph \wedge
\beta)}{3} \\ \nonumber
\pi_{14}(\beta) & = & \frac{2 \beta - \st(\ph \wedge \beta)}{3}
\end{eqnarray}
and
\begin{eqnarray*} \ph \wedge \st \gamma & =
& 2 \pi_7(\gamma) - \pi_{14}(\gamma) \\
 \pi_7(\gamma) & = & \frac{\gamma + \ph \wedge
\st \gamma}{3} \\ 
\pi_{14}(\gamma) & = & \frac{2 \gamma - \ph \wedge \st
\gamma}{3}
\end{eqnarray*}

We now move on to the decompositions for $k=3,4.$
\begin{eqnarray}
\label{wthodesc} \wtho & = & \{ f \ph; f \in C^{\infty}(M)\} \\
\nonumber & = & \{\eta \in \wedge^3; \ph \wedge \left( \st (\stph
\wedge \eta) \right) = 7 \eta \} \\ \label{wthsdesc} \wths & = &
\{\st(\ph \wedge \alpha); \alpha \in \wons\} \\ \nonumber & = & \{
w \hk \stph;  w \in \Gamma(T(M)) \} \\ \nonumber & = &  \{\eta \in
\wedge^3; \st \left( \ph \wedge \st (\ph \wedge \eta) \right) = -4
\eta \} \\ \label{wthtdesc} \wtht & = & \{\eta \in \wedge^3; \ph
\wedge \eta = 0 \text{ and } \stph \wedge \eta = 0\}
\end{eqnarray}

\begin{eqnarray}
\label{wfoodesc} \wfoo & = & \{ f \stph; f \in C^{\infty}(M)\} \\
\nonumber & = & \{\sigma \in \wedge^4; \stph \wedge \left( \st
(\ph \wedge \sigma) \right) = 7 \sigma \} \\ \label{wfosdesc} \wfos
& = & \{\ph \wedge \alpha; \alpha \in \wons\} \\ \nonumber & = & 
\{\sigma \in \wedge^4; \left( \ph \wedge \st (\ph \wedge \st
\sigma) \right) = -4 \sigma \} \\ \label{wfotdesc} \wfot & = &
\{\sigma \in \wedge^4; \ph \wedge \sigma = 0 \text{ and } \ph
\wedge \st \sigma = 0\}
\end{eqnarray}

\subsection{The metric of a \Gs} \label{Gmetricsec}

From Lemma~\ref{onewedgelemma}, we can obtain a formula for
determining the metric
$g$ from the $3$-form $\ph$:
\begin{prop} \label{metricprop}
If $v$ is a vector field on $M$, then
\begin{equation} \label{metriceq}
(v \hk \ph) \wedge ( v \hk \ph ) \wedge \ph = 6 {|v|}^2 \vol
\end{equation}
\end{prop}
\begin{proof}
From Lemma~\ref{intextidentities} and~\eqref{iso2eq} we
have
\begin{equation*}
v \hk \ph = \st ( \vs \wedge \stph )
\end{equation*}
and 
\begin{equation*}
( v \hk \ph ) \wedge \ph = 2 ( \vs \wedge \stph )
\end{equation*}
Thus we obtain
\begin{equation*}
(v \hk \ph) \wedge ( v \hk \ph ) \wedge \ph = 2 {|\vs \wedge
\stph|}^2 \vol = 6 {|v|}^2 \vol
\end{equation*}
where we have used~\eqref{onewedgestpheq}.
\end{proof}
By polarizing~\eqref{metriceq} in $v$, we obtain the relation:
\begin{equation*}
( v \hk \ph ) \wedge (w \hk \ph) \wedge \ph =
6 {\langle v, w \rangle} \vol 
\end{equation*}
From this equation we can obtain the metric.

\begin{lemma} \label{Gmetricprelemma}
Fix a vector field $v = v^k e_k$, where $e_1, e_2, \ldots, e_7$ is
an oriented local frame of vector fields. The expression obtained
from $v$ by
\begin{equation} \label{Gmetricprelemmaeq}
\frac{ { \left( (v \hk \ph) \wedge (v \hk \ph) \wedge \ph
\right) (e_1, e_2, \ldots, e_7) } }{ {\left( \det{\left(
\left( (e_i \hk \ph) \wedge (e_j \hk \ph) \wedge \ph
\right) (e_1, e_2, \ldots, e_7) \right)} \right)}^{\frac{1}{9}}}
\end{equation}
is homogeneous of order $2$ in $v$, and is independent of the
choice of $e_1, e_2, \ldots, e_7$. As shown in the next theorem, up
to a constant this is ${|v|}^2$.
\end{lemma}
\begin{proof}
The homogeneity of~\eqref{Gmetricprelemmaeq} of order $2$ in $v$
is clear. Now suppose we choose a different oriented basis
$e'_1, e'_2, \ldots, e'_7$. Then we have
\begin{equation*}
e'_i = P_{ij} e_j
\end{equation*}
and hence
\begin{equation*}
(e'_i \hk \ph) \wedge (e'_j \hk \ph) \wedge \ph = P_{ik} P_{jl}
(e_k \hk \ph) \wedge (e_l \hk \ph) \wedge \ph 
\end{equation*}
Hence in the new basis the denominator
of~\eqref{Gmetricprelemmaeq} changes by a factor of
\begin{equation*}
{\left( \det(P)^2 \det(P)^7 \right)}^{\frac{1}{9}} = \det(P)
\end{equation*}
and the numerator also changes by a factor of $\det(P)$, leaving
the quotient invariant.
\end{proof}

We can now give the expression for the metric in terms of the
$3$-form $\ph$.
\begin{thm} \label{Gmetrictheorem}
Let $v$ be a tangent vector at a point $p$ and let $e_1, e_2,
\ldots, e_7$ be any basis for $T_p M$. Then the length $|v|$ of
$v$ is given by
\begin{equation} \label{Gmetricthmeq}
{|v|}^2 = 6^{-\frac{2}{9}} \frac{ { \left( (v \hk \ph) \wedge (v \hk
\ph) \wedge \ph \right) (e_1, e_2, \ldots, e_7) } }{ {\left(
\det{\left( \left( (e_i \hk \ph) \wedge (e_j \hk \ph) \wedge \ph
\right) (e_1, e_2, \ldots, e_7) \right)} \right)}^{\frac{1}{9}}}
\end{equation}
\end{thm}
\begin{proof}
We work in local coordinates at the point $p$. In this notation
$g_{ij} = {\langle e_i, e_j \rangle}$ with $1 \leq i,j \leq 7$.
Let $\det (g)$ denote the determinant of $(g_{ij})$. We have
from~\eqref{metriceq} that
\begin{eqnarray*}
\left( (e_i \hk \ph) \wedge (e_j \hk \ph) \wedge \ph
\right) & = & 6 g_{ij} \vol \\ & = & 6 g_{ij} \sqrt{\det (g)} e^1
\wedge e^2 \wedge \ldots \wedge e^7 \\ \det{ \left( \left( (e_i \hk
\ph) \wedge (e_j \hk \ph) \wedge \ph \right) (e_1, e_2, \ldots,
e_7) \right)} & = & 6^7 \det (g) {\det (g)}^{\frac{7}{2}} \\ &
= & 6^7 {\det (g)}^{\frac{9}{2}}
\end{eqnarray*}
and since
\begin{eqnarray*}
(v \hk \ph) \wedge (v \hk \ph) \wedge \ph & = & 6 {|v|}^2 \vol
\\ & = & 6 {|v|}^2  \sqrt{\det (g)} e^1 \wedge e^2 \wedge \ldots
\wedge e^7 \\ \left( (v \hk \ph) \wedge (v \hk \ph) \wedge \ph
\right) (e_1, e_2, \ldots, e_7) & = & 6 {|v|}^2 {\det
(g)}^{\frac{1}{2}}
\end{eqnarray*}
these two expressions can be combined to
yield~\eqref{Gmetricthmeq}.
\end{proof}

\subsection{The cross product of a \Gs}

In this section we will describe the cross product operation on a
manifold with a \Gs\ in terms of the $3$-form $\ph$, and present
some useful relations.

\begin{defn} \label{Gcrossproductdefn}
Let $u$ and $v$ be vector fields on $M$. The {\em cross product},
denoted $u \times v$, is a vector field on $M$ whose associated
$1$-form under the metric isomorphism satisfies:
\begin{equation} \label{crossproducteq}
{(u \times v)}^{\flat} = v \hk u \hk \ph
\end{equation}
Notice that this immediately yields the relation between $\times$,
$\ph$, and the metric $g$:
\begin{equation} \label{gphirelationeq}
g ( u \times v , w) = {( u \times v )}^{\flat} (w) = w \hk v \hk u
\hk \ph = \ph (u,v,w).
\end{equation}
\end{defn}

Another characterization of the cross product is obtained from this
one using Lemma~\ref{intextidentities}:
\begin{eqnarray} \label{Gcrossproductchareq}
{(u \times v)}^{\flat} & = & v \hk u \hk \ph \\ \nonumber & = & -
\st ( \vs \wedge \st ( u \hk \ph ) ) \\ \nonumber & = & - \st ( \vs
\wedge \us \wedge \stph ) \\ \nonumber & = & \st ( \us \wedge \vs
\wedge \stph )
\end{eqnarray}
Now since $\us \wedge \vs$ is a $2$-form, we can write it as
$\beta_7 + \beta_{14}$, with $\beta_j \in \wedge^2_j$. Then we
have, using~\eqref{wtwsdesc} and~\eqref{wtwfdesc}: 
\begin{eqnarray} \label{Gbeta7eq}
{(u \times v)}^{\flat} \wedge \stph & = & \st ( \beta_7 \wedge
\stph ) \wedge \stph \\ \nonumber & = & 3 \st \beta_7
\end{eqnarray}
Taking the norm of both sides, and using~\eqref{onewedgestpheq}:
\begin{equation*}
{| {(u \times v)}^{\flat} \wedge \stph |}^2 = 3 {| {(u \times
v)}^{\flat} |}^2 = 3 {| {u \times v} |}^2 = 9 {| \beta_7 |}^2
\end{equation*}
from which we obtain
\begin{equation} \label{Gbeta7normeq}
{| \beta_7 |}^2 = \frac{1}{3} {| {u \times v} |}^2
\end{equation}
\begin{lemma} \label{Gcrossproductnormlemma}
Let $u$ and $v$ be vector fields. Then
\begin{equation} \label{Gcrossproductnormeq}
{| u \times v |}^2 = { | u \wedge v |}^2
\end{equation}
\end{lemma}
\begin{proof}
With $\beta = \us \wedge \vs$, we have from~\eqref{wtwsdesc}
and~\eqref{wtwfdesc}:
\begin{eqnarray*}
\beta \wedge \ph & = & 2 \st \beta_7 - \st \beta_{14} \\ 
\beta \wedge \beta \wedge \ph & = & 2 {|\beta_7|}^2 \vol -
{|\beta_{14}|}^2 \vol \\ & = & 0
\end{eqnarray*}
since $\beta = \us \wedge \vs$ is decomposable. So
${|\beta_{14}|}^2 = 2{|\beta_7|}^2$ and finally we obtain
from~\eqref{Gbeta7normeq}:
\begin{equation*}
{| u \times v |}^2 = 3 {|\beta_7|}^2 = {|\beta_7|}^2 +
{|\beta_{14}|}^2 = {|\beta|}^2 = {| u \wedge v |}^2   
\end{equation*}
\end{proof}
The next lemma describes the iteration of the cross product.
\begin{lemma} \label{cross7iteratelemma}
Let $u$, $v$, and $w$ be vector fields. Then we have
\begin{equation} \label{cross7iterateeq}
u \times \left( v \times w \right) = - \langle u, v \rangle w +
\langle u, w \rangle v + {\left( u \hk v \hk w \hk \stph
\right) }^{\sharp}
\end{equation}
\end{lemma}
\begin{proof}
From~\eqref{Gcrossproductchareq}, we have
\begin{equation*}
{\left( u \times \left( v \times w \right) \right) }^{\flat} = \st
\left( \us \wedge \st \left( \vs \wedge \ws \wedge \stph \right)
\wedge \stph \right)
\end{equation*}
Now since $\beta \wedge \stph = 0$ for $\beta \in \wtwf$, we can
replace $\vs \wedge \ws$ by $\pi_7 ( \vs \wedge \ws ) =
\frac{1}{3} \left( \vs \wedge \ws + \st \left( \ph \wedge \vs
\wedge \ws \right) \right)$. Then using~\eqref{wtwsdesc}, we have
\begin{eqnarray*}
\st \left( \pi_7 (\vs \wedge \ws) \wedge \stph \right) \wedge \stph
& = & 3 \st \pi_7 (\vs \wedge \ws)  \\ & = & \st \left( \vs
\wedge \ws + \st \left( \ph \wedge \vs \wedge \ws \right) \right)
\\ & = & \st \left( \vs \wedge \ws - v \hk w \hk \stph \right)
\end{eqnarray*}
which we substitute back above to obtain
\begin{eqnarray*}
{\left( u \times \left( v \times w \right) \right) }^{\flat} & = &
\st
\left( \us \wedge \st \left( \vs \wedge \ws - v \hk w \hk \stph
\right) \right) \\ & = & - u \hk \left( \vs \wedge \ws - v \hk w
\hk \stph \right) \\ & = & - \langle u,v \rangle \ws + \langle u,w
\rangle \vs + u \hk v \hk w \hk \stph
\end{eqnarray*}
which completes the proof.
\end{proof}
\begin{cor} \label{cross7iteratespecialcor}
For the special case $u = v$, we obtain the useful relation
\begin{equation} \label{cross7iteratespecialeq}
u \times ( u \times w ) = - {|u|}^2 w + \langle u, w \rangle u
\end{equation}
\end{cor}
\begin{rmk}
Note that if we fix a {\em unit vector field} $u$, then the
operator $J w = u \times w$ when restricted to the space of vector
fields orthogonal to $u$, satisfies $J^2 = -I$, and hence is an
almost complex structure on the orthogonal complement of $u$.
\end{rmk}

The following lemma will be used in
Section~\ref{vfdeformsec} to determine how the metric changes
under a deformation in the $\wths$ direction.
\begin{lemma} \label{Gcrossidentitylemma}
The following identity holds for $v$ and $w$ vector fields:
\begin{equation} \label{Gcrossidentityeq}
( v \hk w \hk \stph ) \wedge (v \hk w \hk \ph) \wedge \stph = 2
{|v \wedge w|}^2 \vol
\end{equation}
\end{lemma}
\begin{proof}
We start with Lemma~\ref{intextidentities} to rewrite
\begin{eqnarray*}
v \hk w \hk \stph & = & \st ( \vs \wedge \st (w \hk \stph) ) \\ &
= & - \st ( \vs \wedge \ws \wedge \ph ) \\ & = & -2 \beta_7 +
\beta_{14}
\end{eqnarray*}
using the notation as above. From equations~\eqref{crossproducteq}
and~\eqref{Gbeta7eq} we have
\begin{equation*}
(v \hk w \hk \ph) \wedge \stph = -3 \st \beta_7
\end{equation*}
Combining these two equations and~\eqref{Gbeta7normeq},
\begin{eqnarray*}
( v \hk w \hk \stph ) \wedge (v \hk w \hk \ph) \wedge \stph  & =
& (-2 \beta_7 + \beta_{14} ) \wedge (-3 \st \beta_7 ) \\ & = & 6
{|\beta_7|}^2 \vol \\ & = & 2 {| u \wedge v |}^2 \vol
\end{eqnarray*}
which completes the proof.
\end{proof}

Finally, we prove a theorem which will be useful in
Section~\ref{vfdeformsec} where we will use it to show that to
first order, deforming a \Gs\ by an element of $\wths$ does not
change the metric.
\begin{thm} \label{w223zerothm}
Let $u$, $v$, $w$ be vector fields. Then
\begin{equation*}
\left(u \hk \ph \right) \wedge \left( v \hk \ph \right) \wedge
\left( w \hk \stph \right) = 0.
\end{equation*}
Note that in terms of the decompositions in~\eqref{wtwsdesc}
and~\eqref{wthsdesc}, this theorem says that the wedge product map
\begin{equation*}
\wtws \times \wtws \times \wths \to \wseo
\end{equation*}
is the zero map.
\end{thm}
\begin{proof}
Since it is an $8$-form,
\begin{equation*}
\left( u \hk \ph \right) \wedge \left( v \hk \ph \right) \wedge
\stph = 0.
\end{equation*}
Taking the interior product with $w$ and rearranging,
\begin{eqnarray*}
\left(u \hk \ph \right) \wedge \left( v \hk \ph \right) \wedge
\left( w \hk \stph \right)  & = & -\left( w \hk u \hk \ph \right)
\wedge \left( v \hk \ph \right) \wedge \stph \\ & & {}- \left( u
\hk \ph \right) \wedge \left( w \hk v \hk \ph \right) \wedge
\stph
\end{eqnarray*}
Now using~\eqref{hk3eq}, we get
\begin{equation*}
\left(u \hk \ph \right) \wedge \left( v \hk \ph \right) \wedge
\left( w \hk \stph \right) = -3 \left( w \hk u \hk \ph \right)
\wedge \st \vs -3 \left( w \hk v \hk \ph \right) \wedge \st \us
\end{equation*}
Finally, from~\eqref{usefuleq}, we have
\begin{eqnarray*}
\left(u \hk \ph \right) \wedge \left( v \hk \ph \right) \wedge
\left( w \hk \stph \right) & = & -3 \left( u \hk \ph \right)
\wedge \st \left( \ws \wedge \vs \right)  -3 \left( v \hk \ph
\right) \wedge \st \left( \ws \wedge \us \right) \\ & = & -3 \ph
\wedge \st \left( \us \wedge \ws \wedge \vs \right) -3 \ph
\wedge \st \left( \vs \wedge \ws \wedge \us \right) \\ & = & 0. 
\end{eqnarray*}
\end{proof}

\subsection{The 16 classes of \Gs s}

According to the classification of Fern\'andez and Gray
in~\cite{FG}, a manifold with a \Gs\ has holonomy a subgroup of
$\G$ if and only if $\nabla \ph = 0$, which they showed to be
equivalent to
\begin{equation*}
d \ph = 0 \qquad \text{ and } \qquad d \stph = 0.
\end{equation*}
They established this equivalence by decomposing the space $W$ that
$\nabla \ph$ belongs to into irreducible $\G$-representations, and
identifying the invariant subspaces of $W$ with isomorphic
subspaces of $\wedge^* (M)$. This space $W$ decomposes as
\begin{equation*}
W = W_1 \oplus W_7 \oplus W_{14} \oplus W_{27}
\end{equation*}
where the subscript $k$ denotes the dimension of the irreducible
representation $W_k$. Now $d \ph \in \wtho \oplus \wths \oplus
\wtht$ and $d \stph \in \wfis \oplus \wfif$. Up to isomorphism,
the projections $\pi_k(d \ph)$ and $\pi_k(d \stph)$ are non-zero
constant multiples of $\pi_k (\nabla \ph)$. Therefore in the
following we will consider $d \ph$ and $d \stph$ instead of
$\nabla \ph$. Since both of these have a component in a
$7$-dimensional representation, they are multiples:
\begin{lemma}
The following identity holds:
\begin{equation} \label{definemueq}
\mu = \st d \ph \wedge \ph = - \st d \stph \wedge \stph
\end{equation}
where we have defined the $6$-form $\mu$ by the above two
equal expressions. They are the components $\pi_7(d\ph)$
and $\pi_7(d \stph)$ transferred to the isomorphic space
$\wsis$.
\end{lemma}
\begin{proof}
See ~\cite{Br2} for a proof of this fact.
\end{proof}
We prefer to work with the associated $1$-form, $\theta = \st
\mu$. We will see that in some subclasses this $1$-form is
closed or at least ``partially closed.''

Now we say a \Gs\ is in the  class $W_i \oplus W_j
\oplus W_k$ with $i,j,k$ distinct where $\{i,j,k\} \subset
\{1,7,14,27\}$ if only the component of $d\ph$ or $d\stph$ in the
$l$-dimensional representation vanishes. Here $\{l\} =
\{1,7,14,27\} \setminus \{i,j,k\}$. Similarly the \Gs\
is in the class $W_i \oplus W_j$ if the $k$ and $l$-dimensional
components vanish, and in the class $W_i$ if the other three
components are zero. In this way we arrive at $16$ classes of \Gs
s on a manifold. In Table~\ref{Gtable} we describe the classes in
terms of differential equations on the form $\ph$. This
classification first appeared in~\cite{FG} and then in essentially
this form in~\cite{C}.
\begin{table}[h]
\begin{center}
\begin{tabular}{|c|c|l|c|} 
\hline Class & Defining Equations & Name & $d \theta$ \\
\hline $W_1 \oplus W_7 \oplus W_{14} \oplus W_{27}$ & no relation
on $d \ph$, $d\stph$. & & \\ \hline $W_7 \oplus W_{14} \oplus
W_{27}$ & $d \ph \wedge \ph = 0$ & & $d \theta = ?$ \\ \hline
$W_1 \oplus W_{14} \oplus W_{27}$ & $\theta = 0$ &  &
$\theta = 0$ \\ \hline $W_1 \oplus W_7 \oplus W_{27}$ & $d \stph +
\frac{1}{3} \theta \wedge \stph = 0$ & ``integrable'' &
$\pi_7(d\theta) = 0$ \\ & or \ $\ph \wedge ( \st d \stph) = - 2
d \stph$ & & \\ \hline $W_1 \oplus W_7 \oplus W_{14}$ &
$d\ph + \frac{1}{4} \theta \wedge \ph - h \stph = 0$ &  &
$d \theta = ?$ \\ \hline $W_{14} \oplus W_{27}$ & $d\ph
\wedge \ph = 0$ and $\theta = 0$ & & $\theta = 0$ \\ \hline $W_7
\oplus W_{27}$ & $d\ph \wedge \ph = 0$ and & & $\pi_7(d \theta)$ =
0 \\ & $d \stph + \frac{1}{3} \theta \wedge \stph = 0$ & & \\
\hline $W_7 \oplus W_{14}$ & $d\ph + \frac{1}{4} \theta \wedge \ph
= 0$ & LC almost \G & $d \theta = 0$ \\ \hline $W_1 \oplus W_{27}$
& $d \stph = 0$ & semi-\G & $\theta = 0$ \\ \hline $W_1 \oplus
W_{14}$ & $d \ph - h \stph = 0$ & & $\theta = 0$ \\ \hline $W_1
\oplus W_7$ & $d \ph + \frac{1}{4} \theta \wedge \ph - h
\stph = 0$ & LC nearly \G & $d \theta = 0$ \\ & and \ $d \stph +
\frac{1}{3} \theta \wedge \stph = 0$ & & \\ \hline $W_{27}$ & $d\ph
\wedge \ph = 0$ and $d\stph = 0$ & & $\theta = 0$ \\ \hline
$W_{14}$ & $d \ph = 0$ & almost \G & $\theta = 0$ \\ \hline $W_7$
&  $d \stph + \frac{1}{3} \theta \wedge \stph = 0$ & LC \G & $d
\theta = 0$ \\ & and \ $d \ph + \frac{1}{4} \theta \wedge \ph
= 0$ & & \\ \hline $W_1$ & $d \ph - h \stph = 0$ and $d\stph = 0$
& nearly \G & $\theta = 0$ \\ \hline $\{0\}$ & $d\ph = 0$ and
$d\stph = 0$ & \G & $\theta = 0$ \\ \hline
\end{tabular}
\end{center}
\caption{The $16$ classes of \Gs s \label{Gtable}}
\end{table}

In Table~\ref{Gtable}, the function $h = \frac{1}{7} \st ( \ph
\wedge d \ph )$ is the image of $\pi_1(d \ph)$ in $\wedge^0$
under the isomorphism $\wfoo \cong \wzeo$. The abbreviation ``LC''
stands for {\em locally conformal to} and means that for those
classes, we can (at least locally) conformally change the metric
to enter a strictly smaller subclass. This will be explained
in Section~\ref{Gconformaldeformsec}.

We now prove the closedness or partial closedness of
$\theta$ in the various classes as given in the final column of
Table~\ref{Gtable}. The closedness of $\theta$ in the classes $W_1
\oplus W_7$ and $W_7 \oplus W_{14}$ was originally shown using a
different approch by Cabrera in~\cite{C}.
\begin{lemma} \label{dthetalemma}
If $\ph$ is in the classes $W_7$, $W_7 \oplus W_{14}$, or $W_1
\oplus W_7$, then $d\theta = 0$. Furthermore, if $\ph$ is in the
classes $W_7 \oplus W_{27}$ or $W_1 \oplus W_7 \oplus W_{27}$ then
$\pi_7 (d\theta) = 0$.
\end{lemma}
\begin{proof}
We begin by showing that if $\ph$ satisfies $d \ph + \frac{1}{4} \theta \wedge \ph
= 0$, then $d \theta = 0$, and if $\ph$ satisfies $d
\stph + \frac{1}{3} \theta \wedge \stph = 0$, then $\pi_7(d
\theta) = 0$.

Suppose $d \ph + \frac{1}{4} \theta \wedge \ph
= 0$. We differentiate this equation to obtain:
\begin{equation*}
d \theta \wedge \ph = \theta \wedge d \ph = \theta \wedge
\left( -\frac{1}{4} \theta \wedge \ph \right) = 0
\end{equation*}
But wedge product with $\ph$ is an isomorphism from $\wedge^2$ to
$\wedge^5$, so $d \theta = 0$. Now suppose $d \stph + \frac{1}{3}
\theta \wedge \stph = 0$. Differentiating this equation yields
\begin{equation*}
d \theta \wedge \stph = \theta \wedge d \stph = \theta \wedge
\left (-\frac{1}{3} \theta \wedge \stph \right) = 0
\end{equation*}
But wedge product with $\stph$ is an isomorphism from
$\wtws$ to $\wsis$, so $\pi_7(d \theta) = 0$. 

Thus by comparing with Table~\ref{Gtable}, we have shown that in
the classes $W_7 \oplus W_{14}$ and $W_7$, we have $d \theta = 0$.
Also, in the classes $W_1 \oplus W_7 \oplus W_{27}$ and $W_7
\oplus W_{27}$ we have $\pi_7( d\theta ) = 0$. We still have to
show that $\theta$ is closed in the class $W_1 \oplus W_7$. We
already have that $\pi_7(d\theta)=0$, so we need only show that
$\pi_{14}(d \theta) = 0$ in this case. We differentiate $d \ph +
\frac{1}{4} \theta \wedge \ph - h \stph = 0$ to obtain
\begin{eqnarray*}
0 & = & \frac{1}{4} d\theta \wedge \ph - \frac{1}{4} \theta \wedge
d \ph - dh \wedge \stph - h d \stph \\ & = & \frac{1}{4}d\theta
\wedge \ph - \frac{1}{4} \theta \wedge \left(-\frac{1}{4} \theta
\wedge \ph + h \stph \right) - dh \wedge \stph - h
\left( -\frac{1}{3} \theta \wedge \stph \right) \\ & = &
\frac{1}{4} d\theta \wedge \ph + \alpha \wedge \stph
\end{eqnarray*}
for some $1$-form $\alpha$, where we have used the fact that $d
\stph + \frac{1}{3} \theta \wedge \stph = 0$ in this class. But
$\alpha \wedge \stph$ is in $\wfis$, and since wedge product with
$\ph$ is an isomorphism from $\wedge^2_k$ to $\wedge^5_k$ for
$k=7,14$, this shows that $\pi_{14}(d\theta) = 0$.
\end{proof}

The inclusion relations among these various subclasses are
analyzed in~\cite{FG, C, CMS, F2, Br2, BS, J1, J2, Sa}.
For all but one case, examples can be found of manifolds which
are in a particular class but {\em not} in a strictly smaller
subclass. For example, a manifold in the class $W_{14}$ which does
{\em not} have holonomy \G\ appears in~\cite{F2}. There is one
case of an inclusion in Table~\ref{Gtable} which is {\em
not} strict. This is given by the following result, which first
appeared in~\cite{C}.
\begin{prop} \label{nonstrictinclusionprop}
The class $W_1 \oplus W_{14}$ equals $W_1 \cup W_{14}$ exactly. 
\end{prop}
\begin{proof}
In the class $W_1 \oplus W_{14}$, we have $d \ph - h \stph = 0$
(and by consequence $\theta = 0$). Differentiating this equation,
\begin{equation*}
d h \wedge \stph = - h d \stph
\end{equation*}
If $h \neq 0$, then by dividing by $h$ and using
Proposition~\ref{isomorphisms}, we see that $d \stph \in
\wfis$, so $\pi_{14} (d \stph) = 0$. But since we already have
that $\theta = 0$, this means $d \stph = 0$ and hence $\ph$ is
actually of class $W_1$ (nearly \G). If $h = 0$ then $d \ph = 0$
and $\ph$ is of class $W_{14}$ (almost \G). 
\end{proof}

\begin{rmk} \label{nearlyGrmk}
Note that in the proof of the above proposition, we see that if
$\ph$ is of class $W_1$ (nearly \G), then $d h \wedge \stph = 0$,
and so $d h = 0$ by Proposition~\ref{isomorphisms}. Therefore in
the nearly \G\ case, the function $h$ is {\em locally
constant}, or constant if the manifold $M$ is connected.
In~\cite{Gr6} Gray showed that all nearly \G\ manifolds are
actually {\em Einstein}.
\end{rmk}

In~\cite{FU1, FU2}, Fern\'andez and Ugarte show that for
manifolds with a \Gs\ in the classes $W_1 \oplus W_7 \oplus W_{27}$
(``integrable'') or $W_7 \oplus W_{14}$, there exists a subcomplex
of the deRham complex. They then show how to define analogues of
Dolbeault cohomology of complex manifolds in these two cases,
including analogues of $\bar \partial$-harmonic forms. They derive
properties of these cohomology theories and topological
restrictions on the existence of \Gs s in some strictly smaller
subclasses.

\section{Deformations of a fixed \Gs} \label{deformsec}

Let us begin with a fixed \Gs\ on a manifold $M$ in a certain
class. We are interested in how deforming the form $\ph$ affects
the class. In other words, we are interested in what kinds of
deformations preserve which classes of \Gs s. Now since $\ph \in
\wtho \oplus \wths \oplus \wtht$, there are three canonical ways
to deform $\ph$. For example, since $\wtho = \{ f \ph \}$, adding
to $\ph$ an element of $\wtho$ amounts to conformally scaling
$\ph$. This preserves the decomposition into irreducible
representations in this case. However, since the decomposition
does depend on $\ph$ (unlike the decomposition of forms into
$(p,q)$ types on a K\"ahler manifold) in general if we add an
element of $\wths$ or $\wtht$ the decomposition does change. So
deforming in those two directions really only makes sense
infinitesmally. However, we shall attempt to get as far as we can
with an actual deformation before we restrict to infinitesmal
deformations.

\subsection{Conformal Deformations of \Gs s}
\label{Gconformaldeformsec}

Let $f$ be a smooth, nowhere vanishing function on $M$. For
notational convenience, which will become evident, we will
conformally scale $\ph$ by $f^3$. Let the new form $\nph = f^3
\oph$. We first compute the new metric $\tilde g$ and the new
volume form $\nvol$ in the following lemma.

\begin{lemma} \label{conformaldeformnewmetric}
The metric $\og$ on vector fields, the metric $\og^{-1}$ on one
forms, and the volume form $\ovol$ transform as follows:
\begin{eqnarray*}
\nvol & = & f^7 \ovol \\ \tilde g & = & f^2 \og \\ \tilde g^{-1}
& = & f^{-2} \og^{-1} 
\end{eqnarray*}
\end{lemma}
\begin{proof}
Using Proposition~\ref{metricprop}, we have in a local
coordinate chart:
\begin{eqnarray*}
\tilde g (u,v) \nvol & = & \frac{1}{6}(u \hk \nph) \wedge (v \hk
\nph) \wedge \nph \\ & = & f^9 \og (u,v) \ovol \\ \tilde g (u,v)
\sqrt{\det(\tilde g)} d x^1\ldots d x^7 & = & f^9 \og (u,v)
\sqrt{\det(\og)} d x^1 \ldots d x^7.
\end{eqnarray*}
Thus, taking determinants of the coefficients of both sides,
\begin{eqnarray*}
{\det(\tilde g)}^{\frac{7}{2}} \det(\tilde g) & = & f^{63}
{\det(\og)}^{\frac{7}{2}} \det(\og) \\ \sqrt{\det(\tilde g)} & = &
f^7 \sqrt{\det(\og)}
\end{eqnarray*}
This gives $\nvol = f^7 \ovol$, from which we can immediately see
that $\tilde g = f^2 \og$ and $\tilde g^{-1} = f^{-2} \og$.
\end{proof}
Notice that the new metric $\tilde g$ is always a positive
definite metric as long as $f$ is non-vanishing, even if $f$ is
negative. However, the {\em orientation} of $M$ changes for
negative $f$ since the volume form changes sign. Now we can
determine the new Hodge star $\nst$ in terms of the old $\ost$.
\begin{lemma} \label{conformaldeformnewstar}
If $\alpha$ is a $k$-form, then $\nst \alpha = f^{7 - 2 k} \ost
\alpha$.
\end{lemma}
\begin{proof}
Let $\alpha$, $\beta$ be $k$-forms. Then from
Lemma~\ref{conformaldeformnewmetric} the new metric on $k$-forms
is $\ltn{  \left \langle  \ ,\ \right \rangle } = f^{-2 k} { \left
\langle  \ ,\ \right \rangle}_{\! \text{o}}$. From this we compute:
\begin{eqnarray*}
\beta \wedge \nst \alpha & = & \ltn{ \left \langle \beta, \alpha
\right \rangle } \nvol \\ & = & f^{-2k} \left \langle \beta,
\alpha \right \rangle _{\! \text{o}} f^7 \ovol \\ & = & f^{7 - 2k}
\beta \wedge \ost \alpha.
\end{eqnarray*}
\end{proof}
\begin{cor} \label{conformaldeformstphi}
The new form $\nph$ satisfies $\nstph = f^4 \ostph$.
\end{cor}
\begin{proof}
This follows immediately from Lemma~\ref{conformaldeformnewstar}.
\end{proof}

Combining all our results so far yields:

\begin{lemma} \label{conformaldeformprelemma}
We have the following relations:
\begin{eqnarray*}
d\nph & = & 3 f^2 d f \wedge \oph + f^3 d \oph \\ d \nstph & = & 4
f^3 d f \wedge \ostph + f^4 d \ostph \\ \nst d \nph & = & 3 f
\ost( d f \wedge \oph) + f^2 \ost d \oph \\ \nst d \nstph & = & 4
\ost ( d f \wedge \ostph) + f \ost(d \ostph)
\end{eqnarray*}
\end{lemma}
\begin{proof}
This follows from Lemma~\ref{conformaldeformnewstar} and
Corollary~\ref{conformaldeformstphi}.
\end{proof}
Using these results, we can determine which classes of \Gs s are
conformally invariant. We can also determine what happens to the
$6$-form $\mu$ from equation~\eqref{definemueq} as well as the
associated $1$-form $\theta = \st \mu$. This is all given in the
following theorem:

\begin{thm} \label{conformaldeformresults}
Under the conformal deformation $\nph = f^3 \oph$, we have:
\begin{eqnarray}
\label{conformaldeformintegeq} \ninteg & = & f^4 \left( \ointeg
\right) \\ \label{conformaldeformintdualeq} \nintdual & = & f^3
\left( \ointdual \right) \\ \label{conformaldeformoneeq} d \nph
\wedge \nph & = & f^6 \left( d \oph \wedge \oph \right) \\
\label{conformaldeform27eq} \nintdual - \nh \nstph & = & f^3
\left( \ointdual - \oh \ostph \right)  \\
\label{conformaldeformmueq} \nmu & = & -12 f^4 \ost df + f^5 \omu
\\ \label{conformaldeformthetaeq} \ntheta & = & -12 d(\log(f)) +
\otheta
\end{eqnarray}
Hence, we see (from Table~\ref{Gtable}) that the classes which are
conformally invariant are exactly $W_7 \oplus W_{14} \oplus
W_{27}$, $W_1 \oplus W_7 \oplus W_{27}$, $W_1 \oplus W_7 \oplus
W_{14}$, $W_7 \oplus W_{27}$, $W_7 \oplus W_{14}$, $W_1 \oplus
W_7$, and $W_7$. These are precisely the classes which have a
$W_7$ component. (This conclusion was originally observed
in~\cite{FG} using a different method.)

Additionally, \eqref{conformaldeformthetaeq} shows that since
$\theta$ changes by an exact form, in the classes where $d \theta
= 0$, we have a well defined cohomology class $[\theta ]$ which is
unchanged under a conformal scaling. These are the classes $W_7
\oplus W_{14}$, $W_1 \oplus W_7$, and $W_7$.
\end{thm}

\begin{proof}
We begin by using Lemma~\ref{conformaldeformprelemma}
and~\eqref{definemueq} to compute $\nmu$ and $\ntheta$:
\begin{eqnarray*}
\nmu & = & \nst d \nph \wedge \nph \\ & = & \left( 3 f
\ost( d f \wedge \oph) + f^2 \ost d \oph \right) \wedge
f^3 \oph \\ & = & 3 f^4 \oph \wedge \ost \left( \oph \wedge
d f \right) + f^5 \omu \\ & = & - 12 f^4 \ost d f + f^5 \omu
\end{eqnarray*}
where we have used~\eqref{isopheq} in the last step. Now
from Lemma~\ref{conformaldeformnewstar}, we get:
\begin{equation*}
\ntheta = \nst \nmu = -12 f^{-1} df + \otheta =  -12 d(\log(f)) +
\otheta.
\end{equation*}
Now using the above expression for $\ntheta$, we have:
\begin{eqnarray*}
\ninteg & = & 4 f^3 d f \wedge \ostph + f^4 d \ostph +
\frac{1}{3} \left( -12 f^{-1} df + \otheta \right) \wedge f^4
\ostph \\ & = & f^4 \left( \ointeg \right) \\ \nintdual & = &  3
f^2 d f \wedge \oph + f^3 d \oph + \frac{1}{4} \left( -12 f^{-1}
df + \otheta \right) \wedge f^3 \oph \\ & = & f^3 \left( \ointdual
\right)
\end{eqnarray*}
and finally, since $\oph \wedge \oph = 0$, 
\begin{equation*}
d\nph \wedge \nph = \left( 3 f^2 d f \wedge \oph + f^3 d \oph
\right) \wedge f^3 \oph = f^6 \left( d\oph \wedge \oph \right).
\end{equation*}
Finally, since $h = \frac{1}{7} \st ( \ph \wedge d \ph )$, we have
\begin{eqnarray*}
\nh \nstph & = & \frac{1}{7} \nst ( \nph \wedge d \nph ) f^4
\ostph \\ & = & \frac{1}{7} f^{-7} \ost ( f^6 \oph \wedge d \oph )
f^4 \ostph \\ & = & f^3 \oh \ostph
\end{eqnarray*}
which yields~\eqref{conformaldeform27eq} when combined
with~\eqref{conformaldeformintdualeq}. This completes the proof.
\end{proof}

These results now enable us to give necessary and sufficient
conditions for obtaining a closed or co-closed $\nph$ by
conformally scaling the original $\oph$.
\begin{thm} \label{Gconformaldeformfinal}
Let $\oph$ be a positive $3$-form (associated to a \Gs). Under the
conformal deformation $\nph = f^3 \oph$, the new $3$-form $\nph$
satisfies
\begin{itemize}
\item $d \nph = 0$ $\iff$ $\oph$ is at least class $W_7
\oplus W_{14}$ and $12 d \log(f) = \otheta$. \item $d \ostph = 0$
$\iff$ $\oph$ is at least class $W_1 \oplus W_7 \oplus
W_{27}$ and $12 d \log(f) = \otheta$.
\end{itemize}
Note that in both cases, in order to have $\nph$ be closed or
co-closed after conformal scaling, the original $1$-form $\otheta$
has to be {\em exact}. In particular if the manifold is
simply-connected or more generally $H^1(M) = 0$ then this will
always be the case if $\oph$ is in the classes $W_7 \oplus
W_{14}$, $W_1 \oplus W_7$, or $W_7$, where $d \otheta = 0$. 
\end{thm}
\begin{proof}
From Lemma~\ref{conformaldeformprelemma}, for $d \nph = 0$, we need
\begin{eqnarray*}
d \nph & = & 3 f^2 d f \wedge \oph + f^3 d \oph = 0 \\ \Rightarrow
\qquad d \oph & = & -3 d \log (f) \wedge \oph
\end{eqnarray*}
which says that $d \oph \in \wfos$ by
Proposition~\ref{isomorphisms}. Hence $\pi_1(d \oph)$ and
$\pi_{27}(d \oph)$ both vanish and $\oph$ must be already at least
of class $W_7 \oplus W_{14}$. Then to make $d \nph = 0$, we need to
eliminate the $W_7$ component, which requires $12 d \log(f) =
\otheta$ by Theorem~\ref{conformaldeformresults}. Similarly, to
make $d \nstph = 0$, Lemma~\ref{conformaldeformprelemma} gives
\begin{eqnarray*}
d \nstph & = & 4 f^3 d f \wedge \ostph + f^4 d \ostph = 0 \\
\Rightarrow \qquad d \ostph & = & -4 d \log(f) \wedge \ostph
\end{eqnarray*}
which says $d \ostph \in \wfis$ and $\pi_{14} (d \ostph) = 0$ by
Proposition~\ref{isomorphisms}. Thus $\oph$ must already be at
least class $W_1 \oplus W_7 \oplus W_{27}$ and we need to choose
$f$ by $12 d \log(f) = \otheta$ to scale away the $W_7$ component.
\end{proof}

\begin{rmk} \label{Gconformaldeformrmk}
We have shown that the transformation $\nph = f^3 \oph$ stays in
a particular subclass as long as there is a $W_7$ component to that
class. If there is, and the original $\otheta$ is exact, then we
can choose $f$ to scale away the $W_7$ component and enter a
stricter subclass. Conversely, Theorem~\ref{conformaldeformresults}
shows that a conformal scaling by a non-constant $f$ will always
generate a non-zero $W_7$ component if we started with none. Hence,
if we are trying to constuct metrics of holonomy \G\ on a
simply-connected manifold, it is enough to construct a metric in
the class $W_7$, since we can then conformally scale (uniquely) to
obtain a metric of holonomy \G. This is why the class $W_7$ is
called {\em locally conformal} \G.
\end{rmk}

\subsection{Deforming $\ph$ by an element of $\wths$}
\label{vfdeformsec}

The type of deformation of $\ph$ that is next in line in terms of
increasing complexity is to add an element of $\wths$. This space
is isomorphic to $\wons \cong \Gamma (T(M))$, so we can think of
this process as deforming $\ph$ by a vector field. In fact, an
element $\eta \in \wths$ is of the form $w \hk \stph$ for some
vector field $w$, by~\eqref{wthsdesc}. Let $\nph = \oph +
t w \hk \ostph$, for $t \in \mathbb R$. We will develop formulas
for the new metric $\tilde g$, the new Hodge star $\nst$, and other
expressions entirely in terms of the old $\oph$, the old $\ost$,
and the vector field $w$. Note in this case the background
decomposition into irreducible $\G$-representations changes, and
in Section~\ref{infvfdeformsec} we will linearize by taking $\left.
\frac{d}{dt} \right|_{t=0}$ of our results.

\begin{lemma} \label{vfdeformmetriclemma}
In the expression
\begin{equation*}
6 \lt{|v|}^2 \nvol = \left( v \hk \nph \right) \wedge \left( v \hk
\nph \right) \wedge \nph
\end{equation*}
which is a cubic polynomial in $t$, the linear and cubic terms
both vanish, and the coefficient of the quadratic term is
\begin{equation*}
6 {|v \wedge w|}^2_{\! \text{o}} \ovol
\end{equation*}
\end{lemma}
\begin{proof}
The coefficient of $t^3$ is:
\begin{equation*}
\left(v \hk w \hk \ostph \right) \wedge \left(v \hk w \hk \ostph
\right) \wedge \left( w \hk \ostph \right) 
\end{equation*}
This expression is zero because it arises by taking the interior
product with $w$ of the $8$-form
\begin{equation*}
\left(v \hk w \hk \ostph \right) \wedge \left(v \hk w \hk \ostph
\right) \wedge \ostph = 0.
\end{equation*}
The coefficient of $t$ is:
\begin{equation*}
\left(v \hk \oph \right) \wedge \left(v \hk \oph \right) \wedge
\left( w \hk \ostph \right) + 2 \left(v \hk \oph \right)
\wedge \left(v \hk w \hk \ostph \right) \wedge \oph
\end{equation*}
Using~\eqref{usefuleq2} on the second term and rearranging, this
coefficient becomes
\begin{equation*}
3 \left(v \hk \oph \right) \wedge \left(v \hk \oph \right) \wedge
\left( w \hk \ostph \right)
\end{equation*}
which vanishes by Theorem~\ref{w223zerothm}.

The coefficient of $t^2$ is:
\begin{equation} \label{quadtermeq}
(v \hk w \hk \ostph) \wedge \left( (v \hk w \hk \ostph) \wedge \oph
+ 2 (v \hk \oph) \wedge (w \hk \ostph) \right)
\end{equation}
Applying~\eqref{usefuleq2} twice and rearranging, this
coefficient becomes 
\begin{equation*}
3 \left(v \hk w \hk \ostph \right) \wedge \left(v \hk w \hk \oph
\right) \wedge \ostph
\end{equation*}
The statement now follows from Lemma~\ref{Gcrossidentitylemma}.
This completes the proof.
\end{proof}

Before we can use Lemma~\ref{vfdeformmetriclemma} to obtain the
new metric, we have to extract the new volume form.
\begin{prop} \label{vfdeformnvolprop}
With $\nph = \oph + w \hk \ostph$, the new volume form is
\begin{equation} \label{vfdeformnvoleq}
\nvol = \opw^{\frac{2}{3}} \ovol
\end{equation}
\end{prop}
\begin{proof}
We work in local coordinates. Let $e_1, e_2, \ldots, e_7$ be a
basis for the tangent space, with $w = w^j e_j$, $g_{ij} =
{\langle e_i, e_j \rangle}_{\! \text{o}}$ and $\tilde{g}_{ij} =
\ltn{ \langle e_i, e_j \rangle}$. Then
Lemma~\ref{vfdeformmetriclemma} says that
\begin{equation*}
\lt{|v|}^2 \sqrt{\det(\tilde{g})} = \left( {|v|}^2_{\text{o}}
+ {|v \wedge w|}^2_{\text{o}} \right) \sqrt{\det{(g)}}
\end{equation*}
Polarizing this equation, we have:
\begin{eqnarray*}
\ltn{ \left \langle v_1,v_2 \right \rangle }
\sqrt{\det(\tilde{g})} & = & ( { \left \langle v_1,v_2 \right
\rangle }_{\! \text{o}} + { \left \langle v_1,v_2 \right \rangle
}_{\! \text{o}}{|w|}^2_{\text{o}} - { \left \langle v_1,w \right
\rangle }_{\! \text{o}} { \left \langle v_2,w \right \rangle
}_{\! \text{o}} ) \sqrt{\det{(g)}} \\ \tilde{g}_{ij}
\sqrt{\det(\tilde{g})} & = & \left( g_{ij} + {\langle
e_i \wedge w , e_j \wedge w \rangle }_{\! \text{o}} \right)
\sqrt{\det{(g)}}
\end{eqnarray*}
with $v_1 = e_i$ and $v_2 = e_j$. Now substituting $w = w^k e_k$
in the second term,
\begin{equation*}
{\langle e_i \wedge w , e_j \wedge w \rangle }_{\! \text{o}} = 
{|w|}^2_{\text{o}} g_{ij} - w_i w_j
\end{equation*}
Thus we have
\begin{equation*}
\tilde{g}_{ij} \sqrt{\det(\tilde{g})} = \left( g_{ij}(1 +
{|w|}^2_{\text{o}}) - w_i w_j \right) \sqrt{\det{(g)}}
\end{equation*}
We take determinants of both sides of this equation, and use the
fact that they are $7 \times 7$ matrices, to obtain
\begin{equation} \label{vfnvoltempeq}
{\left( \det{(\tilde{g})} \right)}^{\frac{9}{2}} =
{\left( \det{(g)} \right)}^{\frac{7}{2}} \det{\left( g_{ij}(1 +
{|w|}^2_{\text{o}}) - w_i w_j \right)}.
\end{equation}
Using Lemma~\ref{specialdetlemma}, the determinant on the right is
\begin{equation} \label{vfnvoltemp2eq}
{(1 + {|w|}^2_{\text{o}})}^7 \det (g) - {|w|}^2_{\text{o}} {(1
+ {|w|}^2_{\text{o}})}^6 \det (g) = {(1 + {|w|}^2_{\text{o}})}^6
\det (g)
\end{equation}
Substituting this result into equation~\eqref{vfnvoltempeq}, we
obtain
\begin{eqnarray*}
{\left( \det{(\tilde{g})} \right)}^{\frac{9}{2}} & = & 
{\left( \det{(g)} \right)}^{\frac{7}{2}} {(1
+ {|w|}^2_{\text{o}})}^6 \det (g) \\ \sqrt{ \det{(\tilde{g})}}
& = &  {(1 + {|w|}^2_{\text{o}})}^{\frac{2}{3}} \sqrt{
\det{(g)}}
\end{eqnarray*}
which completes the proof.
\end{proof}

Now letting $t=1$, with $\nph = \oph
+w \hk \ostph$, Lemma~\ref{vfdeformmetriclemma} and
Proposition~\ref{vfdeformnvolprop} yield
\begin{eqnarray*}
\lt{|v|}^2 \nvol & = & \left( {|v|}^2_{\text{o}} + {|v \wedge
w|}^2_{\text{o}} \right) \ovol \\
\ltn{ \left \langle v,v \right \rangle } & = & 
\frac{1}{{\opw}^{\frac{2}{3}}} ( { \left \langle v,v \right
\rangle }_{\! \text{o}} + {|v|}^2_{\text{o}}{|w|}^2_{\text{o}} - {
\left \langle v,w \right \rangle }^2_{\! \text{o}} )
\end{eqnarray*}
Polarizing this equation, we obtain:
\begin{equation} \label{Gvfdeformnewmetriceq}
\ltn{ \left \langle v_1,v_2 \right \rangle } =
\frac{1}{{\opw}^{\frac{2}{3}}} ( { \left \langle v_1,v_2 \right
\rangle }_{\! \text{o}} + { \left \langle v_1,v_2 \right \rangle
}_{\! \text{o}}{|w|}^2_{\text{o}} - { \left \langle v_1,w \right
\rangle }_{\! \text{o}} { \left \langle v_2,w \right \rangle
}_{\! \text{o}} )
\end{equation}
which by Lemma~\ref{Gcrossproductnormlemma} can also be written as
\begin{equation} \label{Gvfdeformnewmetriceq2}
\lt{\langle v_1, v_2 \rangle} = \frac{1}{{\opw}^{\frac{2}{3}}}
\left( {\langle v_1, v_2 \rangle}_{\! \text{o}} + {\langle w
\times v_1, w \times v_2 \rangle}_{\! \text{o}} \right)
\end{equation}
Note that in the above expression $\times$ refers to the vector
cross product associated to the initial \Gs\ $\oph$. Later we will
describe this metric geometrically.
In local coordinates with $w = w^i e_i$, $g_{ij} = { \left
\langle e_i,e_j \right \rangle  }_{\! \text{o}}$, and $\ws = w_i
e^i$, we see that
\begin{equation} \label{vfdeformnewmetriceq}
\tilde g_{ij} = \frac{1}{{\opw}^{\frac{2}{3}}} ( g_{ij} \opw - w_i
w_j )
\end{equation}
\begin{prop} \label{vfdeformnewinvmetricprop}
In local coordinates, the metric $\tilde g^{ij}$ on
$1$-forms is given by:
\begin{equation*}
\tilde g^{ij} = \frac{1}{{\opw}^{\frac{1}{3}}}(g^{ij} + w^i w^j )
\end{equation*}
\end{prop}
\begin{proof}
We compute:
\begin{eqnarray*}
\tilde g_{ij} \tilde g^{jk} & = & \frac{1}{{\opw}^{\frac{2}{3}}} (
g_{ij} \opw - w_i w_j )  \frac{1}{{\opw}^{\frac{1}{3}}}(g^{jk} +
w^j w^k ) \\ & = &  \frac{1}{\opw}\left( (g_{ij} g^{jk} +
g_{ij} w^j w^k) \opw - g^{jk} w_i w_j - w_i w_j w^j w^k \right) \\
& = &  \frac{1}{\opw}\left( (\delta_i^k + w_i w^k) \opw - w_i
w^k - |w|^2_{\text{o}} w_i w^k \right) \\ & = & \delta_i^k
\end{eqnarray*}
which completes the proof.
\end{proof}

Now with $\alpha = \alpha_i e^i$ and $\beta = \beta_j e^j$ two
$1$-forms, their new inner product is
\begin{eqnarray}
\nonumber \ltn{ \left \langle \alpha, \beta \right \rangle } =
\alpha_i \beta_j \tilde g^{ij} & = & 
\frac{1}{{\opw}^{\frac{1}{3}}}(\alpha_i \beta_j g^{ij} +
\alpha_i w^i \beta_j w^j ) \\ \label{vfdeformnewinvmetriceq} & =
& \frac{1}{{\opw}^{\frac{1}{3}}}( { \left \langle \alpha, \beta
\right \rangle }_{\! \text{o}} + (w \hk \alpha ) ( w \hk \beta ) )
\end{eqnarray}
From this expression we can derive a formula for the new metric
$\ltn{ \left \langle \ ,\  \right \rangle }$ on $k$-forms:
\begin{thm} \label{vfdeformnewformmetricthm}
Let $\alpha$, $\beta$ be $k$-forms. Then
\begin{equation} \label{vfdeformnewformmetriceq}
\ltn{ \left \langle \alpha, \beta \right \rangle } =
\frac{1}{{\opw}^{\frac{k}{3}}}(  { \left \langle \alpha, \beta
\right \rangle }_{\! \text{o}} + { \left \langle w \hk \alpha, w
\hk \beta \right \rangle }_{\! \text{o}} )
\end{equation}
\end{thm}
\begin{proof}
We have already established it for the case $k=1$
in~\eqref{vfdeformnewinvmetriceq}, and the case $k=0$ is trivial.
For the general case, we will prove the statement on {\em
decomposable} forms and it follows in general by linearity. Let
$\alpha = e^{i_1} \wedge e^{i_2} \wedge \ldots \wedge e^{i_k}$ and
$\beta = e^{j_1} \wedge e^{j_2} \wedge \ldots \wedge e^{j_k}$. Then
by the definition of the metric on $k$-forms,
\begin{equation*}
\ltn{ \left \langle \alpha, \beta \right \rangle } = \det
{\begin{pmatrix} 
\ltn{ \left \langle e^{i_1},e^{j_1} \right \rangle } & \ltn{
\left \langle e^{i_1},e^{j_2} \right \rangle } & \ldots &
\ltn{ \left \langle e^{i_1},e^{j_k} \right \rangle } \\ \ltn{
\left \langle e^{i_2},e^{j_1} \right \rangle } & \ltn{ \left
\langle e^{i_2},e^{j_2} \right \rangle } & \ldots & \ltn{ \left
\langle e^{i_2},e^{j_k} \right \rangle } \\ \vdots & & \ddots &
\vdots \\ \ltn{ \left \langle e^{i_k},e^{j_1} \right \rangle } &
\ltn{ \left \langle e^{i_k},e^{j_2} \right \rangle } & \ldots &
\ltn{ \left \langle e^{i_k},e^{j_k} \right \rangle }
\end{pmatrix} }
\end{equation*}
Now from equation~\eqref{vfdeformnewinvmetriceq} each entry in the
above matrix is of the form
\begin{equation*}
\ltn{ \left \langle e^{i_a},e^{j_b} \right \rangle } =
\frac{1}{{\opw}^{\frac{1}{3}}}(g^{i_a j_b} + w^{i_a} w^{j_b} ) 
\end{equation*}
and we have
\begin{equation*}
\ltn{ \left \langle \alpha, \beta \right \rangle } =
\frac{1}{{\opw}^{\frac{k}{3}}} \det {\begin{pmatrix} g^{i_1 j_1} +
w^{i_1} w^{j_1} & \ldots & g^{i_1 j_k} + w^{i_1} w^{j_k} \\ \vdots
& \ddots & \vdots \\ g^{i_k j_1} + w^{i_k} w^{j_1} & \ldots &
g^{i_k j_k} + w^{i_k} w^{j_k} 
\end{pmatrix}}
\end{equation*}
Now we apply Lemma~\ref{detlemma} to obtain
\begin{equation*}
{ \left \langle \alpha, \beta \right \rangle }_{\! \text{o}} +
\sum_{l,m=1}^k {(-1)^{l+m} w^{i_l} w^{j_m} { \left \langle 
e^{i_1} \wedge \ldots \widehat {e^{i_l}}  \ldots \wedge e^{i_k},
e^{j_1} \wedge \ldots \widehat {e^{j_m}} \ldots \wedge e^{j_k}
\right \rangle _{\! \text{o}} } }
\end{equation*}
for the determinant above. Now with $w = w^i e_i$, we can take the
interior product with both $\alpha$ and $\beta$:
\begin{eqnarray*}
w \hk \alpha & = & \sum_{l=1}^k (-1)^{l-1} w^{i_l} e^{i_1} \wedge
\ldots \widehat {e^{i_l}} \ldots \wedge e^{i_k} \\ w \hk \beta & =
 & \sum_{m=1}^k (-1)^{m-1} w^{j_m} e^{j_1} \wedge
\ldots \widehat {e^{j_m}} \ldots \wedge e^{j_k}
\end{eqnarray*} 
and hence the sum over $l$ and $m$ above is just ${ \left \langle
w \hk \alpha, w \hk \beta \right \rangle }_{\! \text{o}} $. Putting
everything together, we  arrive at~\eqref{vfdeformnewformmetriceq}:
\begin{equation*}
\ltn{ \left \langle \alpha, \beta \right \rangle } =
\frac{1}{{\opw}^{\frac{k}{3}}}(  { \left \langle \alpha, \beta
\right \rangle }_{\! \text{o}} + { \left \langle w \hk \alpha, w
\hk \beta \right \rangle }_{\! \text{o}} )
\end{equation*}
\end{proof}

To continue our analysis of the new \Gs \ $\nph$, we now need to
compute the new Hodge star $\nst$.
\begin{thm} \label{vfdeformnewstarthm}
The Hodge star for the new metric on a $k$-form $\alpha$ is given
by:
\begin{eqnarray}
\label{vfdeformnewstareq} \nst \alpha & = & {\opw}^{\frac{2-k}{3}}
\left( \ost \alpha + (-1)^{k-1} w \hk \left( \ost ( w \hk \alpha )
\right) \right) \\ \nonumber & = & {\opw}^{\frac{2-k}{3}} \left(
\ost \alpha + w \hk ( \ws \wedge \ost \alpha ) \right)
\end{eqnarray}
\end{thm}
\begin{proof}
The second form follows from the first from~\eqref{ident1}.
Although it looks a little more cluttered, we prefer to use the
first form for $\nst$. Notice that up to a scaling factor, the new
star is given by `twisting by $w$', taking the old star, then
`untwisting by $w$', and adding this to the old star. To establish
this formula, let $\beta$ be an arbitrary $k$-form and compute:
\begin{eqnarray*}
\beta \wedge \nst \alpha & = & \ltn{ \left \langle \beta, \alpha
\right \rangle } \nvol \\ & = &  \frac{1}{{\opw}^{\frac{k}{3}}}( 
{ \left \langle \alpha, \beta \right \rangle }_{\! \text{o}} + {
\left \langle w \hk \alpha, w \hk \beta \right \rangle }_{\!
\text{o}} ) \opw^{\frac{2}{3}} \ovol \\ & = &
{\opw}^{\frac{2-k}{3}} \left( \beta \wedge \ost \alpha + (w \hk
\beta ) \wedge \ost (w \hk \alpha) \right) 
\end{eqnarray*}
Now if we take the interior product with $w$ of the $8$-form
\begin{equation*}
\beta \wedge \ost (w \hk \alpha) = 0
\end{equation*}
we obtain
\begin{equation*}
(w \hk \beta) \wedge \ost ( w \hk \alpha ) = (-1)^{k-1} \beta
\wedge \left( w \hk \left( \ost ( w \hk \alpha ) \right) \right)
\end{equation*}
and this completes the proof, since $\beta$ is arbitrary.
\end{proof}

At this point before continuing it is instructive to observe
directly that $\nst^2 = 1$, as of course it should on a
$7$-manifold. It clarifies the necessity of all the factors of
$\opw$. Let $\alpha$ be a $k$-form:
\begin{eqnarray*}
\nst \alpha & = & {\opw}^{\frac{2-k}{3}}
\left( \ost \alpha + (-1)^{k-1} w \hk \left( \ost ( w \hk \alpha )
\right) \right) \\ \nst ( \nst \alpha ) & = & 
{\opw}^{\frac{2-(7 -k)}{3}}
\left( \ost (\nst \alpha) + (-1)^{(7-k)-1} w \hk \left( \ost ( w
\hk (\nst \alpha) ) \right) \right) \\ & = & {\opw}^{\frac{(k - 5)
+ (2 - k)}{3}} \left[ \ost (\ost \alpha) + (-1)^{k-1} \ost \left(w
\hk \left( \ost ( w \hk ( \alpha) ) \right) \right) \right. \\ & &
\hskip 0.5in \left. {} + (-1)^k w \hk \left( \ost ( w \hk ( \ost
\alpha) ) \right) + (-1) w \hk \left( \ost ( w \hk w\hk \ldots )
\right) \right]
\end{eqnarray*}
Now the last term is zero because of the two successive interior
products with $w$. Using $\ost^2 = 1$ we now have:
\begin{eqnarray*}
\nst ( \nst \alpha ) & = & {\opw}^{-1}
\left[ \alpha + (-1)^{k-1} \ost \left(w
\hk \left( \ost ( w \hk ( \alpha) ) \right) \right)  
  {} + (-1)^k w \hk \left( \ost ( w \hk ( \ost
\alpha) ) \right) \right]
\end{eqnarray*}
Now using equation~\eqref{ww2eq}, (with $n=7$) we get
\begin{equation*}
\nst (\nst \alpha) =  {\opw}^{-1}
\left( \alpha + {|w|}_{\text{o}}^2 \alpha \right) = \alpha
\end{equation*}

We now give a geometric interpretation of the transformation $\oph
\mapsto \oph + w \hk \ostph$. From~\eqref{Gvfdeformnewmetriceq}
for the new metric $\tilde g$, with $v_1 = v$ and $v_2 = w$, we
have
\begin{eqnarray*} 
\ltn{ \left \langle v,w \right \rangle } & = &
\frac{1}{{\opw}^{\frac{2}{3}}} ( { \left \langle v,w \right
\rangle }_{\! \text{o}} + { \left \langle v,w \right \rangle
}_{\! \text{o}}{|w|}^2_{\text{o}} - { \left \langle v, w \right
\rangle }_{\! \text{o}} { \left \langle w,w \right \rangle
}_{\! \text{o}} ) \\ & = & \frac{1}{{\opw}^{\frac{2}{3}}} {
\left \langle v,w \right \rangle }_{\! \text{o}} 
\end{eqnarray*}
Hence we see that all the distances are {\em shrunk} by a factor
of ${\opw}^{-\frac{2}{3}}$ in the direction of the vector field
$w$. On the other hand, if either $v_1$ or $v_2$ is orthogonal to
$w$ in the old metric, then~\eqref{Gvfdeformnewmetriceq}
gives
\begin{eqnarray*}
\ltn{ \left \langle v_1,v_2 \right \rangle } & = &
\frac{1}{{\opw}^{\frac{2}{3}}} ( { \left \langle v_1,v_2 \right
\rangle }_{\! \text{o}} + { \left \langle v_1,v_2 \right \rangle
}_{\! \text{o}}{|w|}^2_{\text{o}} - 0 ) \\ & = &
{\opw}^{\frac{1}{3}} { \left \langle v_1,v_2 \right \rangle }_{\!
\text{o}} 
\end{eqnarray*}
Thus in the directions perpendicular to the vector field $w$, the
distances are {\em stretched} by a factor of
${\opw}^{\frac{1}{3}}$. Therefore this new metric is expanded in
the $6$ directions perpendicular to $w$ and is compressed in the
direction parallel to $w$. Of course, the situation is more
complicated if neither $v_1$ nor $v_2$ is parallel or
perpendicular to $w$. This produces a {\em tubular} manifold. For
example in the case of $M = N \times S^1$, where
$N$ is a Calabi-Yau $3$-fold and the metric on $M$ is the product
metric, if we take $w = \frac{\partial}{\partial \theta}$ where
$\theta$ is a coordinate on $S^1$, then the Calabi-Yau manifold
$N$ is expanded and the circle factor $S^1$ is compressed under
$\oph \mapsto \oph + w \hk \ostph$. By replacing $w$ by $t w$ and
letting $t \to \infty$, we can make this ``tube'' as long and thin
as we want. The total volume, however, always increases by
${\opw}^{\frac{2}{3}}$ by Proposition~\ref{vfdeformnvolprop}.

In general, determining the class of \Gs\ that $\nph$ belongs to
for $\nph = \oph +w \hk \oph$ involves some very complicated
differential equations on the vector field $w$. However, since
$\nph$ is {\em always} a positive $3$-form for any $w$, it may be
interesting to study some of these differential equations in the
simplest cases to determine if one can choose $w$ to produce a
$\nph$ in a strictly smaller subclass. From
Theorem~\ref{vfdeformnewstarthm}, we have
\begin{eqnarray}
\nonumber \nstph & = & {\opw}^{-\frac{1}{3}} \left( \ost \nph + w
\hk ( \ost ( w \hk \nph ) ) \right) \\ \label{vfdeformnstarphieq} &
= & {\opw}^{-\frac{1}{3}} \left( \ostph + \ost ( w \hk \ostph ) + w
\hk \ost ( w \hk \oph ) \right)
\end{eqnarray}
For example this transformation will yield a manifold of holonomy
\G\ if $w$ satisfies the system
\begin{eqnarray*}
0 & = & d( \oph + w \hk \ostph) \\ 0 & = & d \left(
{\opw}^{-\frac{1}{3}} \left( \ostph + \ost ( w \hk \ostph ) + w \hk
\ost ( w \hk \oph ) \right) \right)
\end{eqnarray*}
The ellipticity and other properties of this system under certain
hypotheses is currently being investigated~\cite{K1}.

\subsection{Infinitesmal deformations in the $\wths$ direction}
\label{infvfdeformsec}

Now since the decomposition of the space of differential forms
corresponding to the \Gs\ $\ph$ {\em changes} when we add
something in $\wths$, we could also consider a
one-parameter family $\ph_t$ of \Gs s satisfying
\begin{equation} \label{infvfdeformeq}
\frac{\partial}{\partial t} \ph_t = w \hk \st_t \ph_t
\end{equation}
for a fixed vector field $w$. That is, at each time $t$, we move in
the direction $w \hk \st_t \ph_t$ which is a $3$-form in
$\wedge^3_{7_t}$, the decomposition depending on $t$. Since the
Hodge star $\st_t$ is also changing in time, this is {\em a priori}
a nonlinear equation. However, our first observation is that this
is in fact not the case:
\begin{prop} \label{infvfdeformlinprop}
Under the flow described by equation~\eqref{infvfdeformeq}, the
metric $g$ {\em does not change}. Hence the volume form and
Hodge star are also constant.
\end{prop}
\begin{proof}
From~\eqref{metriceq} which gives the metric from the
$3$-form, we have:
\begin{equation*}
g_t(u,v) \vol_t = \frac{1}{6} (u \hk \ph_t ) \wedge ( v \hk \ph_t
) \wedge \ph_t
\end{equation*}
Differentiating with respect to $t$, and using the differential
equation~\eqref{infvfdeformeq}, 
\begin{eqnarray*}
6 \frac{\partial}{\partial t} \left( g_t(u,v) \vol_t \right) & =
& (u \hk w \hk \st_t \ph_t ) \wedge (v \hk \ph_t) \wedge \ph_t +
(u \hk \ph_t) \wedge (v \hk w \hk \st_t \ph_t ) \wedge \ph_t \\ &
& {} + (u \hk \ph_t) \wedge (v \hk \ph_t) \wedge (w \hk \st_t
\ph_t)
\end{eqnarray*}
Now from the proof of Lemma~\ref{vfdeformmetriclemma} (the linear
term) we see that this expression is zero, by polarizing. From
this it follows easily by taking determinants that $\vol_t$ is
constant and thus so is $g_t$ and $\st_t$.
\end{proof}
Therefore we can replace $\st_t$ by $\st_0 = \st$ and
equation~\eqref{infvfdeformeq} is actually {\em linear}. Moreover,
the flow determined by this linear equation gives a one-parameter
family of \Gs s each yielding the {\em same} metric $g$. Our
equation is now
\begin{equation*}
\frac{\partial}{\partial t} \ph_t = w \hk \st \ph_t = A \ph_t
\end{equation*}
where $A$ is the linear operator $\alpha \mapsto A \alpha = w \hk
\st \alpha$ on $\wedge^3$.
\begin{prop} \label{skewopprop}
The operator $A$ is {\em skew-symmetric}. Furthermore, the
eigenvalues $\lambda$ of $A$ and their multiplicities
$N_{\lambda}$ are:
\begin{eqnarray*}
\lambda = 0 & & N_{\lambda} = 21 \\ \lambda = i |w| & &
N_{\lambda} = 7 \\ \lambda = - i |w| & & N_{\lambda} = 7
\end{eqnarray*}
\end{prop}
\begin{proof}
Let $e^1, e^2, \ldots e^{35}$ be a basis of $\wedge^3$. Then
\begin{eqnarray*}
A_{ij} \vol & = &  \left \langle e^i, A e^j \right \rangle \vol \\
& = & e^i \wedge \st (w \hk \st e^j)
\\ & = & - e^i \wedge \ws \wedge e^j \\ & = & \ws \wedge e^i
\wedge e^j \\ & = & - A_{ji} \vol
\end{eqnarray*}
since $3$-forms anti-commute. Therefore $A$ is diagonalizable over
$\mathbb C$. Suppose now that $\alpha \in \wedge^3$ is an
eigenvector with eigvenalue $\lambda = 0$. Then
\begin{eqnarray*}
A \alpha & = & w \hk \st \alpha \\ & = & - \st (\ws \wedge
\alpha)  \\ & = & 0  
\end{eqnarray*}
so $\ws \wedge \alpha = 0$ and hence $\alpha = \ws \wedge \beta$
for some
$\beta \in \wedge^2$. Therefore the multiplicity of $\lambda = 0$
is $\dim (\wedge^2) = 21$. If $A \alpha = \lambda \alpha$ for
$\lambda \neq 0$, then $\alpha = \frac{1}{\lambda}(w \hk \st
\alpha)$ and $w \hk \alpha = 0$. Then we can write~\eqref{ww2eq} as
\begin{equation*}
|w|^2 \alpha =  - w \hk \st ( w \hk \st \alpha ) = - A^2 \alpha
= - {\lambda}^2 \alpha
\end{equation*}
and hence $\lambda = \pm i |w|$. Since the eigenvalues come in
complex conjugate pairs and there are $35-21 = 14$ remaining, there
must be $7$ of each. This completes the proof.
\end{proof}
Now if $\alpha$ is an eigenvector for $\frac{\partial}{\partial
t} \alpha_t = A \alpha_t = \lambda \alpha_t$, then $\alpha(t) =
e^{\lambda t} \alpha(0)$. Let $u_1, u_2, \ldots, u_{21}$ be a
basis for the $\lambda = 0$ eigenspace, and $v_1, \ldots, v_7$
and $\bar v_1, \ldots, \bar v_7$ be bases of {\em complex}
eigenvectors corresponding to the $\lambda = +i |w|$ and $\lambda =
- i |w|$ eigenspaces, respectively. We can write 
\begin{eqnarray*}
\phi_0 & = & \sum_{k=1}^7 c_k v_k + \sum_{k=1}^7 \bar c_k \bar v_k
+ \sum_{k=1}^{21} h_k u_k \\ & = & \sum_{k=1}^7 c_k v_k +
\sum_{k=1}^7 \bar c_k \bar v_k + \eta_0
\end{eqnarray*}
where $\eta_0$ as defined by the above equation is the part of
$\ph_0$ in the kernel of $A$. Then the solution is given by
\begin{eqnarray} \nonumber
\ph_t & = & \sum_{k=1}^7 c_k e^{i |w| t} v_k + \sum_{k=1}^7 \bar
c_k e^{- i |w| t} \bar v_k + \eta_0 \\ \nonumber & = & \cos(|w|t)
\sum_{k=1}^7 (c_k v_k + \bar c_k \bar v_k ) + \sin(|w|t)
\sum_{k=1}^7 i(c_k v_k - \bar c_k \bar v_k ) + \eta_0 \\
\label{infvfdeformprelimeq} & = & \cos(|w|t) \beta_0 + \sin(|w|t)
\gamma_0 + \eta_0
\end{eqnarray}
All that remains is to determine $\beta_0$, $\gamma_0$, and
$\eta_0$ in terms of the initial condition $\ph_0$. Substituting
$t=0$ into~\eqref{infvfdeformprelimeq}, we have
\begin{equation*}
\ph_0 = \beta_0 + \eta_0
\end{equation*}
Differentiating, we have
\begin{eqnarray*}
\frac{\partial}{\partial t} \ph_t & = & - |w|\sin(|w|t) \beta_0 +
|w|\cos(|w|t) \gamma_0 \\ A \ph_t & = & \cos(|w|t) A \beta_0 +
\sin(|w|t) A \gamma_0 + A \eta_0
\end{eqnarray*}
Comparing coefficients, we have
\begin{eqnarray*}
A \beta_0 & = & |w| \gamma_0 \\ A \gamma_0 & = & -
|w| \beta_0 \\ A \eta_0 & = & 0
\end{eqnarray*}
From $\beta_0 = \ph_0 - \eta_0$ and the equations above, we get
\begin{equation*}
\gamma_0 = \frac{1}{|w|} (A \ph_0 )
\end{equation*}
and substituting this into the second equation, we obtain
\begin{equation*}
\beta_0 = -\frac{1}{{|w|}^2}(A^2 \ph_0 )
\end{equation*}
Finally, we can state the general solution:
\begin{thm} \label{infvfdeformsolutionthm}
The solution to the differential equation
\begin{equation*}
\frac{\partial}{\partial t} \ph_t = w \hk \st_t \ph_t
\end{equation*}
is given by
\begin{equation} \label{infvfdeformsolutioneq}
\ph (t) = \ph_0 + \frac{1 - \cos(|w|t)}{{|w|}^2}(w \hk \st (w \hk
\st \ph_0 ) ) + \frac{\sin(|w|t)}{|w|} (w \hk \st \ph_0 )
\end{equation}
The solution exists for all time and is closed curve in
$\wedge^3$. Also, the path only depends on $\pm \frac{w}{|w|}$, and
the norm $|w|$ only affects the speed of travel along this curve.
\end{thm}
\begin{proof}
This is all immediate from the above discussion.
\end{proof}

\begin{rmk}
In~\cite{BS}, it is shown that the set of \Gs s on $M$ which
correspond to the same metric as that of a fixed \Gs\ $\oph$ is an
$\mathbb R \mathbb P^7$-bundle over the manifold $M$. The above
theorem gives an explicit formula~\eqref{infvfdeformsolutioneq} for
a path of \Gs s all corresponding to the same metric $g$ starting
from an arbitrary vector field $w$ on $M$.
\end{rmk}
\begin{rmk}
This can also be compared to the K\"ahler case. Since the metric
and the almost complex structure $J$ are independent in this case,
for a fixed metric $g$, the family of $2$-forms $\omega( \cdot,
\cdot) = g( J \cdot, \cdot)$ for varying $J$'s are all K\"ahler
forms corresponding to the same metric.
\end{rmk}
\begin{rmk}
Even though the metric is unchanged under an infinitesmal
deformation in the $\wths$ direction, the class of \Gs\ {\em
can} change. Therefore simply knowing that a metric on a
$7$-manifold arises from a \Gs\ and knowing the metric explicitly
does not determine the class.
\end{rmk}
\begin{rmk}
\textcolor{red}{Removed after publication. This was incorrect.}
\end{rmk}

We now apply this theorem to two specific examples, where we
will reproduce known results. 

\begin{ex} \label{infvfdeformCYex}
Let $N$ be a Calabi-Yau threefold, with K\"ahler form $\omega$ and
holomorphic $(3,0)$ form $\Omega$. The complex coordinates will
be denoted by $z^j = x^j + i y^j$. Then there is a natural
\Gs\ 
$\ph$ on the product $N \times S^1$ given by
\begin{equation} \label{cy3phieq}
\ph = \real (\Omega) + d\theta \wedge \omega
\end{equation}
where $\theta$ is the coordinate on the circle $S^1$. This induces
the product metric on $N \times S^1$, with the flat metric on
$S^1$. With the orientation on $N \times S^1$ given by
$(x^1,x^2,x^3,\theta,y^1,y^2,y^3)$, it is easy to check that
\begin{equation*}
\stph = - d\theta \wedge \imag ( \Omega ) + \frac{\omega^2}{2}
\end{equation*}
Now let $w = \frac{\partial}{\partial \theta}$ be a globally
defined non-vanishing vector field on $S^1$ with $|w|= 1$. Then we
have
\begin{eqnarray*}
w \hk \stph & = & - \imag (\Omega) \\ \st \left( w \hk \stph
\right) & = & - d\theta \wedge \real (\Omega) \\ w \hk \st \left( w
\hk \stph \right) & = & - \real (\Omega) 
\end{eqnarray*}
Thus for this choice of vector field $w$, the flow
in~\eqref{infvfdeformsolutioneq} is given by
\begin{eqnarray*}
\ph_t & = & \real (\Omega) + d\theta \wedge \omega - {(1 -
\cos(t))} \real (\Omega) - {\sin(t)} \imag (\Omega) \\ & = &
\cos(t) \real (\Omega) - \sin(t) \imag (\Omega) + d\theta \wedge
\omega \\ & = & \real( e^{i t} \Omega ) + d\theta \wedge \omega
\end{eqnarray*}
which is the canonical \G\ form on $N \times S^1$ where now the
Calabi-Yau structure on $N$ is given by $e^{i t} \Omega$ and
$\omega$. It is well-known that we have this freedom of changing
the holomorphic volume form $\Omega$ by a phase and preserving the
Ricci-flat metric. Here it arises naturally using the flow
described by~\eqref{infvfdeformsolutioneq} and the
canonical vector field $w = \frac{\partial}{\partial \theta}$.
\end{ex}

\begin{ex} \label{infvfdeformHKex}
As a second example, let $W$ be a $K3$ surface, which is
hyperK\"ahler with hyperK\"ahler triple $\omega_1$, $\omega_2$,
and $\omega_3$. With local complex coordinates $z^1 = y^0 + i y^1$
and $z^2 = y^2 + i y^3$ these forms can be written as:
\begin{eqnarray*}
\omega_1 & = & dy^0 \wedge dy^1 + dy^2 \wedge dy^3 \\ \omega_2 &
= & dy^0 \wedge dy^2 + dy^3 \wedge dy^1 \\ \omega_3 & = & dy^0
\wedge dy^3 + dy^1 \wedge dy^2 
\end{eqnarray*}
The volume form $\vol_W$ on $W$ is given by
$\frac{\omega_j^2}{2}$ for any $j = 1,2,3$. There is a natural
\Gs\ $\ph$ on the product $W \times T^3$ given by
\begin{equation} \label{k3phieq}
\ph = d\theta^1 \wedge d\theta^2 \wedge d\theta^3 -  d\theta^1
\wedge \omega_1 - d\theta^2 \wedge \omega_2 - d\theta^3 \wedge
\omega_3
\end{equation}
where $\theta^1,\theta^2,\theta^3$ are coordinates on the torus
$T^3$. This induces the product metric on $W \times T^3$, with the
flat metric on
$T^3$. With the orientation on $W \times T^3$ given by
$(\theta^1,\theta^2,\theta^3,y^0,y^1,y^2,y^3)$, it is easy to
check that
\begin{equation*}
\stph = \vol_W - d\theta^2 \wedge d\theta^3 \wedge \omega_1 -
d\theta^3 \wedge d\theta^1 \wedge \omega_2 - d\theta^1 \wedge
d\theta^2 \wedge \omega_3
\end{equation*}
Now let $w = \frac{\partial}{\partial \theta^1}$ be one of the 
globally defined non-vanishing vector fields on $T^3$ with
$|w|= 1$. Then we have
\begin{eqnarray*}
w \hk \stph & = & d \theta^3 \wedge \omega_2 - d \theta^2 \wedge
\omega_3 \\ \st \left( w \hk \stph \right) & = & d \theta^1 \wedge
d \theta^2 \wedge \omega_2 - d \theta^3 \wedge d\theta^1
\wedge \omega_3
\\ w \hk \st \left( w \hk \stph \right) & = & d \theta^2 \wedge
\omega_2 + d \theta^3 \wedge \omega_3
\end{eqnarray*}
Thus for this choice of vector field $w$, the flow
in~\eqref{infvfdeformsolutioneq} is given by
\begin{eqnarray*}
\ph_t & = & d\theta^1 \wedge d\theta^2 \wedge d\theta^3 - d\theta^1
\wedge \omega_1 - d\theta^2 \wedge \omega_2 - d\theta^3 \wedge
\omega_3 \\ & & {}+ {(1 - \cos(t))} \left( d \theta^2 \wedge
\omega_2 + d \theta^3 \wedge \omega_3 \right) + \sin(t) \left( d
\theta^3 \wedge \omega_2 - d \theta^2 \wedge \omega_3 \right) \\ &
= &  d\theta^1 \wedge d\theta^2 \wedge d\theta^3 - d\theta^1
\wedge \omega_1 \\ & & {}- d \theta^2 \wedge \left( \cos(t)
\omega_2 + \sin(t) \omega_3 \right) - d \theta^3 \wedge \left(
-\sin(t) \omega_2 + \cos(t) \omega_3 \right)
\end{eqnarray*}
which is the canonical \G\ form on $W \times T^3$ where now the
hyperK\"ahler triple on $W$ is given by $\tilde \omega_1 =
\omega_1$, $\tilde \omega_2 = \cos(t) \omega_2 + \sin(t)
\omega_3$, and $\tilde \omega_3 =  -\sin(t) \omega_2 + \cos(t)
\omega_3$. This is just a restatement of the fact that on a
hyperK\"ahler manifold, we have an $S^2$ worth of complex
structures, and we can choose any triple $I,J,K$ such that $I J =
K$ to obtain the three K\"ahler forms. The above construction
corresponds to a {\em hyperK\"ahler rotation} where $J \mapsto J
\cos(t) + K \sin(t)$ and $K \mapsto -J \sin(t) + K \cos(t)$. This
is rotation by an angle $t$ around the axis in $S^2$ that
represents the complex structure $I$. It is clear that we can
hyperK\"ahler rotate around any axis by taking $w =
a^1 \frac{\partial}{\partial \theta^1} + a^2
\frac{\partial}{\partial \theta^2} + a^3 \frac{\partial}{\partial
\theta^3}$ where $(a^1)^2 + (a^2)^2 + (a^3)^2 = 1$ in the flow
described by~\eqref{infvfdeformsolutioneq}. All these
hyperK\"ahler structures on $W$ yield the same metric, and hence
determine the same metric on $W \times T^3$, as expected.
\end{ex}

\section{Manifolds with a \SPs} \label{Spin7manifoldssec}

\subsection{\SPs s} \label{SPstructuressec}

Let $M$ be an oriented $8$-manifold with a global $3$-fold cross
product structure. Such a structure will henceforth be called
a \SPs. Its existence is also given by (this time more
complicated) topological conditions. (Again see~\cite{Gr5, J4, Sa}
for details.) Similarly to the \G\ case, this cross product $X(
\cdot, \cdot, \cdot )$ gives rise to an associated {\em Riemannian
metric}
$g$ and an alternating $4$-form $\Ph$ which are related by:
\begin{equation} \label{SPcompatibleeq}
\Ph(a,b,c,d) = g(X(a,b,c), d).
\end{equation}
As in the \G\ case, the metric and the cross product
structure cannot be prescribed independently. We will see in
Section~\ref{SPmetricsec} how the $4$-form $\Ph$ determines the
metric $g( \cdot, \cdot )$. For a \SPs\ $\Ph$, near a point $p \in
M$ we can choose local coordinates $x^0, x^1 \ldots, x^7$ so that
{\em at the point $p$}, we have:
\begin{eqnarray} \label{Phicoordinateseq}
\nonumber \Ph_p & = & dx^{0123} - dx^{0167} - dx^{0527} - dx^{0563}
+ dx^{0415} + dx^{0426} + dx^{0437} \\ & & {} + dx^{4567}
- dx^{4523} - dx^{4163} - dx^{4127} + dx^{2637} + dx^{1537} +
dx^{1526}
\end{eqnarray}
where $dx^{ijkl} = dx^i \wedge dx^j \wedge dx^k \wedge dx^l$. In
these coordinates the metric at $p$ is the standard Euclidean
metric \begin{equation*}
g_p = \sum_{k = 1}^8 dx^k \otimes dx^k
\end{equation*}
and $\st \Ph = \Ph$, so $\Ph$ is {\em self-dual}.

The $4$-forms that arise from a \SPs\ are called {\em positive} or
{\em non-degenerate}, and this set is denoted $\wfopos$. The
subgroup of $\operatorname{SO(8)}$ that preserves $\Ph_p$ is $\SP$.
(see~\cite{Br2}.) Hence at each point $p$, the set of \SPs s at
$p$ is isomorphic to $\operatorname{GL(8, \mathbb R)} / \SP$, which
is $64 - 21 = 43$ dimensional. This time, however, in contrast to
the \G\ case, since $\wedge^4 (\Rei)$ is $70$ dimensional, the set
$\wfopos (p)$ of positive $4$-forms at $p$ is {\em not} an open
subset of $\wedge^4_p$. One of the consequences of this is that
the analogous non-infinitesmal deformation in the \SP\ case will
{\em not} work. This is discussed in Section~\ref{w47deformsec}.

\subsection{Decompostion of $\bigwedge^{*}(M)$ into irreducible \SP
-representations} \label{Spin7reptheorysec}
The facts collected in this section about the decomposition of the
space of forms in the \SP\ case can also be found in~\cite{F1, J4,
Sa}.

There is an action of the group $\SP$ on $\Rei$, and hence on the
spaces $\wedge^\st$ of differential forms on $M$. We can decompose
each space $\wedge^k$ into irreducible $\SP$-representations. The
results of this decomposition are presented below. As before, the
notation $\wedge^k_l$ refers to an $l$-dimensional irreducible
$\SP$-representation which is a subspace of $\wedge^k$, $w$ is a
vector field on $M$ and $\vol$ is the volume form.

\begin{eqnarray*} 
\wzeo = \{f \in C^{\infty}(M) \} & \qquad \qquad & \wone =
\{\alpha \in \Gamma(\wedge^1 (M) \} \\ \wedge^2 = \wtws \oplus
\wtwt & \qquad \qquad & \wedge^3 = \wthe \oplus \wthfo \\
\wedge^4 = \wfoo \oplus \wfos \oplus \wfot \oplus \wfoth & \quad
\qquad & \\ \wedge^5 = \wfie \oplus \wfifo & \qquad \qquad &
\wedge^6 = \wsis \oplus \wsit \\ \wsee = \{ w \hk \vol \} & \qquad
\qquad & \weio = \{f \vol; f \in C^{\infty}(M) \}
\end{eqnarray*}

This decomposition respects the {\em Hodge star} $\st$ operator
since $\SP \subset \operatorname{SO(8)}$, so $\st \wedge^k_l =
\wedge^{8-k}_l$. Again, we will give explicit descriptions of the
remaining cases.

Before we describe $\wedge^k_l$ for $k = 2,3,4,5,6$, there are some
isomorphisms between these subspaces:
\begin{prop} \label{SPisomorphisms}
The map $\alpha \mapsto \Ph \wedge \alpha$ is an isomorphism
between the following spaces:
\begin{eqnarray*}
\wzeo \cong \wfoo & \qquad \qquad & \wone \cong \wfie \\ \wtws
\cong \wsis & \qquad \qquad & \wtwt \cong \wsit \\ \wthe \cong
\wsee & \qquad \qquad & \wfoo \cong \weio
\end{eqnarray*}
In addition, if $\alpha$ is a $1$-form, we have the following
identity:
\begin{equation} \label{isoPheq} \st \left( \Ph \wedge \st (\Ph
\wedge \alpha) \right) = - 7 \alpha 
\end{equation}
\end{prop}
\begin{proof}
These statements can be easily checked pointwise using local
coordinates and~\eqref{Phicoordinateseq}.
\end{proof}

From this identity, we can prove the following lemma:
\begin{lemma} \label{SPonewedgelemma}
If $\alpha$ is a $1$-form on $M$, then we have:
\begin{equation} \label{SPonewedgePheq} {|\Ph \wedge \alpha|}^2 =
7 {| \alpha |}^2
\end{equation}
\end{lemma}
\begin{proof}
From~\eqref{isoPheq}, we have:
\begin{eqnarray*}
\Ph \wedge \st ( \Ph \wedge \alpha ) & = & 7 \st \alpha \\ \alpha
\wedge \Ph \wedge \st ( \Ph \wedge \alpha) & = & 7 \alpha \wedge
\st \alpha \\ {| \Ph \wedge \alpha |}^2 \vol & = & 7 {|\alpha
|}^2 \vol \\ 
\end{eqnarray*}
which proves~\eqref{SPonewedgePheq}.
\end{proof}

We have some relations between $\Ph$ and an arbitrary vector field
$w$:
\begin{lemma} \label{isomorphismsv2}
The following relations hold for any vector field $w$, where $\ws$
is the associated $1$-form:
\begin{eqnarray}
\label{hkPheq} \st \left( \Ph \wedge \ws \right) & = & w \hk
\Ph \\ \label{hk7eq} \Ph \wedge \left( w \hk \Ph \right) & =
& 7 \st \ws
\end{eqnarray}
\end{lemma}
\begin{proof}
Since on an $8$-manifold $\st^2 = (-1)^k$ on $k$-forms, these
results follow from Lemma~\ref{intextidentities} and
Proposition~\ref{SPisomorphisms}.
\end{proof}

We now explicitly describe the decomposition of the space
of forms, beginning with $k=2,6$. These should be compared to the
\G\ case which were given in~\eqref{wtwsdesc}~--~\eqref{wfifdesc}.
\begin{eqnarray}
\label{wtws2desc} \wtws & = & \{\beta \in \wedge^2; \st (\Ph
\wedge \beta) = 3 \beta \} \\ \label{wtwtdesc} \wtwt & = & \{\beta
\in \wedge^2; \st (\Ph \wedge \beta) = -\beta \} \\ \nonumber & = &
\{\sum a_{ij} e^i \wedge e^j ; (a_{ij}) \in \liesp \}
\end{eqnarray}

\begin{eqnarray}
\label{wsiedesc} 
\wsis & = & \{\mu \in \wedge^6; \Ph \wedge \st \mu = 3
\mu \} \\ \label{wsitdesc} \wsit & = & \{ \mu \in \wedge^6 ; \Ph
\wedge \st \mu = -\mu \}
\end{eqnarray}

Note that these subspaces are $+3$ and $-1$ eigenspaces of the
operators $G(\beta) = \st( \Ph \wedge \beta)$ on $\wedge^2$ and
$M(\mu) = \Ph \wedge \st \mu$ on $\wedge^6$. From this fact we get
the following useful formulas for the projections $\pi_k$ onto the
$k$-dimensional representations, for $\beta \in \wedge^2$ and $\mu
\in \wedge^6$:

\begin{eqnarray} \label{SPprojectioneq} \st (\Ph \wedge \beta) & =
& 3 \pi_7(\beta) - \pi_{21}(\beta) \\
\nonumber \pi_7(\beta) & = & \frac{\beta + \st(\Ph \wedge
\beta)}{4} \\ \nonumber
\pi_{21}(\beta) & = & \frac{3 \beta - \st(\Ph \wedge \beta)}{4}
\end{eqnarray}
and
\begin{eqnarray*} \Ph \wedge \st \mu & = & 3 \pi_7(\mu) -
\pi_{21}(\mu) \\ \pi_7(\mu) & = & \frac{\mu + \Ph \wedge  \st
\mu}{4} \\ \pi_{21}(\mu) & = & \frac{3 \mu - \Ph \wedge \st
\mu}{4}
\end{eqnarray*}

We now move on to the decompositions for $k=3,4,5.$ For $k=3$, we
have:
\begin{eqnarray}
\label{wthedesc} \wthe & = & \{\st(\Ph \wedge \alpha); \alpha \in
\wone\} \\ \nonumber & = & \{ w \hk \Ph;  w \in \Gamma(T(M)) \}
\\ \nonumber & = &  \{\eta \in \wedge^3; \st \left( \Ph \wedge \st
(\Ph \wedge \eta) \right) = -7 \eta \} \\ \label{wthfdesc} \wthfo
& = & \{ \eta \in \wedge^3; \Ph \wedge \eta = 0 \}
\end{eqnarray}

For $k=5$, the decomposition is:
\begin{eqnarray}
\label{wfiedesc} \wfie & = & \{\Ph \wedge \alpha ; \alpha \in
\wone \} \\ \label{wfifodesc} \wfifo & = & \{ \mu \in \wedge^5 ;
\Ph \wedge \st \mu = 0 \}
\end{eqnarray}

And finally, the middle dimension $k=4$ decomposes as:
\begin{eqnarray}
\label{wfoodesc2} \wfoo & = & \{ f \Ph; f \in C^{\infty}(M)\} \\
\nonumber & = & \{\eta \in \wedge^4; \Ph \wedge \left( \st (\Ph
\wedge \eta) \right) = 14 \eta \} \\ \label{wfosdesc2} \wfos & = &
\{ \beta_i^{\ j} e^i \wedge ( e_j \hk \Ph ) - \beta^i_{\ j} e^j
\wedge ( e_i \hk \Ph ); \beta_{ij} e^i \wedge e^j \in \wtws \} \\
\label{wfotdesc2} \wfot & = & \{ \sigma \in
\wedge^4; \st \sigma = \sigma, \sigma \wedge \Ph = 0, \sigma \wedge
\tau = 0 \quad \forall \tau \in \wfos \} \\ \label{wfothdesc}
\wfoth & = & \{ \sigma \in \wedge^4; \st \sigma = - \sigma \}
\end{eqnarray}

\subsection{The metric of a \SPs} \label{SPmetricsec}

Here the situation differs significantly from the \G\ case.
Because $\Ph$ is self-dual equation~\eqref{isoPheq} gives us
only one useful identity rather than the four identities in
equations~\eqref{isopheq}~--~\eqref{iso2eq}. In particular it was
equation~\eqref{iso2eq} which enabled us to prove
Proposition~\ref{metricprop} to obtain a formula for the metric
from the $3$-form $\ph$ in the \G\ case.

The prescription for obtaining the metric from the $4$-form $\Ph$
in the \SP\ case is much more complicated. Before we can do this,
we need to collect some facts about various $2$-forms which
can be constructed from pairs of vector fields, as these facts
will be used both to determine the metric and later to analyze how
it changes under a $\wfos$ deformation in
Section~\ref{w47deformsec}.

\begin{prop} \label{SP2formsprop}
Let $a$, $b$, $c$, and $d$ be vector fields. Define the
$2$-forms $\beta = \as \wedge \bs = \beta_7 + \beta_{21}$ and $\mu
= \cs \wedge \ds = \mu_7 + \mu_{21}$. Then we can construct other
$2$-forms $a \hk b \hk \Ph$ and $\st \left( (a \hk \Ph)
\wedge (b \hk \Ph) \right)$ from $a$ and $b$ and these are related
to $\beta$ by
\begin{eqnarray}
\label{vhkwhkPheq} a \hk b \hk \Ph & = & -3 \beta_7 + \beta_{21} \\
\label{vhkPhwhkPheq} \st \left( (a \hk \Ph) \wedge (b \hk \Ph)
\right) & = & 2 \beta_7 - 6 \beta_{21}
\end{eqnarray}
Furthermore, if we define
\begin{eqnarray} \label{twoformsAeq}
A & = & \langle a \wedge b, c \wedge d \rangle = \langle a, c
\rangle \langle b, d \rangle - \langle a, d \rangle \langle b, c
\rangle \\ B & = & \Ph (a,b,c,d) 
\end{eqnarray}
then the following relations hold between these
$2$-forms:
\begin{eqnarray}
\label{SP2formdot1eq}  (a \hk b \hk \Ph) \wedge (\cs \wedge \ds)
\wedge \Ph & = & \left(-3 A -2 B \right) \vol \\
\label{SP2formdot2eq} (\as \wedge \bs) \wedge (c \hk \Ph)
\wedge (d \hk \Ph) & = & \left(-4 A + 2 B \right) \vol \\
\label{SP2formdot5eq} (a \hk b \hk \Ph) \wedge (c \hk d \hk \Ph)
\wedge \Ph & = & \left( 6 A + 7 B \right) \vol 
\end{eqnarray}
\end{prop}
\begin{proof}
Let $\beta = \as \wedge \bs = \beta_7 + \beta_{21}$ using the
decompositions in~\eqref{wtws2desc} and~\eqref{wtwtdesc}. From
Lemma~\ref{intextidentities} we can write
\begin{eqnarray*}
a \hk b \hk \Ph & = & \st ( \as \wedge \st (b \hk \Ph) ) \\ &
= & - \st ( \as \wedge \bs \wedge \Ph ) \\ & = & -3 \beta_7 +
\beta_{21}
\end{eqnarray*}
where we have used the self-duality $\st \Ph = \Ph$ and the
characterizations of $\wtws$ and $\wtwt$. Now since $\Ph \wedge \Ph
= 14 \vol$, we have
\begin{equation*}
(w \hk \Ph) \wedge \Ph = 7 w \hk \vol = 7 \st \ws
\end{equation*}
where we have used~\eqref{wvoleq}. Taking the interior product on
both sides with $v$, 
\begin{eqnarray*}
(v \hk w \hk \Ph) \wedge \Ph - (w \hk \Ph) \wedge ( v \hk \Ph) & =
& 7 v \hk \st \ws \\ & = & - 7 \st (\vs \wedge \ws) \\ (-3 \beta_7
+ \beta_{21}) \wedge \Ph + (v \hk \Ph) \wedge ( w \hk \Ph) & = & -
7 \st \beta_7 - 7 \st \beta_{21} \\ - 9 \st \beta_7 - \st
\beta_{21} + (v \hk \Ph) \wedge ( w \hk \Ph) & = & -
7 \st \beta_7 - 7 \st \beta_{21}
\end{eqnarray*}
which can be rearranged to give~\eqref{vhkPhwhkPheq}.
We also have
\begin{eqnarray*}
B \vol = \Ph (a,b,c,d) \vol & = & \as \wedge \bs \wedge \cs \wedge
\ds \wedge \Ph \\ & = & \left( \beta_7 + \beta_{21} \right) \wedge
\left( 3 \st \mu_7 - \st \mu_{21} \right) \\ & = & \left( 3 \langle
\beta_7, \mu_7 \rangle - \langle \beta_{21} , \mu_{21} \rangle
\right) \vol
\end{eqnarray*}
and 
\begin{equation*}
A = \langle \beta, \mu \rangle = \langle \beta_7 , \mu_7 \rangle +
\langle \beta_{21}, \mu_{21} \rangle
\end{equation*}
which together give that
\begin{equation*}
\langle \beta_7 , \mu_7 \rangle = \frac{A + B}{4} \qquad \qquad
\langle \beta_{21} , \mu_{21} \rangle = \frac{3 A - B}{4}
\end{equation*}
Hence, for example
\begin{eqnarray*}
(a \hk b \hk \Ph) \wedge (c \hk d \hk \Ph)
\wedge \Ph & = & (-3 \beta_7 + \beta_{21} ) \wedge (-9 \st \mu_7 -
\st \mu_{21} ) \\ & = &  27 \left( \frac{A + B}{4} \right) \vol -
\left( \frac{3 A - B}{4} \right) \vol \\ & = & \left( 6 A + 7 B
\right) \vol 
\end{eqnarray*}
which is~\eqref{SP2formdot5eq}. The other two are obtained
similarly.
\end{proof}

Proposition~\ref{SP2formsprop} immediately yields the following
corollary, which is analogous to Proposition~\ref{metricprop} in
the \G\ case.
\begin{cor} \label{SPmetricidentitycor}
The following identity holds for $v$ and $w$ vector fields:
\begin{equation} \label{SPmetricidentityeq}
( v \hk w \hk \Ph ) \wedge (v \hk w \hk \Ph) \wedge \Ph = 6
{|v \wedge w|}^2 \vol
\end{equation}
\end{cor}
\begin{proof}
This follows from~\eqref{SP2formdot5eq}.
\end{proof}
If we polarize~\eqref{SPmetricidentityeq} in $w$, we obtain the
useful equation:
\begin{eqnarray} \nonumber
( v \hk w_1 \hk \Ph ) \wedge (v \hk w_2 \hk \Ph) \wedge \Ph & = &
6 {\langle v \wedge w_1 , v \wedge w_2 \rangle} \vol \\ & =
\label{SPmetricidentity2eq} & 6 {\left( {|v|}^2 {\langle w_1, w_2
\rangle} - {\langle v, w_1 \rangle}{\langle v, w_2 \rangle}
\right)} \vol
\end{eqnarray}
From this equation we can obtain the metric.

\begin{lemma} \label{SPmetricprelemma}
Fix a {\em non-zero} vector field $v = v^k e_k$, and extend to
$v, e_1, e_2, \ldots, e_7$ an oriented local frame of vector
fields. The expression obtained from $v$ by
\begin{equation} \label{SPmetricprelemmaeq}
\frac{ {\left( \det{\left( \left( (e_i \hk v \hk \Ph) \wedge (e_j
\hk v \hk \Ph) \wedge (v \hk \Ph) \right) (e_1, e_2, \ldots, e_7)
\right)} \right)}^{\frac{1}{3}} }{{ \left( \left( (v \hk \Ph)
\wedge \Ph \right) (e_1, e_2, \ldots, e_7) \right)}^3}
\end{equation}
is homogeneous of order $4$ in $v$, and is independent of the
choice of extension to a basis $e_1, \ldots, e_7$. We will see
in the next theorem that up to a constant, this is ${|v|}^4$.
\end{lemma}
\begin{proof}
The expression $(e_i \hk v \hk \Ph) \wedge (e_j \hk v \hk \Ph)
\wedge (v \hk \Ph)$ is cubic in $v$, so after taking the $7 \times
7$ determinant and the cube root, the numerator
of~\eqref{SPmetricprelemmaeq} is of order $7$ in $v$. Since the
denominator is cubic in $v$, the whole expression is homogeneous
of degree $4$ in $v$. Now suppose we extend $v$ to an oriented
basis $v, e'_1, e'_2, \ldots, e'_7$ in a different way. Then
\begin{equation*}
e'_i = P_{ij} e_j + Q_i v 
\end{equation*}
and we have
\begin{equation*}
(e'_i \hk v \hk \Ph) \wedge (e'_j \hk v \hk \Ph) \wedge (v \hk
\Ph) = P_{ik} P_{jl} (e_k \hk v \hk \Ph) \wedge (e_l \hk v \hk
\Ph) \wedge (v \hk \Ph) 
\end{equation*}
Hence in the new basis the numerator of~\eqref{SPmetricprelemmaeq}
changes by a factor of
\begin{equation*}
{\left( \det(P)^2 \det(P)^7 \right)}^{\frac{1}{3}} = \det(P)^3
\end{equation*}
and the denominator also changes by a factor of $\det(P)^3$, so
the quotient is invariant.
\end{proof}
\begin{rmk}
Note how different from Lemma~\ref{Gmetricprelemmaeq} this is. In
this case, we have a determinant in the numerator rather than the
denominator. Also, the numerator and denominator are totally
different in the \SP\ case whereas one is the determinant of the
other in the \G\ case.
\end{rmk}

We now derive the expression for the metric in terms of the
$4$-form $\Ph$ in the \SP\ case.
\begin{thm} \label{SPmetrictheorem}
Let $v$ be a non-zero tangent vector at a point $p$ and let $e_0,
e_1, \ldots, e_7$ be any oriented basis for $T_p M$, so that
$\vol (e_0, e_1,\ldots, e_7) > 0 $. Assume without loss of
generality that $v^0 \neq 0$. Then the length $|v|$ of $v$ is given
by
\begin{equation*}
{|v|}^4 = \frac{(7)^3}{(6)^{\frac{7}{3}}} \frac{ {\left(
\det{\left( \left( (e_i \hk v \hk \Ph) \wedge (e_j
\hk v \hk \Ph) \wedge (v \hk \Ph) \right) (e_1, e_2, \ldots, e_7)
\right)} \right)}^{\frac{1}{3}} }{{ \left( \left( (v \hk \Ph)
\wedge \Ph \right) (e_1, e_2, \ldots, e_7) \right)}^3}
\end{equation*}
\end{thm}
\begin{proof}
We work in local coordinates at the point $p$. In this notation
$g_{ij} = {\langle e_i, e_j \rangle}$ with $0 \leq i,j \leq 7$.
Let $\det_8 (g)$ denote the $8 \times 8$ determinant of $(g_{ij})$
and let $\det_7(g)$ denote the $7 \times 7$ determinant of the
submatrix where $1 \leq i,j \leq 7$. Using the fact that $\Phi^2 =
14 \vol = 14 \sqrt{\det_8 (g)} e^0 \wedge e^1 \ldots \wedge e^7$,
and writing $v = v^k e_k$, we compute
\begin{eqnarray} \nonumber
A(v) & = & \left( (v \hk \Ph) \wedge \Ph \right) (e_1, e_2, \ldots,
e_7) \\ \label{SPmetricAeq} & = & 7 v^0 \sqrt{{\det}_8(g)}
\end{eqnarray}
Now $\langle v, e_j \rangle = v^k g_{kj} = v_j$. We also have the
$7 \times 7$ matrix (for $1 \leq i,j \leq 7$)
\begin{eqnarray} \nonumber
B_{ij}(v) & = & \left( (e_i \hk v \hk \Ph) \wedge (e_j \hk v \hk
\Ph) \wedge (v \hk \Ph) \right) (e_1, e_2, \ldots, e_7) \\ & = &
\label{SPmetricBeq} 6 \left( {|v|}^2 g_{ij} - v_i v_j \right) v^0
\sqrt{{\det}_8 (g)}
\end{eqnarray}
where we have used Corollary~\ref{SPmetricidentitycor}. Now
consider the $7 \times 7$ matrix $( {|v|}^2 g_{ij} - v_i v_j)$. By
examining the proof of Lemma~\ref{specialdetlemma}, we see that its
determinant is
\begin{equation} \label{SPmetrictempdeteq}
{|v|}^{14} {\det}_7{ \left( g \right) } - {|v|}^{12} {| \vs
\wedge e^0 |}^2 {\det}_8 (g)
\end{equation}
Now from Cramer's rule ${\det}_7 (g) = g^{00} {\det}_8 (g)$
and we also have ${| \vs \wedge e^0 |}^2 = {|v|}^2 g^{00} - v^0
v^0$. Hence~\eqref{SPmetrictempdeteq} becomes
\begin{equation*}
{|v|}^{12} v^0 v^0 {\det}_8 (g)
\end{equation*}
Returning to~\eqref{SPmetricBeq}, we have now shown that
\begin{eqnarray*} 
\det { B_{ij}(v) } & = & 6^7 {|v|}^{12} {(v^0)}^2 {\det}_8
(g) {(v^0)}^7 {({\det}_8 (g))}^{\frac{7}{2}} \\ & = & 6^7
{|v|}^{12} {(v^0)}^9 {({\det}_8 (g))}^{\frac{9}{2}}
\end{eqnarray*}
and hence
\begin{equation*}
{(\det { B_{ij}(v) })}^{\frac{1}{3}} = 6^{\frac{7}{3}} {|v|}^4 
{(v^0)}^3 {({\det}_8 (g))}^{\frac{3}{2}}.
\end{equation*}
Finally, since from~\eqref{SPmetricAeq} we have
\begin{equation*}
{(A(v))}^3 = (7)^3 {(v^0)}^3 {({\det}_8 (g))}^{\frac{3}{2}}
\end{equation*}
these two expressions can be combined to yield
\begin{equation} \label{SPmetricproofeq}
{|v|}^4 = \frac{(7)^3}{(6)^{\frac{7}{3}}} \frac{{(\det { B_{ij}(v)
})}^{\frac{1}{3}}}{{(A(v))}^3}
\end{equation}
which completes the proof.
\end{proof}

\subsection{The triple cross product of a \SPs}

In this section we will describe the triple cross product operation
on a manifold with a \SPs\ in terms of the $4$-form $\Ph$, and
present some useful relations.

\begin{defn} \label{SPcrossproductdefn}
Let $u$, $v$, and $w$ be vector fields on $M$. The {\em triple
cross product}, denoted $X( u, v, w)$, is a vector field on $M$
whose associated $1$-form under the metric isomorphism satisfies:
\begin{equation} \label{SPcrossproducteq}
{\left( X(u, v, w) \right)}^{\flat} = w \hk v \hk u \hk \Ph
\end{equation}
This immediately yields the relation between $X$, $\Ph$, and the
metric $g$:
\begin{equation} \label{gPhirelationeq}
g ( X(u, v, w) , y) = {\left( X(u, v, w) \right)}^{\flat} (y) = y
\hk w \hk v \hk u \hk \Ph = \Ph (u,v,w,y).
\end{equation}
\end{defn}

We can obtain another useful characterization of the triple cross
product from this one using Lemma~\ref{intextidentities}:
\begin{eqnarray} \label{SPcrossproductchareq}
{\left( X(u, v, w) \right)}^{\flat} & = & w \hk v \hk u \hk \Ph \\
\nonumber & = & \st ( \ws \wedge \st ( v \hk u \hk \Ph ) ) \\
\nonumber & = & \st ( \ws \wedge \vs \wedge \st ( u \hk \Ph) ) \\
\nonumber & = & - \st ( \ws \wedge \vs \wedge \us \wedge \st \Ph )
\\ \nonumber & = & \st( \us \wedge \vs \wedge \ws \wedge \Ph ) 
\end{eqnarray}
Note the similarity to the \G\ case given
by~\eqref{Gcrossproductchareq}.

Since $\us \wedge \vs \wedge \ws$ is a $3$-form, we can write
it as $\gamma_8 + \gamma_{48}$, with $\gamma_j \in \wedge^3_j$.
Using~\eqref{wthedesc} and~\eqref{wthfdesc} we see that: 
\begin{eqnarray} \label{SPgamma7eq}
{\left( X(u, v, w) \right)}^{\flat} \wedge \Ph & = & \st ( \gamma_7
\wedge \Ph ) \wedge \Ph \\ \nonumber & = & 7 \st \gamma_7
\end{eqnarray}
Taking the norm of both sides, and using~\eqref{SPonewedgePheq}:
\begin{equation*}
{|{\left( X(u, v, w) \right)}^{\flat} \wedge \Ph |}^2 = 7 {|
{\left( X(u, v, w) \right)}^{\flat} |}^2 = 7 {|  X(u, v, w) |}^2 =
49 {|
\gamma_7 |}^2
\end{equation*}
from which we obtain
\begin{equation} \label{SPgamma7normeq}
{| \gamma_7 |}^2 = \frac{1}{7} {| X(u, v, w) |}^2
\end{equation}
We can now establish the following Lemma.
\begin{lemma} \label{SPcrossproductnormlemma}
Let $u$, $v$, and $w$ be vector fields. Then
\begin{equation} \label{SPcrossproductnormeq}
{| X(u, v, w) |}^2 = { | u \wedge v \wedge w |}^2
\end{equation}
\end{lemma}
\begin{proof}
First we note that
\begin{eqnarray*}
{ | u \wedge v \wedge w |}^2 & = & \det {\begin{pmatrix} {|u|}^2 &
{\langle u, v \rangle} & {\langle u, w \rangle} \\ {\langle u, v
\rangle} & {|v|}^2 & {\langle v, w \rangle} \\ {\langle u, w
\rangle} & {\langle v, w \rangle} & {|w|}^2
\end{pmatrix}} \\ & = & {|u|}^2 {|v|}^2 {|w|}^2 + 2 {\langle u, v
\rangle} {\langle v, w \rangle} {\langle u, w \rangle} \\ & & {} -
{|u|}^2 {\langle v, w \rangle}^2 - {|v|}^2 {\langle u, w
\rangle}^2 - {|w|}^2 {\langle u, v \rangle}^2
\end{eqnarray*}
as we will have to identify an expression of this form in what
follows. Starting from $\Ph^2 = 14 \vol$ and taking the interior
product with $u$, $v$, and $w$, we obtain
\begin{eqnarray} \label{cross8usefulinteq}
\left( w \hk v \hk u \hk \Ph \right) \wedge \Ph & - & \left( u \hk
v \hk \Ph \right) \wedge \left( w \hk \Ph \right) \\ \nonumber & {}
- & \left( w \hk u \hk \Ph \right) \wedge \left( v \hk \Ph \right)
\\ \nonumber & {} - &  \left( v \hk w \hk \Ph \right) \wedge \left(
u \hk \Ph \right) = 7 \st \left( \us \wedge \vs \wedge \ws \right)
\end{eqnarray}
where we have interchanged interior and exterior multiplication
through the $\st$ operator using the identities in
Appendix~\ref{identitiessec}. Now from~\eqref{SPcrossproducteq}
and~\eqref{SPcrossproductchareq}, if we take the wedge product of
the above expression with $\us \wedge \vs \wedge \ws$, we have
\begin{eqnarray*}
{| X(u,v,w) |}^2 \vol & - & \us \wedge \vs \wedge \ws \wedge \left(
u \hk v \hk \Ph \right) \wedge \left( w \hk \Ph \right) \\ & {} -
& \us \wedge \vs \wedge \ws \wedge \left( w \hk u \hk \Ph \right)
\wedge \left( v \hk \Ph \right) \\ & {} - & \us \wedge \vs
\wedge \ws \wedge \left( v \hk w \hk \Ph \right) \wedge \left( u
\hk \Ph \right) = 7 {| u \wedge v \wedge w |}^2 \vol
\end{eqnarray*}
Consider the second term:
\begin{equation*}
\us \wedge \vs \wedge \ws \wedge \left( u \hk v \hk \Ph
\right) \wedge \left( w \hk \Ph \right)
\end{equation*}
This can be rewritten as
\begin{eqnarray*}
& & \Ph \wedge \left( w \hk \left( \us \wedge \vs \wedge \ws
\wedge \left( u \hk v \hk \Ph \right) \right) \right) \\ & = & 
\langle w, u \rangle \Ph \wedge \vs \wedge \ws \wedge \left (u
\hk v \hk \Ph \right)  - \langle w, v \rangle \Ph \wedge \us \wedge
\ws \wedge \left (u \hk v \hk \Ph \right) \\ & & {} + {|w|}^2
\Ph \wedge \us \wedge \vs \wedge \left (u \hk v \hk \Ph \right) -
\Ph \wedge \us \wedge \vs \wedge \ws \wedge \left( w \hk u
\hk v \hk \Ph \right)
\end{eqnarray*}
Using~\eqref{SP2formdot1eq}, this is equal to
\begin{eqnarray*}
& & \left( 3 \langle w , u \rangle \langle v \wedge w, v \wedge u
\rangle + 3 \langle w , v \rangle \langle u \wedge w, u \wedge v
\rangle - 3 {|w|}^2 {|u \wedge v|}^2 \right) \vol + {|X(u,v,w)|}^2
\vol \\ & & = - 3 {|u \wedge v \wedge w|}^2 \vol + {|X(u,v,w)|}^2
\vol
\end{eqnarray*}
The other two terms are identical by cyclically permuting $u$,
$v$, and $w$. Thus we have
\begin{eqnarray*}
{|X(u,v,w)|}^2 \vol + 9 {|u \wedge v \wedge w|}^2 \vol  - 3
{|X(u,v,w)|}^2 \vol & = & 7 {|u \wedge v \wedge w|}^2 \\ 
{|X(u,v,w)|}^2 & = & {|u \wedge v \wedge w|}^2
\end{eqnarray*}
and the lemma is proved.
\end{proof}

Iterating the $3$-fold cross product is described in the next
lemma.
\begin{lemma} \label{cross8iteratelemma}
Let $a$, $b$, $u$, $v$, and $w$ be vector fields. Then we have
\begin{eqnarray*}
& & X\left( a, b, X( u,v,w) \right) = \\ & &
-\langle a \wedge b, u \wedge v \rangle w - \Ph(a,b,u,v)
w + \langle b, w \rangle X(a,u,v) - \langle a, w \rangle
X(b,u,v) \\ & & -\langle a \wedge b, v \wedge w \rangle u
- \Ph(a,b,v,w) u + \langle b, u \rangle X(a,v,w) -
\langle a, u \rangle X(b,v,w) \\ & & -\langle a \wedge b,
w \wedge u \rangle v - \Ph(a,b,w,u) v + \langle b, v
\rangle X(a,w,u) - \langle a, v \rangle X(b,w,u)
\end{eqnarray*}
\end{lemma}
\begin{proof}
We begin with~\eqref{cross8usefulinteq},
\begin{eqnarray*}
{X(u,v,w)}^{\flat} \wedge \Ph = \left( w \hk v \hk u \hk \Ph
\right) \wedge \Ph & = & \left( u \hk v \hk \Ph \right) \wedge
\left( w \hk \Ph \right) \\ & {} + & \left( w \hk u \hk \Ph
\right) \wedge \left( v \hk \Ph \right) \\ & {} + &  \left( v \hk w
\hk \Ph \right) \wedge \left( u \hk \Ph \right) + 7 \st \left( \us
\wedge \vs \wedge \ws \right)
\end{eqnarray*}
We take the wedge product of both sides of this equation with $\as
\wedge \bs$, and obtain
\begin{eqnarray} \label{cross8iteratet2eq}
\as \wedge \bs \wedge {X(u,v,w)}^{\flat} \wedge \Ph & = & \as
\wedge \bs \wedge \left( u \hk v \hk \Ph \right) \wedge
\left( w \hk \Ph \right) \\ \nonumber & & {}+ \as \wedge \bs \wedge
\left( v \hk w \hk \Ph \right) \wedge \left( u \hk \Ph \right) \\
\nonumber & & {}+ \as \wedge \bs \wedge \left( w \hk u \hk \Ph
\right) \wedge \left( v \hk \Ph \right) \\ \nonumber & & {}+ 7 \as
\wedge \bs \wedge \st \left( \us \wedge \vs \wedge \ws \right)
\end{eqnarray}
Note that up to the metric isomorphism, the desired iterated cross
product we are looking for is the Hodge star of the above
expression. It is easy to check that the last term can be
simplified to
\begin{equation*}
7 \left( \langle a \wedge b, v \wedge w \rangle \st \us + \langle
a \wedge b, w \wedge u \rangle \st \vs + \langle a \wedge b, u
\wedge v \rangle \st \ws  \right)
\end{equation*}
The first three terms are identical under cyclic permutation of
$u,v,w$. Let us consider the first term.
From~\eqref{SP2formdot1eq} we have
\begin{equation*}
\as \wedge \bs \wedge \left( u \hk v \hk \Ph \right) \wedge \Ph =
-3 \langle a \wedge b, u \wedge v \rangle \vol - 2 \Ph(a,b,u,v)
\vol
\end{equation*}
Taking the interior product with $w$ and rearranging,
\begin{eqnarray*}
\as \wedge \bs \wedge \left( u \hk v \hk \Ph \right) \wedge
\left( w \hk \Ph \right) & = & -3 \langle a \wedge b, u \wedge v
\rangle \st \ws - 2 \Ph(a,b,u,v) \st \ws \\ & & {}+ \langle w, b
\rangle \as \wedge \left( u \hk v \hk \Ph \right) \wedge \Ph -
\langle w, a \rangle \bs \wedge \left( u \hk v \hk \Ph \right)
\wedge \Ph \\ & & {}+ \as \wedge \bs \wedge {X(u,v,w)}^{\flat}
\wedge \Ph
\end{eqnarray*}
Adding up the four expressions in~\eqref{cross8iteratet2eq}, we
have
\begin{eqnarray*} 
-2 \as \wedge \bs \wedge {X(u,v,w)}^{\flat} \wedge \Ph & = & 4
\langle a \wedge b, u \wedge v \rangle \st \ws - 2 \Ph(a,b,u,v) \st
\ws \\ & & {}+ 4 \langle a \wedge b, v \wedge w \rangle \st \us - 2
\Ph(a,b,v,w) \st \us \\ & & {}+ 4 \langle a \wedge b, w \wedge u
\rangle \st \vs - 2 \Ph(a,b,w,u) \st \vs \\ & & {}+ \langle w, b
\rangle \as \wedge \left( u \hk v \hk \Ph \right) \wedge \Ph -
\langle w, a \rangle \bs \wedge \left( u \hk v \hk \Ph \right)
\wedge \Ph \\ & & {}+ \langle u, b \rangle \as \wedge \left( v \hk
w \hk \Ph \right) \wedge \Ph - \langle u, a \rangle \bs \wedge
\left( v \hk w \hk \Ph \right) \wedge \Ph \\ & & {}+ \langle v, b
\rangle \as \wedge \left( w \hk u \hk \Ph \right) \wedge \Ph -
\langle v, a \rangle \bs \wedge \left( w \hk u \hk \Ph \right)
\wedge \Ph
\end{eqnarray*}
Now consider the first of the last six terms:
\begin{equation*}
\langle w, b \rangle \as \wedge \left( u \hk v \hk \Ph \right)
\wedge \Ph
\end{equation*}
with $\beta = \beta_7 + \beta_{21} = \us \wedge \vs$, exactly as
in the derivation of~\eqref{SP2formdot1eq}, one can show easily
that
\begin{equation*}
\bs \wedge \left( u \hk v \hk \Ph \right) \wedge \Ph = -3 A - 2 B
\end{equation*}
where $A = \bs \wedge \st ( \us \wedge \vs ) = \langle b, v
\rangle \st \us - \langle b, u \rangle \st \vs$ and $B = \bs
\wedge \us \wedge \vs \wedge \Ph = - \st {X(b,u,v)}^{\flat}$.
Substituting this expression in the last six terms above,
simplifying, and taking the Hodge star of both sides yields the
desired result.
\end{proof}
\begin{cor} \label{cross8iteratespecial1cor}
For the special case $v = b$, we obtain the relation
\begin{eqnarray*}
& & X\left( a, b, X( u,b,w) \right) = \\ & &
-\langle a \wedge b, u \wedge b \rangle w - \langle b, w \rangle
X(a,b,u) + \langle a \wedge b, w \wedge b \rangle u + \langle
b, u \rangle X(a,b,w) \\ & & -\langle a \wedge b, w \wedge u
\rangle b - \Ph(a,b,w,u) b + {|b|}^2 X(a,w,u) - \langle a, b
\rangle X(b,w,u)
\end{eqnarray*}
This expression will be used in Section~\ref{w47deformsec} to
study the $\wfos$ deformation. Finally, if we also have
$u = a$ we have the special case
\begin{equation*}
X\left( a, b, X( a,b,w) \right) = 
-{|a \wedge b|}^2 w + \langle a \wedge b, w \wedge b \rangle a
+ \langle a \wedge b, a \wedge w \rangle b
\end{equation*}
which is analogous to~\eqref{cross7iteratespecialeq} in the
\G\ case.
\end{cor}

Finally, the corollary to the next theorem is analogous to
Theorem~\ref{w223zerothm} and we will use it in 
Section~\ref{w47deformsec} to show that to first order, deforming a
\SPs\ by an element of $\wfos$ does not change the metric. Note
first that by~\eqref{wfosdesc}, one can check (by computing
explicitly) that if we start with two vector fields $v$ and $w$,
take the projection $\pi_7 (\vs \wedge \ws)$, and transfer this to
$\wfos$, we get the $4$-form $\vs \wedge (w \hk \Ph) - \ws
\wedge (v \hk \Ph)$ in $\wfos$. Of course, this is not the most
general element of $\wfos$, but we will work with elements of
this form. 
\begin{thm} \label{w224zerothm}
Let $a$, $b$, $c$, and $d$ be vector fields. Let $\sigma \in
\wfos$ be given by $\sigma = \vs \wedge (w \hk \Ph) - \ws
\wedge (v \hk \Ph)$ for two other vector fields $v$ and $w$. Then
\begin{eqnarray*}
\left(a \hk b \hk \Ph \right) \wedge \left( c \hk d \hk \Ph
\right) \wedge \sigma & = & \left( 3 \langle w, c \rangle
\Ph(a,b,v,d) - 3 \langle w, d \rangle \Ph(a,b,v,c) \right. \\ & &
{}+ 3 \langle w, a \rangle \Ph(c,d,v,b) - 3 \langle w, b \rangle
\Ph(c,d,v,a) \\ & & {}- 3 \langle v, c \rangle \Ph(a,b,w,d) + 3
\langle v, d \rangle \Ph(a,b,w,c) \\ & & \left. {}- 3 \langle v, a
\rangle \Ph(c,d,w,b) + 3 \langle v, b \rangle \Ph(c,d,w,a) \right)
\vol
\end{eqnarray*}
\end{thm}
\begin{proof}
We start with the $9$-form
\begin{equation*}
(a \hk b \hk \Ph) \wedge (c \hk d \hk \Ph) \wedge \ws \wedge \Ph =
0
\end{equation*}
and take the interior product with $v$:
\begin{eqnarray*}
(a \hk b \hk \Ph) \wedge (c \hk d \hk \Ph) \wedge v \hk (\ws
\wedge \Ph) & = & - (v \hk a \hk b \hk \Ph) \wedge (c \hk d \hk
\Ph) \wedge \ws \wedge \Ph \\ & & {}- (a \hk b \hk \Ph) \wedge ( v
\hk c \hk d \hk \Ph) \wedge \ws \wedge \Ph
\end{eqnarray*}
Now interchanging the roles of $v$ and $w$ and subtracting, and
using~\eqref{SPcrossproducteq}, we obtain
\begin{eqnarray*}
(a \hk b \hk \Ph) \wedge (c \hk d \hk \Ph) \wedge \sigma & = &
{}- \ws \wedge {X(a,b,v)}^{\flat} \wedge (c \hk d \hk \Ph) \wedge
\Ph \\ & & {}- \ws \wedge {X(c,d,v)}^{\flat} \wedge (a \hk b \hk
\Ph) \wedge \Ph \\ & & {}+ \vs \wedge {X(a,b,w)}^{\flat} \wedge (c
\hk d \hk \Ph) \wedge \Ph \\ & & {}+ \vs \wedge {X(c,d,w)}^{\flat}
\wedge (a \hk b \hk \Ph) \wedge \Ph \end{eqnarray*}
Now from~\eqref{SP2formdot1eq}, this becomes
\begin{eqnarray*}
(a \hk b \hk \Ph) \wedge (c \hk d \hk \Ph) \wedge \sigma & = &
\left( 3 \langle w \wedge {X(a,b,v)}^{\flat}, c \wedge d \rangle +2
\Ph( w, {X(a,b,v)}^{\flat}, c, d) \right) \vol \\ & & {}+ \left( 3 
\langle w \wedge {X(c,d,v)}^{\flat}, a \wedge b \rangle +2
\Ph( w, {X(c,d,v)}^{\flat}, a, b) \right) \vol \\ & & {}- \left( 3
\langle v \wedge {X(a,b,w)}^{\flat}, c \wedge d \rangle +2
\Ph( v, {X(a,b,w)}^{\flat}, c, d) \right) \vol \\ & & {}- \left( 3
\langle v \wedge {X(c,d,v)}^{\flat}, a \wedge b \rangle +2
\Ph( v, {X(c,d,w)}^{\flat}, a, b) \right) \vol 
\end{eqnarray*}
The right hand side can be further simplified to (up to a factor
of $\vol$)
\begin{eqnarray*}
& & 3 \langle w, c \rangle \Ph(a,b,v,d) - 3 \langle w, d \rangle
\Ph(a,b,v,c) + 2 \langle {X(a,b,v)}^{\flat}, {X(c,d,w)}^{\flat}
\rangle \\ & & {}+ 3 \langle w, a \rangle \Ph(c,d,v,b) - 3 \langle
w, b \rangle \Ph(c,d,v,a) + 2 \langle {X(c,d,v)}^{\flat},
{X(a,b,w)}^{\flat} \rangle \\ & & {}- 3 \langle v, c \rangle
\Ph(a,b,w,d) + 3 \langle v, d \rangle \Ph(a,b,w,c) - 2 \langle
{X(a,b,w)}^{\flat}, {X(c,d,v)}^{\flat} \rangle \\ & & {}- 3
\langle v, a \rangle \Ph(c,d,w,b) + 3 \langle v, b \rangle
\Ph(c,d,w,a) - 2 \langle {X(c,d,w)}^{\flat}, {X(a,b,v)}^{\flat}
\rangle   
\end{eqnarray*}
The last terms in each line cancel and the theorem is proved.
\end{proof}
\begin{cor} \label{w224zerocor}
Let $h$, $u_1$, and $u_2$ be vector fields. Let $\sigma \in
\wfos$ be given by $\sigma = \vs \wedge (w \hk \Ph) - \ws
\wedge (v \hk \Ph)$ for two other vector fields $v$ and $w$. Then
\begin{equation*}
\left(h \hk u_1 \hk \Ph \right) \wedge \left( h \hk u_2 \hk \Ph
\right) \wedge \sigma = 0
\end{equation*}
\end{cor}
\begin{proof}
Letting $a=c$ in Theorem~\ref{w224zerothm} yields this result.
\end{proof}
\begin{rmk}
We can actually show the stronger result that in terms of
the decompositions in~\eqref{wtws2desc} and~\eqref{wfosdesc2}, the
wedge product map
\begin{equation*}
\wtws \times \wtws \times \wfos \to \weio
\end{equation*}
is the zero map. This is a direct analogy with
Theorem~\ref{w223zerothm}. However, we will not have occasion to
use this fact.
\end{rmk}

\subsection{The 4 classes of \SPs s}

Similar to the classification of \Gs s by Fern\'andez and Gray
in~\cite{FG}, Fern\'andez studied \SPs s in~\cite{F1}. In this
case, the results are slightly different because a $4$-form $\Ph$
which determines a \SPs\ is self-dual. Such a manifold has
holonomy a subgroup of $\SP$ if and only if $\nabla \Ph = 0$,
which Fern\'andez showed to be equivalent to
\begin{equation*}
d \Ph = 0.
\end{equation*}
Again this equivalence was established by decomposing the space $W$
that $\nabla \Ph$ belongs to into irreducible
$\SP$-representations, and comparing the invariant subspaces of
$W$ to the isomorphic spaces in $\wedge^* (M)$. In the \SP\ case
this space $W$ decomposes as
\begin{equation*}
W = W_8 \oplus W_{48}
\end{equation*}
where again the subscript $k$ denotes the dimension of the
irreducible representation $W_k$. Again in analogy with the \G\
case, we have a canonically defined $7$-form
$\zeta$ and $1$-form $\theta$, given by
\begin{eqnarray} \label{definezetaeq} \zeta & = & \st d \Ph \wedge
\Ph \\ \label{SPoneformeq} \theta & = & \st \zeta = \st \left(  \st
d \Ph \wedge \Ph \right)
\end{eqnarray}
Note that $\theta = 0$ when the manifold has holonomy contained in
\SP, and more generally $\theta$ vanishes if $\pi_8 (d \Phi) =
0$. We will see below that in the case $\pi_{48} (d \Phi) = 0$
the form $\theta$ is closed.

This time we have only $4$ classes of \SPs s: the classes $\{0\}$,
$W_8$, $W_{48}$, and $W = W_8 \oplus W_{48}$. Table~\ref{SPtable} 
describes the classes in terms of differential equations on the
form $\Ph$. Unlike the \G\ case, the inclusions between these
classes are all strict, and this is discussed in~\cite{F1}.
\begin{table}[h]
\begin{center}
\begin{tabular}{|c|c|l|c|}
\hline Class & Defining Equations & Name & $d \theta$ \\
\hline $W_8 \oplus W_{48}$  & no relation on $d \Ph$. & &
\\ \hline $W_8$ &  $d \Ph + \frac{1}{7} \theta \wedge \Ph = 0$ & LC
\SP & $d \theta = 0$ \\ \hline $W_{48}$ & $\theta = 0$ & &
$\theta = 0$ \\ \hline $\{0\}$ & $d\Ph = 0$ & \SP & $\theta = 0$
\\
\hline
\end{tabular}
\end{center}
\caption{The $4$ classes of \SPs s \label{SPtable}}
\end{table}

\begin{rmk}
Note that in the \SP\ case, there is no analogue of
an ``integrable'' structure, nor are there analogues of {\em
almost} or {\em nearly} \SPs\ as there are in the \G\ case. An
almost \SP\ manifold ($d\Ph = 0$) automatically has holonomy \SP.
And $d \Ph$ does not have a one-dimensional component which would
give us the analogue of a nearly \Gs.
\end{rmk}

We now prove the closedness of $\theta$ in the class $W_8$ as given
in the final column of Table~\ref{SPtable}.
\begin{lemma} \label{SPdthetalemma}
If $\Ph$ satisfies $d \Ph + \frac{1}{7} \theta \wedge \Ph
= 0$, then $d \theta = 0$.
\end{lemma}
\begin{proof}
Suppose $d \Ph + \frac{1}{7} \theta \wedge \Ph
= 0$. We differentiate this equation to obtain:
\begin{equation*}
d \theta \wedge \Ph = \theta \wedge d \Ph = \theta \wedge
\left( -\frac{1}{7} \theta \wedge \Ph \right) = 0
\end{equation*}
But wedge product with $\Ph$ is an isomorphism from $\wedge^2$ to
$\wedge^6$, so $d \theta = 0$.
\end{proof}

\section{Deformations of a fixed \SPs} \label{SPdeformsec}

We begin with a fixed \SPs\ on a manifold $M$ in a certain
class. We will deform the form $\Ph$ and see how this affects
the class. This time there are only $4$ classes, and only two
intermediate classes. However, the ways we can deform $\Ph$ in the
\SP\ case are more complicated. Since $\Ph \in \wfoo \oplus \wfos
\oplus \wfot \oplus \wfoth$, there are now four canonical ways to
deform the $4$-form $\Ph$. Again, since $\wfoo = \{ f \Ph \}$,
adding to $\Ph$ an element of $\wfoo$ amounts to conformally
scaling $\Ph$. This preserves the decomposition into irreducible
representations. In all other case, however, since the
decomposition does depend on $\Ph$ the decomposition will change
for the other kinds of deformations

We will see that analogously to the \G\ case, flowing in the
$\wfos$ direction gives us a path in the space of positive
$4$-forms, all corresponding to the same metric. However,
this time simply deforming non-infinitesmally by an element
of $\wfos$ will {\em not} yield a positive $4$-form, in fact
we can show that it {\em never} does. We will explain how much of
the construction does carry over and some reasons why it should
not be a surprise that discovering an analagous construction in
the \SP\ case that works should be considerably more complicated.

\subsection{Conformal Deformations of \SPs
s} \label{SPconformaldeformsec}

Let $f$ be a smooth, nowhere vanishing function on $M$. We
conformally scale $\Ph$ by $f^4$, for notational convenience.
Denote the new form by $\nPh = f^4 \oPh$. We first compute the new
metric $\tilde g$ and the new volume form $\nvol$ in the following
lemma.

\begin{lemma} \label{SPconformaldeformnewmetric}
The metric $\og$ on vector fields, the metric $\og^{-1}$ on one
forms, and the volume form $\ovol$ transform as follows:
\begin{eqnarray*}
 \tilde g & = & f^2 \og \\ \tilde g^{-1}
& = & f^{-2} \og^{-1} \\ \nvol & = & f^8 \ovol
\end{eqnarray*}
\end{lemma}
\begin{proof}
We substitute $\nPh = f^4 \oPh$ into equations~\eqref{SPmetricAeq}
and~\eqref{SPmetricBeq} to obtain
\begin{eqnarray*}
\tilde A (v) & = & \left( (v \hk f^4 \oPh) \wedge f^4 \oPh \right)
(e_1, e_2, \ldots, e_7) \\ & = & f^8 A_{\text{o}} (v) \\ \tilde
B_{ij}(v) & = & \left( (e_i \hk v \hk f^4 \oPh) \wedge (e_j \hk v
\hk f^4 \oPh) \wedge (v \hk f^4 \oPh) \right) (e_1, e_2, \ldots,
e_7) \\ & = & f^{12} (B_{\text{o}})_{ij}(v)
\end{eqnarray*}
Substituting these expressions into~\eqref{SPmetricproofeq} we
compute
\begin{eqnarray*}
\lt{|v|}^4 & = & \frac{(7)^3}{(6)^{\frac{7}{3}}} \frac{
(f^{12})^{\frac{7}{3}}{(\det {(B_{\text{o}})_{ij}
(v)})}^{\frac{1}{3}}} { {(f^8)}^3 {(A_{\text{o}} (v))}^3} \\
& = & f^4 {|v|}^4_{\text{o}}
\end{eqnarray*}
from which we have $\lt{|v|}^2 = f^2 {|v|}^2_{\text{o}}$ and the
remaining conclusions now follow.
\end{proof}

We now determine the new Hodge star $\nst$ in terms of the old
$\ost$.
\begin{lemma} \label{SPconformaldeformnewstar}
If $\alpha$ is a $k$-form, then $\nst \alpha = f^{8 - 2 k} \ost
\alpha$.
\end{lemma}
\begin{proof}
Let $\alpha$, $\beta$ be $k$-forms. Then from
Lemma~\ref{SPconformaldeformnewmetric} the new metric on $k$-forms
is $\ltn{ \left \langle \ ,\  \right \rangle } = f^{-2 k} \left
\langle \ ,\  \right \rangle _{\! \text{o}}$. From this we compute:
\begin{eqnarray*}
\beta \wedge \nst \alpha & = & \ltn{ \left \langle \beta, \alpha
\right \rangle } \nvol \\ & = & f^{-2k} \left \langle \beta,
\alpha \right \rangle _{\! \text{o}} f^8 \ovol \\ & = & f^{8 - 2k}
\beta
\wedge \ost \alpha.
\end{eqnarray*}
\end{proof}

From this we obtain the following:
\begin{lemma} \label{SPconformaldeformprelemma}
The exterior derivative of the new $4$-form $d\nPh$ and $\nst d
\nPh$ are
\begin{eqnarray*}
d\nPh & = & 4 f^3 d f \wedge \oPh + f^4 d \oPh \\ \nst d \nPh & =
& 4 f \ost (df \wedge \oPh) + f^2 \ost d \oPh
\end{eqnarray*}
\end{lemma}
\begin{proof}
This is immediate from $\nPh = f^4 \oPh$ and
Lemma~\ref{SPconformaldeformnewstar}.
\end{proof}

Using these results, we can determine which classes of \SPs s are
conformally invariant. We can also determine what happens to the
$7$-form $\zeta$ and the associated $1$-form $\theta = \st \zeta$.
This is all given in the following theorem:
\begin{thm} \label{SPconformaldeformresults}
Under the conformal deformation $\nPh = f^4 \oPh$, we have:
\begin{eqnarray}
\label{conformaldeformSPclasseq} \nSPclass & = & f^4
\left( \oSPclass \right) \\ \label{SPconformaldeformzetaeq} \nzeta
& = & -28 f^5 \ost df + f^6 \ozeta \\
\label{SPconformaldeformthetaeq}
\ntheta & = & -28 d(\log(f)) + \otheta
\end{eqnarray}
Hence, we see from Table~\ref{SPtable} and
equation~\eqref{conformaldeformSPclasseq} that only the class
$W_8$ is preserved under a conformal deformation of $\Ph$.
(This part was originally proved in~\cite{F1} using a
different method.) Also, \eqref{SPconformaldeformthetaeq} shows
that $\theta$ changes by an exact form, so in the class $W_8$,
where $\theta$ is closed, we have a well defined cohomology
class $[\theta ]$ which is unchanged under a conformal scaling.
\end{thm}
\begin{proof}
We begin by using Lemma~\ref{SPconformaldeformprelemma}
and~\eqref{definezetaeq} to compute $\nzeta$ and $\ntheta$:
\begin{eqnarray*}
\nzeta & = & \nst d \nPh \wedge \nPh \\ & = & \left( 4 f
\ost( d f \wedge \oPh) + f^2 \ost d \oPh \right) \wedge
f^4 \oPh \\ & = & 4 f^5 \oPh \wedge \ost \left( \oPh \wedge
d f \right) + f^6 \ozeta \\ & = & - 28 f^5 \ost d f + f^6 \ozeta
\end{eqnarray*}
where we have used~\eqref{isoPheq} in the last step. Now
from Lemma~\ref{SPconformaldeformnewstar}, we get:
\begin{equation*}
\ntheta = \nst \nzeta = -28 f^{-1} df + \otheta =  -28 d(\log(f)) +
\otheta.
\end{equation*}
Now using the above expression for $\ntheta$, we have:
\begin{eqnarray*}
\nSPclass & = &  4 f^3 d f \wedge \oPh + f^4 d \oPh + \frac{1}{7}
\left( -28 f^{-1} df + \otheta \right) \wedge f^4 \oPh \\ & = &
f^4 \left( \oSPclass \right)
\end{eqnarray*}
which completes the proof.
\end{proof}

The next result gives necessary and sufficient conditions for being
able to achieve holonomy \SP\ by conformally scaling.
\begin{thm} \label{SPconformaldeformfinal}
Let $\oPh$ be a positive $4$-form (associated to a \SPs). Under the
conformal deformation $\nPh = f^4 \oPh$, the new $4$-form $\nPh$
satisfies $d \nPh = 0$ if and only if $\oPh$ is already at least
class $W_8$ and $28 d \log(f) = \otheta$. Hence in order to have
$\nPh$ be closed (and hence correspond to holonomy \SP), the
original $1$-form $\otheta$ has to be {\em exact}. In particular if
the manifold is simply-connected or more generally $H^1(M) = 0$
then this will always be the case if $\oPh$ is in the class $W_8$,
since $d \otheta = 0$. 
\end{thm}
\begin{proof}
From Lemma~\ref{SPconformaldeformprelemma}, for $d \nph = 0$, we
need
\begin{eqnarray*}
d \nPh & = & 4 f^3 d f \wedge \oPh + f^4 d \oPh = 0 \\ d \oPh & =
& -4 d \log (f) \wedge \oPh
\end{eqnarray*}
which says that $d \oPh \in \wfie$ by
Proposition~\ref{SPisomorphisms}. Hence $\pi_{48}(d \oPh) = 0$ so
$\oPh$ must be already of class $W_8$. Then to make $d \nPh = 0$,
we need to eliminate the $W_8$ component, which requires $28 d
\log(f) = \otheta$ by Theorem~\ref{SPconformaldeformresults}.
\end{proof}

\begin{rmk} \label{SPconformaldeformrmk}
Note that if we start with a \SPs\ $\oPh$ that is already holonomy
\SP, then Theorem~\ref{SPconformaldeformresults} shows that a
conformal scaling by a non-constant $f$ will always generate a
non-zero $W_8$ component.
\end{rmk}

\subsection{Deforming $\Ph$ by an element of $\wfos$}
\label{w47deformsec}

We continue our analogy with the \G\ case and now try to deform the
\SP\ $4$-form $\Ph$ by an element of $\wfos$.
Using~\eqref{wfosdesc2}, one can check that if we start with two
vector fields $v$ and $w$, we can construct a special kind of
element $\sigma_7 \in \wfos$ by $\sigma_7 = \vs \wedge (w \hk \oPh)
- \ws \wedge (v \hk \oPh)$. We will consider this type since at
least locally every element in $\wfos$ is a linear combination
of elements of this type. Now let $\nPh = \oPh + t \left( \vs
\wedge (w \hk \oPh) - \ws \wedge (v \hk \oPh) \right)$, for $t \in
\mathbb R$. Using the notation of Theorem~\ref{SPmetrictheorem}, we
first prove the following.

\begin{prop} \label{w47deformnewdenomprop}
Let $\sigma_7 = \left( \vs \wedge (w \hk \oPh) - \ws \wedge (v \hk
\oPh) \right)$. Under the transformation $\nPh = \oPh + 
\sigma_7$, we have
\begin{equation} \label{w47deformnewvoleq}
\nPh^2 = \opvw \oPh^2
\end{equation}
\end{prop}
\begin{proof}
We compute
\begin{eqnarray*}
\nPh^2 & = & {\left( \oPh + \vs \wedge (w \hk \oPh) - \ws \wedge (v
\hk \oPh) \right)}^2 \\ & = & \oPh^2 + 2 \vs \wedge (w \hk \oPh)
\wedge \oPh - 2 \ws \wedge (v \hk \oPh) \wedge \oPh - 2 \vs \wedge
(w \hk \oPh) \wedge \ws \wedge (v \hk \oPh) \\ & = & \oPh^2 + 14
\vs \wedge \ost \ws - 14 \ws \wedge \ost \vs - 2 \vs \wedge \ws
\wedge (v \hk \oPh) \wedge (w \hk \oPh) \\ & = & \oPh^2 + 8 {|v
\wedge w|}^2_{\text o} \ovol
\end{eqnarray*}
where we have used both~\eqref{hk7eq} and~\eqref{SP2formdot2eq}.
Now since $\oph^2 = 14 \ovol$ we have
\begin{equation*}
\nPh^2 = 14 \ovol + 14 \left( \frac{4}{7} \right) \ovw^2 \ovol =
\opvw \oPh^2
\end{equation*}
which completes the proof.
\end{proof}
\begin{cor} \label{w47deformnewdenomcor}
Let $h$ be a vector field. Under the transformation $\nPh = \oPh + 
\sigma_7$, the expression
$A_{\text{o}}(h) = \left( (h \hk \oPh) \wedge \oPh \right) (e_1,
e_2, \ldots, e_7)$ changes by
\begin{eqnarray} \label{w47deformnewdenomeq}
\tilde A (h) & = & \left( (h \hk \nPh) \wedge \nPh \right) (e_1,
e_2, \ldots, e_7) \\ \nonumber & = & \opvw A_{\text{o}}(h)
\end{eqnarray}
\end{cor}
\begin{proof}
This follows from Proposition~\ref{w47deformnewdenomprop}
by taking the interior product of both sides with $h$.
\end{proof}

We continue the computation of the expressions needed to determine
if $\nPh$ is indeed a \SPs\ with the following lemma.
\begin{lemma} \label{w47deformmetriclemma}
With $\nPh = \oPh + t \sigma$, in the expression
\begin{equation*}
\left( e_i \hk u \hk \nPh \right) \wedge \left( e_i \hk u \hk
\nPh \right) \wedge \nPh
\end{equation*}
which is a cubic polynomial in $t$, the linear and cubic terms
both vanish, and the coefficient of the quadratic term is 
\begin{eqnarray*}
& & 6 \left( - {\oPh(v,w,h,e_i)}^2 + {|v \wedge w \wedge h \wedge
e_i |}^2 -2 \langle h \wedge e_i , v \wedge w \rangle
\oPh(v,w,h,e_i) \right) \ovol \\ & & {}+ 6 \left( \langle w \wedge
e_i , w \wedge v \rangle \langle h \wedge e_i, h \wedge v \rangle +
\langle e_i \wedge h, e_i \wedge v \rangle \langle w \wedge h, w
\wedge v \rangle \right) \ovol \\ & & {}+ 6 \left( \langle h \wedge
e_i , h \wedge w \rangle \langle v \wedge e_i, v \wedge w \rangle +
\langle e_i \wedge h, e_i \wedge w \rangle \langle v \wedge h, v
\wedge w \rangle \right) \ovol \\ & & {}-12 {\langle h \wedge
e_i , v \wedge w \rangle}^2 \ovol
\end{eqnarray*}
\end{lemma}
\begin{proof}
We begin with the linear term. The coefficient of $t$ is
\begin{equation*}
2 ( e_i \hk u \hk \sigma ) \wedge ( e_i \hk u \hk \oPh ) \wedge
\oPh + ( e_i \hk u \hk \oPh) \wedge ( e_i \hk u \hk \oPh ) \wedge
\sigma
\end{equation*}
Two applications of~\eqref{usefuleq2} and collecting terms shows
that this coefficient is
\begin{equation*}
3 ( e_i \hk u \hk \oPh ) \wedge ( e_i \hk u \hk \oPh ) \wedge
\sigma
\end{equation*}
which vanishes by Corollary~\ref{w224zerocor}.

Next we consider the cubic term:
\begin{equation*}
\left( e_i \hk u \hk \sigma \right) \wedge \left( e_i \hk u \hk
\sigma \right) \wedge \sigma
\end{equation*}
This is a lengthy computation. We will sketch the steps involved 
and the interested reader can reproduce the details if desired. 
Let $e_i = h$ for notational simplicity. First, one can compute
that
\begin{eqnarray*}
h \hk u \hk \sigma & = & \langle u,v \rangle \left( h \hk w \hk
\oPh \right) - \langle h,v \rangle \left( u \hk w \hk \oPh
\right) - \langle u,w \rangle \left( h \hk v \hk \oPh \right) \\ & & {}+ \langle h,w \rangle \left( u \hk v \hk \oPh
\right) + \vs \wedge \left( h \hk u \hk w \hk \oPh \right) - \ws
\wedge \left( h \hk u \hk v \hk \oPh \right)
\end{eqnarray*}
and also
\begin{eqnarray*}
v \hk h \hk u \hk \sigma & = & \langle u,v \rangle \left( v \hk h
\hk w \hk \oPh \right) - \langle h,v \rangle \left( v \hk u \hk w
\hk \oPh \right) + {|v|}^2 \left( h \hk u \hk w \hk \oPh
\right) \\ & & {}- \langle v,w \rangle \left( h \hk u \hk v \hk
\oPh \right) - \oPh(w,u,h,v) \vs
\end{eqnarray*}
Now starting from the $9$-form $\left( h \hk u \hk \sigma \right)
\wedge \left( h \hk u \hk \sigma \right) \wedge \ws \wedge \oPh =
0$, taking the interior product with $v$, and rearranging, we
obtain
\begin{eqnarray*}
\left( h \hk u \hk \sigma \right) \wedge \left( h \hk u \hk \sigma
\right) \wedge \sigma & = & 2 \left( v \hk h \hk u \hk \sigma
\right) \wedge \ws \wedge \left( h \hk u \hk \sigma \right) \wedge
\oPh \\ & & {}- 2 \left( w \hk h \hk u \hk \sigma \right) \wedge
\vs \wedge \left( h \hk u \hk \sigma \right) \wedge \oPh 
\end{eqnarray*}
Since $\sigma \mapsto - \sigma$ if we interchange $v$ and $w$, we
see that the cubic term will vanish provided we can check that the
first term on the right hand side above is symmetric in $v$ and
$w$. Now one substitutes the expressions above for $h \hk u \hk
\sigma$ and $v \hk h \hk u \hk \sigma$ and
uses~\eqref{SP2formdot1eq} many times, as well as $\as \wedge \bs
\wedge \cs \wedge \ds \wedge \oPh = \oPh(a,b,c,d) \ovol$ to verify
that this is indeed the case.

We move on to the coefficient of $t^2$, which is
\begin{equation*}
2 ( e_i \hk u \hk \oPh ) \wedge ( e_i \hk u \hk \sigma ) \wedge
\sigma + ( e_i \hk u \hk \sigma) \wedge ( e_i \hk u \hk \sigma )
\wedge \oPh
\end{equation*}
As with the linear term, this can be rewritten as
\begin{equation*}
3 ( e_i \hk u \hk \sigma) \wedge ( e_i \hk u \hk \sigma )
\wedge \oPh
\end{equation*}
As in the case of the cubic term, the interested reader can now
apply Proposition~\ref{SP2formsprop} and
Lemma~\ref{cross8iteratelemma} many times to obtain the stated
expression. We will not provide the details since we will now see
that this deformation will not yield a \SPs\ anyway.
\end{proof}

If we now polarize the expression $( h \hk e_i \hk \nPh) \wedge (h
\hk e_i \hk \nPh) \wedge \nPh$, take the interior product with
$h$, and apply this to a basis extension $e_1, e_2, 
\ldots, e_7$, as required by Theorem~\ref{SPmetrictheorem}, one
can check that
\begin{eqnarray*}
\frac{1}{6} \widetilde{B}_{ij} & = & {|h|^2_{\text{o}}}\left(1 +
\ovw^2_{\text{o}} \right) g_{ij} - {\langle w, h \rangle}^2 v_i
v_j - {\langle v, h \rangle}^2 w_i w_j \\ & & {}+ \langle v, h
\rangle \langle w, h \rangle \left( v_i w_j + w_i v_j \right) -
\left(1 + \ovw^2_{\text{o}} \right) h_i h_j - X_i X_j \\ & &
{}-\langle w,h \rangle \left( v_i X_j + X_i v_j \right) +\langle
v,h \rangle \left( w_i X_j + X_i w_j \right) 
\end{eqnarray*}
where $X$ is the vector field $X(v,w,h)$. From this expression the
determinant of $\widetilde{B}_{ij}$ can be computed as
\begin{equation*}
\det \left( \widetilde{B}_{ij} \right) = 6^7 {|h|^{12}_{\text{o}}}
{\left(1 + \ovw^2_{\text{o}} \right)}^6
\end{equation*}
Now if this was indeed a \SPs\ then Theorem~\ref{SPmetrictheorem}
would imply that
\begin{equation*}
\lt{|h|}^{4} = \frac{{\left(1 + \ovw^2_{\text{o}}
\right)}^2}{\opvw^3 } {|h|^{4}_{\text{o}}}
\end{equation*}
This would mean the metric changes conformally, but the conformal
factor is {\em not} compatible with what would be the new
volume form $\nvol = \frac{1}{14} \nPh^2$ from
Proposition~\ref{w47deformnewdenomprop}. Hence this is {\em never}
a \SPs. Note that the construction very closely parallels the \G\
case. Even though the deformation does not yield a \SPs, it is
nevertheless true that $\det(\widetilde{B}_{ij})$ turns out to be a
positive definite quadratic form.

Recall now one major difference between the \G\ and \SP\ cases:
the space $\wthpos$ of \Gs s at a point is an {\em open subset} of
the space $\wedge^3$ of $3$-forms at that point. In constrast, the
space $\wfopos$ of \SPs s at a point is a $43$-dimensional
submanifold of the $70$-dimensional vector space $\wedge^4$ of
$4$-forms at that point, and this submanifold is not linearly
embedded. So we should not expect that moving linearly in
$\wedge^4$ would keep us on this submanifold, in
general~\cite{Br3}. In effect, the \SP\ case is {\em more
non-linear} than the \G\ case. There may still exist, however, some
non-linear way of deforming $\oPh$ by an element of $\wfos$ to
obtain a \SPs. In Section~\ref{conclusionsec} we present an
argument as to why this may be.

\subsection{Infinitesmal deformations in the $\wfos$ direction}
\label{infw47deformsec}

Even though the non-infinitesmal $\wfos$ deformation did not
produce a \SPs, we will see that in analogy with the \G\ case, we
can get a family of \SPs s all corresponding to the same metric by
taking {\em infinitesmal deformations} in the
$\wfos$ direction. Consider a one-parameter family $\Ph_t$ of \SPs
s, satisfying
\begin{equation} \label{infw47deformeq}
\frac{\partial}{\partial t} \Ph_t = w \hk \st_t \left( v \hk \Ph_t
\right) - v \hk \st_t \left( w \hk \Ph_t \right)
\end{equation}
for a pair of vector fields $v$ and $w$. That is, at each time $t$,
we move in the direction of a $4$-form in $\wedge^4_{7_t}$, since
the decomposition of $\wedge^4$ depends on $\Ph_t$ and
hence is changing in time. Since the Hodge star $\st_t$ is also
changing in time, this is again {\em a priori} a nonlinear
equation. However, just like in the \G\ case, it is actually
linear:
\begin{prop} \label{infw47deformlinprop}
Under the flow described by equation~\eqref{infw47deformeq}, the
metric $g$ {\em does not change}. Hence the volume form and
Hodge star are also constant.
\end{prop}
\begin{proof}
From Theorem~\ref{SPmetrictheorem},
Corollary~\ref{w47deformnewdenomcor}, and
Lemma~\ref{w47deformmetriclemma} we see that if we expand the
expression for $|h|^4$ for some vector field $h$, as a power series
in $t$, there is no linear term and hence to first order the
metric does not change.
\end{proof}
Therefore we can replace $\st_t$ by $\st_0 = \st$ and
equation~\eqref{infw47deformeq} is actually {\em linear}. Moreover,
the flow determined by this linear equation gives a one-parameter
family of \SPs s each yielding the {\em same} metric $g$. Our
equation is now
\begin{equation*}
\frac{\partial}{\partial t} \Ph_t = w \hk \st \left( v \hk \Ph_t
\right) - v \hk \st \left( w \hk \Ph_t \right) = B \Ph_t
\end{equation*}
where $B$ is the linear operator $\alpha \mapsto B \alpha = w \hk
\st ( v \hk \alpha) - v \st ( w \hk \alpha )$ on $\wedge^4$.
\begin{prop} \label{SPskewopprop}
The operator $B$ is {\em skew-symmetric}. Furthermore, the
eigenvalues $\lambda$ of $B$ are $\lambda = 0, \pm i \ovw$.
\end{prop}
\begin{proof}
Let $e^1, e^2, \ldots, e^{70}$ be a basis of $\wedge^4$. Then
\begin{eqnarray*}
B_{ij} \vol & = &  \left \langle e^i, B e^j \right \rangle \vol \\
& = & e^i \wedge \st \left( v \hk (\ws \wedge e^j) - w \hk (\vs
\wedge e^j) \right) \\ & = & \vs \wedge e^i \wedge \st (\ws
\wedge e^j) - \ws \wedge e^i \wedge \st ( \vs \wedge e^j) \\ & = &
\ws \wedge e^j \wedge \st ( \vs \wedge e^i) - \ws \wedge e^i \wedge
\st ( \vs \wedge e^j) \\ & = & - B_{ji} \vol
\end{eqnarray*}
and hence $B$ is diagonalizable over $\mathbb C$. In order to find
the eigenvalues of $B$, we first compute some powers of $B$, using
the fact that $B \alpha$ can also be written as $B \alpha = \vs
\wedge (w \hk \alpha) - \ws \wedge (v \hk \alpha)$:
\begin{eqnarray*}
B^2 \alpha & = & v \hk ( \ws \wedge B \alpha ) - w \hk ( \vs
\wedge B \alpha ) \\ & = & v \hk \left( \ws \wedge \vs \wedge (w
\hk \alpha ) \right) + w \hk \left( \vs \wedge \ws \wedge (v
\hk \alpha ) \right) \\ & = & {\langle v, w \rangle} \left( \vs
\wedge ( w \hk \alpha ) + \ws \wedge ( v \hk \alpha ) \right) -
{|v|}^2 \ws \wedge (w \hk \alpha) \\ & & {} - {|w|}^2 \vs \wedge
(v \hk \alpha) - 2 \vs \wedge \ws \wedge (v \hk w \hk \alpha )
\end{eqnarray*}
and thus
\begin{eqnarray*}
B^3 \alpha & = & v \hk ( \ws \wedge B^2 \alpha ) - w \hk ( \vs
\wedge B^2 \alpha ) \\ & = & v \hk \left( {\langle v, w \rangle}
\ws \wedge \vs \wedge (w \hk \alpha) - {|w|}^2 \ws \wedge \vs
\wedge (v \hk \alpha) \right) \\ & & {} - w \hk \left( {\langle v,
w \rangle} \vs \wedge \ws \wedge (v \hk \alpha) - {|v|}^2 \vs
\wedge \ws \wedge (w \hk \alpha) \right) \\ & = & {\langle v, w
\rangle}^2 B \alpha - {|v|}^2 {|w|}^2 B \alpha 
\end{eqnarray*}
after some simplification. Thus we have $\lambda^3 = \left(
{\langle v, w \rangle}^2 - {|v|}^2 {|w|}^2 \right) \lambda$.
Therefore the non-zero eigvenvalues are $\lambda = \pm i \ovw$.
This completes the proof.
\end{proof}
Now we proceed exactly as in the \G\ case. If we replace $A$ by
$B$ and the non-zero eigenvalues by $\pm i \ovw$, then all
the remaining calculations of Section~\ref{infvfdeformsec} carry
through. Therefore we have
\begin{equation*}
\Ph_t = -\frac{1}{\ovw^2} \cos(\ovw t) B^2 \Ph_0 +
\frac{1}{\ovw} \sin(\ovw t) B \Ph_0 + \Ph_0 + \frac{1}{\ovw^2} B^2
\Ph_0
\end{equation*}
which we summarize as the following theorem.
\begin{thm} \label{infw47deformsolutionthm}
The solution to the differential equation
\begin{equation*}
\frac{\partial}{\partial t} \Ph_t =  w \hk \st \left( v \hk \Ph_t
\right) - v \hk \st \left( w \hk \Ph_t \right)
\end{equation*}
is given by
\begin{equation} \label{infw47deformsolutioneq}
\Ph (t) = \Ph_0 + \frac{1 - \cos(\ovw t)}{\ovw^2} B^2 \Ph_0 +
\frac{\sin(\ovw t)}{\ovw} B \Ph_0
\end{equation}
where $B \alpha = v \hk (\ws \wedge \alpha) - w \hk (\vs \wedge
\alpha)$.  The solution exists for all time and is closed curve in
$\wedge^4$.
\end{thm}
\begin{proof}
This follows from the above discussion.
\end{proof}
\begin{rmk}
In~\cite{BS}, it is shown that the set of \SP s on $M$ which
correspond to the same metric as that of a fixed \SPs\ $\oPh$ is an
$\operatorname{O}(8)/ \SP$-bundle (which is rank $7$) over the
manifold $M$. The above theorem gives an explicit
formula~\eqref{infw47deformsolutioneq} for a path of \SPs s all
corresponding to the same metric $g$ starting from two vector
fields $v$ and $w$ on $M$.
\end{rmk}
\begin{rmk}
Again, even though the metric is unchanged under an infinitesmal
deformation in the $\wfos$ direction, the class of \SPs\ {\em
can} change.
\end{rmk}

We now apply this theorem to two examples, where we will again
reproduce known results. 
\begin{ex} \label{infw47deformCYex}
Let $N$ be a Calabi-Yau fourfold, with K\"ahler form $\omega$ and
holomorphic $(4,0)$ form $\Omega$. The complex coordinates will
be denoted by $z^j = x^j + i y^j$. Then $N$ has a natural \SPs\
$\Ph$ on it given by
\begin{equation} \label{cy4phieq}
\Ph = \real (\Omega) + \frac{\omega^2}{2}
\end{equation}
It is easy to check in local coordinates that $\omega \in \wtws$ in
the \SP\ decomposition. Since we are computing pointwise, if we
take two tangent vectoras $v$ and $w$ for which $\pi_7( \vs
\wedge \ws ) = \omega$, then one can compute that
\begin{eqnarray*}
w \hk \st (v \hk \Ph) - v \hk \st (w \hk \Ph) = & B \Ph & = -
\imag (\Omega) \\ & B^2 \Ph & = -\real (\Omega)
\end{eqnarray*}
Thus for the element of $\wfos$ which corresponds to
$\omega$, the flow in~\eqref{infw47deformsolutioneq} is given by
\begin{eqnarray*}
\Ph_t & = & \real (\Omega) + \frac{\omega^2}{2} - {(1 -
\cos(t))} \real (\Omega) - {\sin(t)} \imag (\Omega) \\ & = &
\cos(t) \real (\Omega) - \sin(t) \imag (\Omega) +
\frac{\omega^2}{2} \\ & = & \real( e^{i t} \Omega ) +
\frac{\omega^2}{2}
\end{eqnarray*}
which is the canonical \SP\ form on $N$ where now the
Calabi-Yau structure is given by $e^{i t} \Omega$ and $\omega$.
Thus we arrive at the phase freedom for Calabi-Yau fourfolds. 
\end{ex}
\begin{ex} \label{infw47deformGex}
Consider a $7$-manifold $M$ with a \Gs\ $\ph$. We can put
a \SPs\ $\Ph$ on the product $M \times S^1$ given by
\begin{equation*}
\Ph = d\theta \wedge \ph + \st_7 \ph
\end{equation*}
where $\st_7 \ph$ is the $4$-form dual to $\ph$ on $M$.
This induces the product metric on $M \times S^1$, with the flat
metric on $S^1$. Now let $v = \frac{\partial}{\partial \theta}$ be
a globally defined non-vanishing vector field on $S^1$ with $|v|=
1$. Choose another vector field $w$ on $M$. Then one computes
\begin{eqnarray*}
w \hk \st (v \hk \Ph) - v \hk \st (w \hk \Ph) = & B \Ph & = d
\theta \wedge ( w \hk \st_7 \ph ) + \st_7 ( w \hk \st_7 \ph) \\ &
B^2 \Ph & = d \theta \wedge \left( w \hk \st_7 ( w \hk \st_7 \ph)
\right) + \st_7 \left( w \hk \st_7 ( w \hk \st_7 \ph) \right)
\end{eqnarray*}
The flow in~\eqref{infw47deformsolutioneq} gives
\begin{equation*}
\Ph_t = d \theta \wedge \ph_t + \st_7 \ph_t
\end{equation*}
where $\ph_t$ is the flow given by~\eqref{infvfdeformsolutioneq}
for the vector field $w$. Thus in the product case $M \times S^1$
we recover the results of Section~\ref{infvfdeformsec}.
\end{ex}

\section{Conclusion} \label{conclusionsec}

In the construction of Calabi-Yau manifolds, we start from a
K\"ahler manifold and we reduce the holonomy from
$\operatorname{U}(n)$ to $\operatorname{SU}(n)$, which is a drop
of $1$ in dimension and hence it might be expected that it would
involve the solution of an equation for one function. In
going from $\operatorname{SO}(7)$ to \G\, we have a drop of $7$ in
dimension, so we might expect to need $7$ conditions, which could
involve an equation for a vector field (or equivalently an element
of $\wths$). Similarly the difference in dimension between
$\operatorname{SO}(8)$ and \SP\ is also $7$, which could be related
to an element $\wfos$.

Note that in the \G\ case, since elements of $\wths$ are
canonically identified with vector fields, they are intrinsic to
the manifold without reference to a \Gs. For \SPs s, we need the
$4$-form $\Ph$ to define $\wfos$ and this introduces more
non-linearity. To maintain the analogy with the \G\ case, there
should be some (non-linear) way of transforming a \SPs\ $\Ph$ using
an element of $\wfos$ so that we get a new \SPs\ whose new metric
is related to the old one by
\begin{equation*}
\lt{\langle u_1, u_2 \rangle} = f \left( {\langle u_1, u_2
\rangle}_{\! \text{o}} + {\langle X(v,w,u_1), X(v,w,u_2)
\rangle}_{\! \text{o}} \right)
\end{equation*}
where $v$ and $w$ are vector fields which determine the
corresponding element of $\wfos$ and $f$ is some positive function
of $\ovw^2$.

\appendix

\section{Some Linear Algebra}
\label{identitiessec}

Here we collect together various identities involving the exterior
and interior products and the Hodge star operator. Also some
identities involving determinants are proved. We establish the
results in the general case of a Riemannian manifold
$M$ of dimension $n$, although only the cases $n=7,8$ are used in
the text. Let $ \left \langle \ ,
\  \right \rangle $ denote the metric on $M$, as well as the
induced metric on forms. In all that follows, $\alpha$ and
$\gamma$ are $k$-forms, $\beta$ is a
$(k-1)$-form, $w$ is a vector field, and $\ws$ is the $1$-form dual
to $w$ in the given metric. That is,
\begin{equation*}
|w|^2 =  \langle w,w \rangle  = \ws(w) = 
\langle \ws,\ws \rangle 
\end{equation*}
Now $\st$ takes $k$-forms to $(n-k)$-forms, and is defined by
\begin{equation*}
 \left \langle \alpha, \gamma \right \rangle  \vol = \alpha \wedge
\st \gamma = \gamma \wedge \st \alpha
\end{equation*}
We also have
\begin{equation}
\label{starsquaredeq} \st^2 = (-1)^{k(n-k)}
\end{equation}
on $k$-forms.

\begin{lemma} \label{intextidentities}
We have the following four identities:
\begin{eqnarray} \label{ident1} \st(w \hk \alpha) & = & (-1)^{k+1}
(\ws \wedge \st \alpha) \\ \label{ident2} (w \hk \alpha) & =
& (-1)^{nk + n} \st (\ws \wedge \st \alpha) \\ \label{ident3} \st
( w \hk \st \alpha ) & = & (-1)^{nk + n + 1}(\ws \wedge \alpha
) \\ \label{ident4} (w \hk \st \alpha) & = & (-1)^k \st (\ws
\wedge \alpha) 
\end{eqnarray}
and when $\alpha = \vol$, the special case
\begin{equation} \label{wvoleq}
w \hk \vol = \st \ws
\end{equation}
\end{lemma}
\begin{proof}
We compute:
\begin{eqnarray*}
 \left \langle \beta, w \hk \alpha  \right \rangle  \vol & = &
\beta \wedge \st (w \hk \alpha) \\ & = & (w \hk
\alpha)(\beta^{\sharp}) \vol \\ & = & \alpha(w \wedge
\beta^{\sharp}) \vol \\ & = &  \left \langle \alpha, \ws \wedge
\beta  \right \rangle  \vol \\ & = & (\ws \wedge \beta) \wedge \st
\alpha \\ & = & (-1)^{k-1}\beta
\wedge (\ws \wedge \st \alpha)
\end{eqnarray*}
Since $\beta$ is arbitrary, \eqref{ident1} follows.
Substituting $\st \alpha$ for $\alpha$, and
using~\eqref{starsquaredeq}, we obtain~\eqref{ident3}. The
other two are obtained by taking $\st$ of both sides of the first
two identities. 
\end{proof}
Note that from the above proof we have also the useful relation
\begin{equation} \label{usefuleq} \left( X \hk \alpha \right)
\wedge \st \beta = \alpha \wedge \st \left( X^{\flat} \wedge \beta
\right)
\end{equation}
for any $k$-form $\alpha$, $(k-1)$-form $\beta$, and vector
field $X$. If $\alpha$ is a $k$-form, and $\beta$ is an
$(n+1-k)$-form, then $\alpha \wedge \beta = 0$ and taking the
interior product with a vector field $X$ gives
\begin{equation} \label{usefuleq2}
( X \hk \alpha ) \wedge \beta = {(-1)}^{k+1} \alpha \wedge ( X \hk
\beta )
\end{equation}
which is also used very often.

The next lemma gives further relations between a $k$-form $\alpha$
and a vector field $w$.
\begin{lemma} \label{wwidentities}
With the same notation as above, we have the following three
identities:
\begin{eqnarray} \label{ww1eq} |w|^2 \alpha & = & \ws \wedge (w
\hk \alpha) + w \hk (\ws \wedge \alpha) \\ \label{ww2eq} |w|^2
\alpha & = & (-1)^{nk + 1} \st \left( w \hk \left( \st (w \hk
\alpha) \right) \right) + (-1)^{nk + n + 1}\left( w \hk \st \left(
w \hk \st \alpha \right) \right) \\ \label{ww3eq} |w|^2 |\alpha|^2
& = & |w \hk \alpha|^2 + |w \hk \st \alpha|^2
\end{eqnarray}
\end{lemma}
\begin{proof}
Using the fact that the interior product is an anti-derivation, we
have
\begin{equation*}
w \hk (\ws \wedge \alpha) = (w \hk \ws) \wedge \alpha - \ws \wedge
(w \hk \alpha)
\end{equation*}
which is equation~\eqref{ww1eq} since $w \hk \ws = |w|^2$.
Now~\eqref{ww2eq} follows from this one using the identities
of Lemma~\ref{intextidentities}. The last identity can also be
obtained this way, but it is faster to proceed as follows:
\begin{eqnarray*}
0 & = & \ws \wedge \alpha \wedge \st \alpha \\ 0 & = & w \hk (\ws
\wedge \alpha \wedge \st \alpha) \\  & = & |w|^2 \alpha \wedge \st
\alpha - \ws \wedge (w \hk \alpha) \wedge \st \alpha + (-1)^{k+1}
\ws \wedge \alpha \wedge (w \hk \st \alpha) \\ & = & |w|^2
|\alpha|^2 \vol - (w \hk \alpha) \wedge \st (w \hk \alpha) - (w
\hk \st \alpha) \wedge \st (w \hk \st \alpha)
\end{eqnarray*}
using Lemma~\ref{intextidentities}, which can be rearranged to
give~\eqref{ww3eq}.
\end{proof}

The following lemma about determinants is used many times in the
computation of the metrics and volume forms arising from \G\ and
\SPs s.
\begin{lemma} \label{detlemma}
Let $g_{ij}$ be an $n \times n$ matrix, $v_i$ and $w_j$ be
two $n \times 1$ vectors, and $C, K$ constants. Consider the matrix
\begin{equation*}
B_{ij} = C g_{ij} + K v_i w_j
\end{equation*}
Its determinant is given by
\begin{equation} \label{detlemmaeq}
\det (B) = C^n \det (g) + \sum_{k,l=1}^n (-1)^{k+l} v_k
w_l C^{n-1} K G_{kl}
\end{equation}
where $G_{kl}$ is the ${(k,l)}^{\text{th}}$ {\em minor} of the
matrix $g_{ij}$. That is, it is the determinant of $g_{ij}$ with
the $k^{\text{th}}$ row and $l^{\text{th}}$ column removed
\end{lemma}
\begin{proof}
The determinant of $B_{ij}$ is
\begin{equation*}
\det {\begin{pmatrix} C g_{11} + K v_1 w_1 & \ldots & C g_{1n} + K
v_1 w_n \\ \vdots & \ddots & \vdots \\ C g_{n1} + K v_n w_1 &
\ldots & C g_{nn} + K v_n w_n \end{pmatrix}}
\end{equation*}
Since the determinant is a linear function of the {\em
columns} of a matrix, we can write the above determinant as a
sum of determinants where each column is of one of these two forms:
\begin{equation*}
\begin{pmatrix} C g_{1k} \\ C g_{2k} \\ \vdots \\ C g_{nk}
\end{pmatrix} \qquad \text{ or } \qquad \begin{pmatrix} K v_1 w_k
\\ K v_2 w_k \\ \vdots \\ K v_n w_k \end{pmatrix}
\end{equation*}
Each of the determinants which has at least two columns of
$K v_l w_k$'s will have at least two proportional columns and
hence will vanish. So the only non-zero contributions come from the
determinant with all $C g_{kl}$'s and the $n$ determinants with
only {\em one} column of $K v_l w_k$'s. Therefore we are left
with:
\begin{equation} \label{detlemmatempeq}
\det{ \left( C g \right)} + \sum_{k,l = 1}^n {(-1)^{k+l}
\left( K v_k w_l \right) C^{n-1} (G_{kl})}
\end{equation}
In equation~\eqref{detlemmatempeq} the first term is the
determinant with all columns of $C g_{kl}$'s, the sum over $l$ is
the sum over the $n$ different determinants with a single column of
$K v_k w_l$'s, in the $l^{\text{th}}$ column, and in the sum
over $k$ we expand each of those determinants along the
$l^{\text{th}}$ column. This completes the proof.
\end{proof}

We can obtain a special case of Lemma~\ref{detlemma} when the
matrix $g_{ij}$ is a metric. It is used several times in the text,
most notably in the derivation of the metric from the $4$-form
$\Ph$ in the \SP\ case in Theorem~\ref{SPmetrictheorem}.
\begin{lemma} \label{specialdetlemma}
Let $g_{ij} = \langle e_i, e_j \rangle$ be a Riemannian metric in
local coordinates, and $\vs = v_i e^i$ and $\ws = w_j e^j$ be two
one forms dual to the vector fields $v$ and $w$. Then if we define
\begin{equation*}
B_{ij} = C g_{ij} + K v_i w_j
\end{equation*}
we have
\begin{equation*}
\det (B) = C^n \det (g) + C^{n-1} K \langle v, w \rangle \det (g)
\end{equation*}
\end{lemma}
\begin{proof}
In local coordinates, the volume form is $\vol =
\sqrt{\det(g)} e^1 \wedge e^2 \wedge \ldots \wedge e^n$. Then
its dual $n$-vector is
\begin{equation*}
\vol^{\sharp} = \frac{1}{\sqrt{\det (g)}} e_1 \wedge e_2
\wedge \ldots \wedge e_n
\end{equation*}
The induced inner product on the two $(n-1)$-vectors $\vs
\hk \vol^{\sharp}$ and $\ws \hk \vol^{\sharp}$ is
\begin{eqnarray*}
& & \langle \vs \hk \vol^{\sharp}, \ws \hk \vol^{\sharp} \rangle
\\ & = & \frac{1}{\det (g)} \sum_{k,l = 1}^n {(-1)}^{k + l} v_k w_l
\langle e_1 \wedge \ldots \wedge \widehat{e_k} \wedge \ldots \wedge
e_n, e_1 \wedge \ldots \wedge \widehat{e_l} \wedge \ldots \wedge
e_n \rangle \\ & = & \frac{1}{\det (g)} \sum_{k,l = 1}^n {(-1)}^{k
+ l} v_k w_l G_{kl}
\end{eqnarray*}
Now comparing with~\eqref{detlemmatempeq}, we have
\begin{equation*}
\det(B) = C^n \det(g) + C^{n-1} K \langle \vs \hk \vol^{\sharp},
\ws \hk \vol^{\sharp} \rangle \det(g)
\end{equation*}
The result now follows since $\vs \hk \vol^{\sharp} = \st v$ and
the Hodge star is an isomorphism. (This is the star operator on
the space of $k$-vectors.)
\end{proof}

\end{document}